%% file: spheres.tex
\documentclass[12pt]{amsart}

\usepackage[colorlinks=true,urlcolor=blue,bookmarks=true,bookmarksopen=true,citecolor=blue]{hyperref}

\usepackage{graphicx}
\usepackage[usenames,dvipsnames]{color}
\usepackage[all]{xy}
\usepackage{amsmath}
\usepackage{amsfonts,amssymb}
\usepackage{amsthm}
\usepackage{bbm}
\usepackage{mathrsfs}
\usepackage{marginnote}
\usepackage{pdfsync}
\usepackage{xfrac}
\usepackage{varioref}
\usepackage{tikz-cd}
\usepackage{amscd}
\usepackage{url}
\usepackage{bbm}
\usepackage{upgreek}
\usepackage[all]{xy}
\usepackage{tikz}
\usetikzlibrary{decorations.pathreplacing}
\usepackage{epsf,overpic,amssymb,amscd,times,stmaryrd,euscript}
\usepackage{a4wide}
\usepackage{yfonts}
\usepackage[T1]{fontenc}

\usepackage{parskip}

\newcommand{\id}{id}

\makeatletter
\DeclareFontFamily{OMX}{MnSymbolE}{}
\DeclareSymbolFont{MnLargeSymbols}{OMX}{MnSymbolE}{m}{n}
\SetSymbolFont{MnLargeSymbols}{bold}{OMX}{MnSymbolE}{b}{n}
\DeclareFontShape{OMX}{MnSymbolE}{m}{n}{
    <-6>  MnSymbolE5
   <6-7>  MnSymbolE6
   <7-8>  MnSymbolE7
   <8-9>  MnSymbolE8
   <9-10> MnSymbolE9
  <10-12> MnSymbolE10
  <12->   MnSymbolE12
}{}
\DeclareFontShape{OMX}{MnSymbolE}{b}{n}{
    <-6>  MnSymbolE-Bold5
   <6-7>  MnSymbolE-Bold6
   <7-8>  MnSymbolE-Bold7
   <8-9>  MnSymbolE-Bold8
   <9-10> MnSymbolE-Bold9
  <10-12> MnSymbolE-Bold10
  <12->   MnSymbolE-Bold12
}{}

\let\llangle\@undefined
\let\rrangle\@undefined
\DeclareMathDelimiter{\llangle}{\mathopen}%
                     {MnLargeSymbols}{'164}{MnLargeSymbols}{'164}
\DeclareMathDelimiter{\rrangle}{\mathclose}%
                     {MnLargeSymbols}{'171}{MnLargeSymbols}{'171}
\makeatother

\newcommand{\Z}{\mathbbm{Z}}                     
\newcommand{\R}{\mathbbm{R}}                     
\newcommand{\C}{\mathbbm{C}}                     

\newcommand{\dist}{\mathrm{dist\,}}             
   
\newcommand{\euc}{\mathrm{euc}}                 
\newcommand{\round}{\mathrm{round}}                 
\newcommand{\crit}{\mathrm{Crit\,}}               

\newcommand{\Vol}{\operatorname{Vol}} 
\newcommand{\Cyl}{\operatorname{Cyl}} 
\newcommand{\im}{\mathrm{Im}}

\newcommand{\HW}{\mathrm{HW}}
\newcommand{\CW}{\mathrm{CW}}

\newcommand{\dlim}{\underrightarrow{\lim}\,} 

\newcommand{\Ho}{\mathrm{H}}
\newcommand{\symp}{\mathrm{symp}}
\newcommand{\alg}{\mathrm{alg}}
\newcommand{\can}{\mathrm{can}}
\newcommand{\reg}{\mathrm{reg}}
\newcommand{\regium}{\mathrm{regium}}
\newcommand{\topo}{\mathrm{top}}
\newcommand{\lef}{\mathrm{left}}
\newcommand{\rig}{\mathrm{right}}
\newcommand{\midde}{\mathrm{mid}}

\newtheorem{thm}{Theorem}[section]               
\newtheorem*{thm*}{Theorem}               
\newtheorem{cor}[thm]{Corollary}        
\newtheorem*{cor*}{Corollary}        
\newtheorem{lem}[thm]{Lemma}            
\newtheorem{prop}[thm]{Proposition}     
\newtheorem{conj}[thm]{Conjecture}      
\theoremstyle{definition}
\newtheorem{defn}[thm]{Definition}      
\newtheorem{rem}[thm]{Remark}           
 \newtheorem*{acknowledgement*}{\protect\acknowledgementname}

 \providecommand{\acknowledgementname}{Acknowledgement}

\date{} 

\author{Marcelo R.R. Alves}
\address{Institut de math\'ematiques\\
Universit\'e de Neuch\^atel}
\email{\texttt{marcelo.ribeiro@unine.ch}}

\author{Matthias Meiwes}
\address{Department of Mathematics\\
Universit\"at Heidelberg}
\email{\texttt{mmeiwes@mathi.uni-heidelberg.de}}

\thanks{Marcelo R.R. Alves supported by the Swiss National Foundation. Matthias Meiwes supported by German-Israeli Foundation (GIF).}

\title{Dynamically exotic contact spheres in dimensions~$\geq 7$}

\date{}



\begin{document}

\begin{abstract}
We exhibit the first examples of contact structures on $S^{2n-1}$ with $n\geq 4$ and on $S^3\times S^2$, all equipped with their standard smooth structures, for which every Reeb flow has positive topological entropy. 
As a new technical tool for the study of the volume growth of Reeb flows we introduce the notion of algebraic growth of wrapped Floer homology. Its power stems from its stability under several geometric operations on Liouville domains.
\end{abstract}
\maketitle
\thispagestyle{empty}
\section{Introduction}
\color{black}
 On a contact manifold there exists a natural class of flows, the so-called Reeb flows. Although the dynamics of distinct Reeb flows on the same contact manifold can be very different, there are dynamical properties which are common to all Reeb flows on a given contact manifold. For instance,  the combined works of Hofer \cite{Hofer} and Taubes \cite{Taubes} imply that on a closed contact 3-manifold all Reeb flows have at least one periodic orbit. 
 In this paper we construct a large class of contact manifolds on which all Reeb flows have chaotic dynamics. Surprisingly, some of the contact manifolds we construct have a very simple topology, which contrasts with the complicated dynamics of their Reeb flows.


A contact structure is said to have positive  entropy if all Reeb flows associated to this contact structure have positive topological entropy.
We show that there exist contact structures with positive entropy on spheres of dimension $\geq7$ and on $S^3\times S^2$.  As a consequence we prove that every manifold of dimension $\geq 7$ that admits an exactly fillable contact structure also admits a (possibly different) contact structure with positive  entropy. \color{black} 
Our approach to prove these results is based on wrapped Floer homology and uses in an essential way its product structure. This product structure enables us to define the notion of algebraic growth of wrapped Floer homology, and we relate this growth to the volume growth of Reeb flows.  Even though the richer algebraic structures in Floer homology were studied extensively, so far they lead to only very few applications in dynamics: the ones we are aware of are Viterbo's result \cite{Viterbo1999} on the existence of one closed Reeb orbit on hypersurfaces of restricted contact type in Liouville  domains with vanishing symplectic homology, and Ritter's result \cite{Ritter2013} on the existence of Reeb chords for exactly fillable Legendrian submanifolds on Liouville domains with vanishing symplectic homology.

\subsection{Basic notions}

An important measure of the complexity of a dynamical system on a manifold $M$ is the topological entropy $h_{top}$ which quantifies in a single number the exponential complexity of the system. We refer the reader to \cite{Hasselblat-Katok} for the definition and basic properties of $h_{\topo}$. 
By deep results of Yomdin and Newhouse, $h_{\topo}(\phi)$ for a $C^{\infty}$-flow $\phi = (\phi^t)_{t\in\R}$ equals the exponential growth rate of volume
\begin{equation*}
v(\phi) = \sup_{N\subset M} v(\phi,N), \text{ where }
\end{equation*}
\begin{equation} \label{Yomdin}
v(\phi, N) = \limsup_{t\to \infty} \frac{\log \Vol_g^n(\phi^t(N))}{t}.
\end{equation}
Here, $n = \dim N$, the  supremum is taken over all submanifolds $N \subset M$, and $\Vol^n_g$ is the $n$-dimensional volume with respect to some Riemannian metric $g$ on $M$. 

In this paper we study the topological entropy for Reeb flows of contact manifolds.
Recall that a \textit{(co-oriented) contact manifold} $(\Sigma, \xi)$ is 
a compact odd-dimensional manifold $\Sigma^{2n-1}$ equipped with a contact structure $\xi$, that is, a hyperplane distribution on $\Sigma$ which is given by $\xi = \ker \alpha$ for a $1$-form $\alpha$ with $\alpha \wedge (d\alpha)^{n-1} \neq 0$. Such an $\alpha$ is called a \textit{contact form} on $(\Sigma,\xi)$, and we can associate to it the \textit{Reeb vector field} $X_{\alpha}$ defined by $\iota_{X_{\alpha}}d\alpha = 0$, $\alpha(X_{\alpha})= 1$. Denote the flow of $X_\alpha$, the \textit{Reeb flow} of $\alpha$,  by $\phi_{\alpha}= (\phi_{\alpha}^t)_{t\in \R}$. 
An \textit{isotropic submanifold} of   $\Sigma^{2n-1}$ is one whose tangent space is contained in $\xi$; isotropic submanifolds of dimension $n-1$ are called \textit{Legendrian submanifolds}. 

\subsection{Main results}
The main result of this paper is the existence of contact structures with positive entropy on high dimensional manifolds.

 \begin{thm}\label{spheres}
 \
 
 \begin{itemize} 
 \item [A)]\label{itm:spheres}
 Let $S^{2n-1}$ be the $(2n-1)$ - dimensional sphere with its standard smooth structure.
For $n\geq 4$ there exists a contact structure on $S^{2n-1}$ with positive entropy.
\item [B)] \label{itm:S3S2}
There exists a contact structure on $S^3\times S^2$ with positive entropy.
\end{itemize}
\end{thm}

Recall that a contact manifold is said to be exactly fillable if it is the boundary of a Liouville domain. From Theorem \ref{spheres} and the methods developed in this paper we obtain the following more general result.
\begin{thm}\label{maincorollary}
\
\begin{itemize}
\item [$\clubsuit$] 
\ If  $V$ is a manifold of dimension $2n-1\geq 7$ that admits an exactly fillable contact structure, then  $V$ admits a contact structure with positive  entropy.
\item [$\diamondsuit$]
If $V$ is a $5$-manifold that admits an exactly fillable contact structure, then the connected sum  $V \# (S^3\times S^2)$ admits a contact structure with positive entropy.
\end{itemize}
\end{thm}

Note that the standard contact structure on spheres as well as the canonical contact structure on $S^{*}S^3 \cong S^3 \times S^2$ have a contact form with periodic Reeb flow. In particular these are not diffeomorphic to the contact structures in Theorem \ref{spheres}. Other exotic contact spheres have been constructed by several authors, see \cite{Eliashberg1991, Ustilovsky1999, GeigesDing2004, McLean2011}. The contact spheres constructed in this paper are, from our perspective, the ``most exotic'' ones. From the dynamical point of view they are the most remote from the standard contact spheres since they admit Legendrian submanifolds that have exponential volume growth under every Reeb flow. It would be interesting to relate our examples of  exotic contact spheres to others that were constructed so far. 

In order to explain further the relevance of these results we recall what is known about the topological entropy of Reeb flows.
Motivated by results on topological entropy for geodesic flows (see \cite{Paternain}),  combined with the geometric ideas of \cite{FrauenfelderSchlenk2006}, Macarini and Schlenk proved in \cite{MacariniSchlenk2011} that for various manifolds $Q$ the unit cotangent bundle $(S^{*}Q, \xi)$ equipped with the canonical contact structure $\xi$ has positive entropy\footnote{In a recent work \cite{Dahinden} Dahinden extended the results in \cite{MacariniSchlenk2011} proving that on the unit cotangent bundles $(S^{*}Q, \xi)$ studied in \cite{MacariniSchlenk2011} every positive contactomorphism has positive topological entropy. It would be interested to investigate if Dahinden's result is true for any contact manifold with positive entropy.}.

In previous works of the first author, different examples of contact 3-manifolds with positive entropy were discovered. In \cite{Alves-Cylindrical,Alves-Anosov,Alves-Legendrian} it was shown that contact 3-manifolds with positive entropy exist in abundance: there exist hyperbolic contact 3-manifolds with positive entropy (see also \cite{ACH}), non-fillable contact 3-manifolds with positive  entropy, and even 3-manifolds which admit infinitely many non-diffeomorphic contact structures with positive  entropy. This shows that the class of contact manifolds with positive entropy is much larger than the class of unit cotangent bundles over surfaces with positive entropy, which were studied in \cite{MacariniSchlenk2011}.
One common feature of all known examples of contact 3-manifolds with positive entropy is that the fundamental group of the underlying smooth 3-manifold has exponential growth. We expect this to be always the case:
\begin{conj} \label{conjecture3dim}
If a contact 3-manifold $(\Sigma,\xi)$ has positive entropy, then $\pi_1(\Sigma)$ grows exponentially.
\end{conj}

Already from the unit cotangent bundles of simply connected rationally hyperbolic manifolds, which were considered in \cite{MacariniSchlenk2011}, we know that Conjecture \ref{conjecture3dim} is false in higher dimensions.
However it is natural to ask if there are restrictions on the smooth topology of contact manifolds with positive  entropy in higher dimensions.
Theorem \ref{spheres} shows that in contrast to what happens in dimension three, the phenomenon in higher dimensions is quite flexible from the topological point of view.

\color{black}
\begin{rem} {Examples of contact manifolds of dimension $\geq 9$ which have positive entropy and are not unit cotangent bundles are also constructed using connected sums in an ongoing work of the first author and Macarini \cite{connected}, following an idea of Schlenk. However, these contact manifolds have very complicated smooth topology, in the sense that the underlying smooth manifolds are rationally hyperbolic. For this reason they are much less surprising than the ones obtained in the present paper.} \end{rem}
\color{black}

Let us now explain our approach to establishing these results.

\subsection{Symplectic and algebraic growth}

To establish our results we introduce the notion of algebraic growth of wrapped Floer homology. This notion is useful because, on one hand, it gives a lower bound for the growth rate of wrapped Floer homology defined using its action filtration and, on the other hand, it is stable under several geometric modifications of Liouville domains.

 The contact manifolds  we consider in this paper arise as boundaries  of Liouville domains. Recall that a Liouville domain $M=(Y,\omega, \lambda)$ is an exact symplectic manifold $(Y,\omega)$ with boundary $\Sigma = \partial Y$ and a primitive $\lambda$ of $\omega$ such that $\alpha_M := \lambda|_{\Sigma}$ is a contact form on $\Sigma$: we let $ \xi_{M}= \ker \alpha_M$ be the contact structure induced by $M$ on $\Sigma$. For two exact Lagrangians $L_0$ and $L_1$ in $M$ that are asymptotically conical, i.e. conical near $\partial Y$ with Legendrian boundaries $\Lambda_0$ and $\Lambda_1$ in $(\Sigma,\xi_M)$, we consider the wrapped Floer homology of $(M,L_0,L_1)$ with $\Z_2$-coefficients denoted by $\HW(M,L_0 \to L_1)$, whose underlying chain complex is, informally speaking, generated by Reeb chords from $\Lambda_0$ to $\Lambda_1$ and intersections of $L_0$ and  $L_1$. We write $\HW(M,L)$ for $\HW(M,L \to L)$, see Section \ref{subsec:def}.

Results on positive entropy can be obtained from the exponential  \textbf{ symplectic growth} of wrapped Floer homology, which is defined as follows. \color{black}
By  considering only critical points below an action value $a$, one obtains the filtered Floer homology $\HW^a(M,L_0 \to L_1)$. 
 The homologies $\HW^a(M,L_0 \to L_1)$ form a natural filtration of $\HW(M,L_0 \to L_1)$, and they come with natural maps $\iota_a: \HW^a(M,L_0 \to L_1) \rightarrow \HW(M,L_0 \to L_1)$ into the (unfiltered) Floer homology. The \textit{exponential symplectic growth rate} ${\Gamma}^{\symp}(M,L_0 \to  L_1)$ of $\HW(M,L_0 \to L_1)$ is given by
\begin{align}\label{Gamma_symp}
\Gamma^{\symp}(M,L_0 \to L_1) = \limsup_{a \to \infty} \frac{\log (\dim \im \ \iota_a)}{a};
\end{align}
see Section \ref{subsec:def} and Definition \ref{defi:Gamma_symp}.
Since the generators of  $\HW(M,L_0 \to L_1)$ correspond essentially to Reeb chords from $\Lambda_0$ to $\Lambda_1$, the symplectic growth gives a lower bound on the growth of Reeb chords with respect to their action. \color{black}
Assuming that $\Lambda_1$ is a sphere, we adapt the ideas of the first author in \cite{Alves-Legendrian} to get lower bounds for the volume growth $v(\phi_{\alpha},\Lambda_0)$ in terms of the exponential symplectic growth rate of $\HW(M,L_0,L_1)$ for every contact form $\alpha$ on $\xi_{M}$.

A \textit{topological operation} on a Liouville domain $M$ is a recipe for producing a new Liouville domain $N$ from $M$. 
To obtain  examples of contact manifolds with positive entropy  we  perform certain topological operations on Liouville domains. The operations we consider are: attaching symplectic handles on $M$ and, in the case $M$ is the unit disk bundle of a manifold, plumbing $M$ with the unit disk bundle of  another manifold. 
Although one can understand the change or invariance of the (unfiltered) wrapped Floer homology under these operations, it is often much harder or not even possible to understand the effect of these operations on the symplectic growth. 
For instance, by an adaptation of a theorem of Cieliebak \cite{Cieliebak2001} we show that $\HW({M}',L)$ is isomorphic to $\HW(M,L\cap M)$ if ${M}'$ is obtained by subcritical handle attachment on $M$ (Theorem \ref{Viterbo_iso}). By contrast it is much harder to control the filtered Floer homology under this operation, see \cite{McLean2011} for an approach in the case of symplectic homology. In the case of plumbings of two cotangent bundles the computational results of a relevant part of the unfiltered wrapped Floer homology obtained by \cite{AbouzaidSmith2012} do not carry over to the symplectic growth rates of the plumbing.

To overcome this difficulty we look at a notion of growth that is defined purely in terms the algebraic structure on wrapped Floer homology, the \textbf{algebraic growth}. 
Let us explain this briefly.  
Let $A$ be a (not necessarily unital) $K$-algebra with multiplication $\star$ and $S \subset A$ a finite set of elements of $A$. \color{black} 
Given $j\geq 0$, let $N_S(j) = \{ a \in A \, \mid \, a = s_1\star s_2 \star \cdots \star s_j; \, s_1, \dots, s_j \in S \}$; i.e.\ $N_S(j)$ is the set of elements of $A$ that can be written as a product of $j$ not necessarily distinct elements of $S$.
 We define $W_S(n) \subset A$ to be the smallest $K$-vector space that contains the union $\bigcup_{j=1}^{n} N_S(j)$. \color{black}  
 The \textit{exponential algebraic growth rate} of the pair $(A,S)$ is defined as \color{black}
\[
\Gamma^{\alg}_S(A) = \limsup_{n\to \infty}\frac{1}{n}\log\dim_K W_S(n) \in [0,\infty).
\]
\color{black}
 In case $A = K \langle G \rangle$ is the group algebra over a finitely generated group $(G = \langle S \rangle ,\star) $, 
it is elementary to see that  $\Gamma^{\alg}_S(A)$ coincides with the exponential algebraic growth of $G$ in the usual geometric group theoretical sense.
Now, induced by the triangle product in Floer homology, $\HW(M,L)$ is equipped with a ring structure $\star$ turning it into a $\Z_2$-algebra. \color{black} Given a finite set $S$ of $\HW(M,L)$ we define (cf.\ Definition \ref{defi:growthHW}) $\Gamma_{S}^{\alg}(M,L):= \Gamma_{S}^{\alg}(\HW(M,L)) $.
  We say that $\HW(M,L)$ has \textit{exponential algebraic growth} if there exists a finite subset $S$ of  $\HW(M,L)$ such that $\Gamma_{S}^{\alg}(M,L) > 0$.  \color{black}

Our main motivation for studying the exponential algebraic growth of $\HW$ is  the following
\begin{prop}\label{prop:operations}
Let $M$ be a Liouville domain and $L$ be an asymptotically conical exact Lagrangian in it, and assume that $\HW(M, L)$ has exponential algebraic growth. Then we have:
\begin{itemize}
\item [A)] The Liouville domain ${M}'$ obtained by attaching subcritical handles to $M$ has exponential algebraic growth of $\HW$. More precisely, if the attachments are made away from $L$ (so that $L$ survives as an asymptotically conical exact Lagrangian submanifold of ${M}'$) then $\HW({M}',L)$ has exponential algebraic growth. 
\item [B)] If  $M$ is the unit disk bundle of a closed orientable manifold $Q^n$ whose fundamental group grows exponentially, and ${M}'$ is obtained by a plumbing whose graph is a tree and one of the vertices is $M$, then ${M}'$ has exponential algebraic growth of $\HW$. More precisely, if $L_q$ is a unit disk fibre in $M$ and the plumbing is done away from $L_q$ then $\HW({M}',L_q)$ has exponential algebraic growth.
\end{itemize}
\end{prop}
This result essentially says that plumbing and subcritical surgeries are topological operations that preserve exponential algebraic growth of $\HW$, and will allow us to construct many examples of Liouville domains which admit asymptotically conical exact Lagrangian disks with exponential algebraic growth of $\HW$.

The exponential algebraic growth of our examples stems from the algebraic growth of the homology of the based loop space $H_{*}(\Omega Q)$ equipped with the Pontrjagin product, where $Q$  is a compact manifold. In fact, we will only use the degree $0$ part whose algebraic growth is that of $\pi_1(Q)$. 

\begin{rem}
The exponential algebraic growth of symplectic homology always vanishes since its product is commutative. Thus our approach is specifically designed for the open string case.    
\end{rem}

In order to obtain our main results we will bound the topological entropy of Reeb flows from below in terms of the algebraic growth of $\HW(M,L)$. For that we will use the crucial fact that the spectral number $c: \HW(M,L) \rightarrow \R_+$ defined by $c(x) = \inf\{a \in \R \, |\, x \in \im \, i_a   \}$ is subadditive, i.e. $c(x \star y) \leq c(x) + c(y)$ for all $x,y \in \HW(M,L)$.  \color{black} It follows (see Proposition~\ref{prop:alg_symp}) that for any finite $ S \subset \HW(M,L)$ we have 
\begin{equation*}\label{algebraic_symplectic}
\Gamma^{\symp}(M,L) \geq \frac{1}{\rho(S)}\Gamma^{\alg}_S(M,L),
\end{equation*}
where $\rho(S) = \max_{s\in S} c(s)$. 
By using that $\HW(M,L \to L_1)$ is a module over $\left(\HW(M,L), \star \right)$, this lower bound can be extended to $\Gamma^{\symp}(M,L \to L_1)$ for all $L_1$ that are exact Lagrangian isotopic to $L$, see Lemma \ref{lemmaprelim}.  In other words, exponential algebraic growth of $\HW(M,L)$ implies positive symplectic growth of $\HW(M,L \to L_1)$. This, combined with ideas from \cite{Alves-Legendrian}, leads to  
\begin{thm} \label{theorementropy}
Let $L$ be an asymptotically conical exact Lagrangian on a Liouville domain $M=(Y,\omega,\lambda)$, $\Sigma:= \partial Y$ and $ \alpha_{M}:= \lambda|_\Sigma$. We denote by $\xi_{M}:=\ker  \alpha_{M})$ the contact structure induced by $M$ on $\Sigma$.
Assume that there is a finite set $ S \subset \HW(M, L)$ such that $\Gamma^{\alg}_S(M,L) >0$ and that $\Lambda = \partial L$ is a sphere. Then, for every contact form $\alpha$ on $(\Sigma ,\xi_{M})$ the topological entropy of the Reeb flow $\phi_{\alpha}$ is positive. Moreover, if $\mathsf{f}_{\alpha}$ is the function such that $\mathsf{f}_{\alpha}\alpha_{M} = \alpha$ then
\begin{equation*}
h_{\topo}(\phi_{\alpha}) \geq \frac{\Gamma^{\alg}_S(M,L)}{\rho(S) \max(\mathsf{f}_{\alpha})}.
\end{equation*}
\end{thm}  \color{black}

Our paper is organised as follows.
In section \ref{sec:Floer} we consider the algebraic growth and the growth of filtered directed systems in general, and then we recall the definition of wrapped Floer homology together with its product structure. In section \ref{sec:Viterbo} we present the construction of the Viterbo map and derive some of its properties. Section \ref{sec:entropy} establishes implications of the growth properties of $\HW$ to topological entropy. In section \ref{sec:ring_module} we recall the computation of the algebra structure of the Floer homology of unit disk bundles and in section \ref{sec:top_operations} we give a proof of the invariance of $\HW$ under subcritical handle attachment, recollect a result on $\HW$ of plumbings and prove Proposition \ref{prop:operations}. Finally, in section \ref{sec:constructions}, we construct our examples and prove the main theorems. The Appendix contains a construction of exact Lagrangian cobordisms used in the paper. 

\acknowledgementname: Most of this work was done when the second author visited the Universit\'e of Neuch\^atel supported by the Erasmus mobility program, and the first author visited the Universit\"at M\"unster supported by the SFB/TR 191. This work greatly benefited from discussions with Felix Schlenk and Peter Albers: we thank them for their interest in this work and their many suggestions.  We also thank Lucas Dahinden for carefully reading the manuscript.

\section{Wrapped Floer homology and its growth}\label{sec:Floer}
As explained in the introduction, two features of wrapped Floer homology are crucial in this paper. 

 First, its natural filtration by action gives the wrapped Floer homology $\HW$ the structure of a filtered directed system and allows one to define the spectral value of elements of $\HW$. These give rise to the notion of symplectic growth\footnote{This was explicitly observed in \cite{Mclean2015} although it is implicit in \cite{FrauenfelderSchlenk2006,MacariniSchlenk2011}.} of $\HW$; this is explained in Section \ref{sec:wrapped}.  
 
 Second, the product structure of $\HW$ gives it the structure of an algebra and gives rise to the notion of algebraic growth of $\HW$. This is explained in Section \ref{sec:productonHW}.  
The link between these notions is given by the crucial fact that the spectral number is subadditive with respect to the product structure on $\HW$, see also Section \ref{sec:productonHW}.

 We first recall the relevant algebraic notions and deduce some direct consequences.

\subsection{Algebraic growth and growth of filtered directed systems} 
Fix a field $K$. We use the convention that $\log(0) := 0$. 

\subsubsection{Filtered directed systems and growth}\label{subsubsec:fds}

\begin{defn}
A \textit{filtered directed system} over $\R_{+} = [0,\infty)$ or for short \textit{f.d.s.} is a pair $(V,\pi)$ where
\begin{itemize}
\item $V_t$, $t \in [0,\infty)$, are finite dimensional $K$-vector spaces.
\item $\pi_{s \rightarrow t}: V_s \rightarrow V_t$, for $s\leq t$ are homomorphisms (\textit{persistence homomorphisms}), such that 
$\pi_{s \rightarrow t}\circ\pi_{r \rightarrow s}=\pi_{r \rightarrow t}$ for $r\leq s\leq t$, and $\pi_{t\rightarrow t} = \id_{V_t}$ for all $t \in \R_{+}$.
\end{itemize}
\end{defn}

\color{black}


\color{black}
Let $\mathfrak{J}$ be the smallest vector space of $\bigoplus_{t\in \R_+} V_t$ containing $\bigcup_{s \leq t} \{\pi_{s\rightarrow t}(x_s) - x_s\}$. 
The \textit{direct limit} $\dlim V$ of $V$ is defined by $\dlim V := \bigoplus_{t\in \R_+} V_t/\mathfrak{J}$. The inclusions $V_t \hookrightarrow  \bigoplus_{t\in \R_+} V_t$ induce maps to $\dlim V$ which we denote by $i_t$.  
The \textit{spectral number} $c_V$, or just $c$ if the context is clear, of an element $x \in \dlim V$ is
\[ 
c_V(x) := \inf\{t \in [0,\infty) \, | \, \exists x_t \in V_t \text{ such that } i_t(x_t) = x \}.
\]
It is clear from the definition of $c_V$ that if $x_1,..,x_n \in V$ and $k_1,...,k_n\in K$ we have
\begin{equation} \label{spectralsuminequality}
c_V\bigg( \sum_{i=1}^n k_i x_i  \bigg) \leq \max_{1 \leq i \leq n}{c_V(x_i)}.
\end{equation}

\begin{defn}\label{fds_growth}
Let $d_t^V := \dim \{x \, | \, c_V(x) \leq t\}$.
The \textit{exponential growth rate} of the f.d.s. $V$ is
\[
\widetilde{\Gamma}(V) := \limsup_{t \rightarrow \infty}\frac{1}{t} \log d^V_t.
\]
We say that $V$ has  exponential growth if $0 < \widetilde{\Gamma}(V) < \infty$.
\end{defn}

\begin{defn}
A \textit{morphism} between f.d.s. $(V,\pi)$ and $(V^{'},\pi^{'})$ is a collection of homomorphisms 
$f = (f_t)_{t\in[0, \infty)}$, $f_t : V_t \rightarrow V^{'}_t$, that are compatible with respect to the persistence homomorphisms: 
\begin{equation}\label{morph}
f_t \circ \pi_{s \rightarrow t} = \pi^{'}_{s\rightarrow t}\circ f_s.
\end{equation}
An \textit{asymptotic morphism} is a collection of homomorphisms $f_t : V_t \rightarrow V^{'}_t$, $t\in (K, \infty)$, for some $K > 0$ such that \eqref{morph} holds for $K<s<t$. 
\end{defn}

\color{black}


\color{black}

Let $(V,\pi)$ be a f.d.s. and \color{black} $\eta \geq 1$. \color{black}
   We can dilate $V$ by $\eta$ to a filtered directed system $(V(\eta), \pi(\eta))$ given by $V({\eta})_t = V_{\eta t}$, $\pi(\eta)_{s\rightarrow t}= \pi_{\eta s \rightarrow \eta t}$. It follows that
$\pi$ gives rise to a canonical morphism $\pi[\eta]:V \rightarrow V({\eta})$ by $\pi[\eta]_t = \pi_{t \to \eta t}$. For a morphism $f: V\rightarrow W$ we get a dilated morphism $f(\eta):V(\eta) \rightarrow W(\eta)$ by setting $f(\eta)_t = f_{\eta t}$.
 
\begin{defn}\label{interl}
Let $(V,\pi_V)$ and $(W,\pi_W)$ be f.d.s. We call them  \textit{$(\eta_1,\eta_2)$-interleaved}, or \textit{interleaved}, if there are asymptotic morphisms
$f:V\rightarrow W(\eta_1)$ and $g:W \rightarrow V(\eta_2)$ for two real numbers $\eta_1, \eta_2 \geq 1$ such that  
\[f(\eta_2)\circ g = \pi_{W}[\eta_1 \eta_2] \ \text{   and   }  \ g(\eta_1)\circ f = \pi_{V}[\eta_1 \eta_2].\]
\end{defn}

The direct limits of interleaved f.d.s. are isomorphic. It is also easy to see the following
\begin{lem}\label{interl_growth}
Let $V$ and $W$ be $(\eta_1, \eta_2)$-interleaved for some $\eta_1,\eta_2 \geq 1$. Then
\begin{equation*}
\widetilde{\Gamma}(V) \leq \eta_1 \widetilde{\Gamma}(W) \  \mbox{ and } \  \widetilde{\Gamma}(W) \leq \eta_2 \widetilde{\Gamma}(V).
\end{equation*}
\end{lem}



\begin{rem}

The notion of interleaving comes from the theory of persistence modules (see \cite{PolterovichShelukhin2016} for applications of persistence modules and interleaving distance in symplectic geometry).
\end{rem}

\subsubsection{Algebras and their algebraic growth} \label{subsec:algebra} 

\

\color{black}
We recall from the introduction the definition of the algebraic growth of a $K$-algebra $A$ and a finite subset $S \subset A$.
Given $j\geq 0$ let $N_S(j) = \{ a \in A \, \mid \, a = s_1\star s_2 \star \cdots \star s_j; \, s_1, \dots, s_j \in S \}$; i.e.\ $N_S(j)$ is the set of elements of $A$ that can be written as a product of $j$, not necessarily distinct, elements of $S$.
 We define $W_S(n) \subset A$ to be the smallest $K$-vector space that contains the union $\bigcup_{j=1}^{n} N_S(j)$.
 The exponential algebraic growth rate of the pair $(A,S)$ is defined as
\[
\Gamma^{\alg}_S(A) = \limsup_{n\to \infty}\frac{1}{n}\log\dim_K W(n) \in [0,\infty).
\]
\color{black}

We will need the following definition.
\begin{defn} \label{streched}
Let $M$ be a module over an algebra $A$ with scalar multiplication denoted by~$\ast$. The module $M$ is called \textit{stretched} if there exists an element $m_0 \in M$ such that for all elements $a \neq 0  \in A$ we have $a\ast m_0 \neq 0$. An element $m_0 \in M$ satisfying this condition is called a \textit{stretching} element.
\end{defn}

\color{black}
In the following let $V$ be a filtered directed system and assume that the vector space $A= \dlim V$ has a $K$-algebra structure with multiplication $\star$. We do not assume that $A$ is finitely generated. Furthermore, let $W$ be a filtered directed system, such that $M= \dlim W$ is a module over $A$ with multiplication $\ast$,  i.e. a module over $(A,\star)$ with scalar multiplication $\ast$ which is compatible with the $K$-vector space structure of $A$ and $M$. \color{black}
\color{black}

Furthermore assume that the spectral numbers $c_V$ and $c_W$ are \textit{subadditive} with respect to $\star$ and $\ast$, i.e.
\begin{align}\label{subadd_ring}
c_V(a \star b) \leq c_V(a) + c_V(b), \text{ for all } a,b \in A,
\end{align}
and
\begin{align}\label{subadd_module}
c_W(a \ast m) \leq c_V(a) + c_W(m), \text{ for all }a\in A \text{ and } m\in M.
\end{align}


\color{black}

\begin{lem}\label{finite_growth}
Let $V$ be a f.d.s. such that $A= \dlim V$ has a $K$-algebra structure with multiplication~$\star$, and assume that $c_V$ is subadditive with respect to $\star$. Then for every finite subset $S \subset A$ we have 
\[
\widetilde{\Gamma}(V) \geq \frac{1}{\rho(S)}\Gamma^{\alg}_S(A),
\]
where $\rho(S) = \max_{x\in S}c_V(x)$.
\end{lem}

\textit{Proof: }
From the subadditivity of $c_V$ with respect to $\star$ it follows that if $a = s_1 \star s_2 \star \cdots \star s_n$, $s_i \in S$, we have
\[ 
c_V(a) = c_V(s_1 \star \cdots \star s_n) \leq c_V(s_1) + \cdots + c_V(s_n) \leq \rho(S) n.
\]
It then follows from \eqref{spectralsuminequality} that $W(n) \subset \{x \in A \, \mid c(x) \leq \rho(S)n \}$. We thus conclude that
\begin{equation*}
\begin{split}
\Gamma^{\alg}_S(A) &=  \limsup_{n\to \infty}\frac{1}{n}\log \dim W(n)  \\
&\leq \limsup_{n \to \infty}\frac{1}{n}\log \, \dim \{x  \mid  c(x) \leq \rho(S) n\} \leq 
\rho(S) \widetilde{\Gamma}(V).
\end{split}
\end{equation*}
\qed

\color{black}

\color{black}
\begin{lem}\label{mod_growth} Let $V$ and $W$ be f.d.s. and assume that the vector space $A= \dlim V$ has an $K$-algebra structure with multiplication $\star$, and that $M := \dlim W$ has the structure of a module over $A$ with multiplication $\ast$.
Assume that $c_V$ and $c_W$ are subadditive with respect to $\star$ and $\ast$, respectively, and that $M \neq 0 $ is a stretched module over the algebra $A$.  
Then 
\begin{equation} \label{firsteq}
\widetilde{\Gamma}(W) \geq \widetilde{\Gamma}(V).
\end{equation}
Moreover, for every finite set  $S \subset A$ we have
\begin{equation} \label{seceq}
\widetilde{\Gamma}(W) \geq \frac{1}{\rho(S)}\Gamma^{\alg}_S(A).
\end{equation}
\end{lem} \color{black}
\textit{Proof: }
Take a stretching element $m_0 \neq 0$ in $M$. We have $a \ast m_0 \neq b \ast m_0$ for $a \neq b$, $a,b \in A$. In particular $a \mapsto a \ast m_0$ is an injective homomorphism from $A$ to $M$.    
Therefore, by \eqref{subadd_module}, $d_t^V = \dim \{a \in A \, |\, c_V(a) \leq t \} \leq \dim \{m \in M \, | c_W(m) \leq t + c_W(m_0) \} = d_{t+c_W(m_0)}^W$, for all $t >0$. 
We then get
\begin{equation*} 
\begin{split}
\widetilde{\Gamma}(V) &= \limsup_{t \rightarrow \infty} \frac{\log \, d_t^V}{t} \leq \limsup_{t \rightarrow \infty} 
\frac{\log \, d^W_{t+c_W(m_0)}}{t} \\ &= 
\limsup_{t \rightarrow \infty} \frac{\log \, d^W_{t+c_W(m_0)}}{t +c_W(m_0)} \frac{t+c_W(m_0)}{t} = 
\widetilde{\Gamma}(W).
\end{split}
\end{equation*}
This proves \eqref{firsteq}. Inequality \eqref{seceq} is obtained by combining \eqref{firsteq} with Lemma \ref{finite_growth}.
\qed


In order to get results on entropy, we will need the following notions.

\begin{defn}
Let $\mathcal{W} = W(i)_{i \in I}$ be a family of f.d.s. with direct limits $M(i)$ that are modules over $A:= \dlim V$.
 We say that the family $M(i)_{i \in I}$ is \textit{uniformly stretched} if there exists a constant $B\geq 0$ such that for every $i \in I$ there exists a stretching element $m_i \in M(i)$ with $c_{M(i)}(m_i) \leq B$. 
\end{defn}
 
\begin{defn} \label{defi:fam_growth}
Let $\mathcal{W} = W(i)_{i \in I}$ be a family of filtered directed systems. The \textit{uniform exponential growth rate } of $\mathcal{W}$ is 
\[
\widetilde{\Gamma}_{i\in I}(\mathcal{W}) := \limsup_{t \rightarrow \infty} \frac{1}{t} \log \left(\inf_{I} d_t^{W(i)}\right).
\]
\end{defn}

\color{black}
\begin{lem}\label{mod_fam_growth}
Let $V$ be a f.d.s. such that $A= \dlim V$ has a $K$-algebra structure with multiplication $\star$. 
Let $\mathcal{W} = W(i)_{i \in I}$ be a family of f.d.s. such that for every $i\in I$ the direct limit $M(i)=\dlim W(i)$  is a module over $A$ with multiplication $\ast(i)$. Assume that $c_V$ is  subadditive with respect to $\star$, that  $c_{W(i)}$ is subadditive with respect to $\ast(i)$ for every $i\in I$, and that the family $M(i)_{i \in I}$ is uniformly stretched over the  algebra $A$.
Then
\begin{equation}
\widetilde{\Gamma}_{i\in I}(\mathcal{W}) \geq \widetilde{\Gamma} (V).
\end{equation}
\end{lem} \color{black}
\textit{Proof: }
Since $M(i)_{i \in I}$ is uniformly stretched there exists $B>0$ such that for every $i \in I$, we can find a stretching element $m_i \in M(i)$ with $c_{M(i)}(m_i) \leq B$.  Hence we have by \eqref{subadd_module} that $d_t^V \leq \inf_{I}d^{W(i)}_{t+B}$ and the result is obtained as in the proof of Lemma \ref{mod_growth}. 
\qed

\subsection{Wrapped Floer homology}\label{subsec:def}

In the following we give the definition and conventions for wrapped Floer homology used in this paper. This Floer type homology theory appeared in \cite{AS-iso} for contangent bundles, and the case of general Liouville domains can be found in \cite{AbouzaidSeidel2010}. We refer to these papers and \cite[Section 4]{Ritter2013} for more details.

\subsubsection{Liouville domains and Lagrangians}

A \textit{Liouville domain} $M=(Y,\omega, \lambda)$ is an exact symplectic manifold $(Y, \omega)$ with boundary $\Sigma = \partial Y$ and a primitive $\lambda$ of $\omega$ such that $\alpha_M = \lambda|_{\Sigma}$ is a contact form on $\Sigma$.
 The \textit{Liouville vector field} X, is given by $i_X\omega = \lambda$ and points outwards along~$\Sigma$. Using the flow of X, one can attach an infinite cone to $M$ along $\Sigma$ that gives the \textit{completion} $\widehat{M}:=\left(\widehat{Y}, \widehat{\omega}, \widehat{\lambda}\right)$ of $M$ with $\widehat{Y} = Y \cup_{\Sigma} \left([1, \infty) \times \Sigma \right)$, $\widehat{\lambda}|_{Y} = \lambda$, $\widehat{\lambda}|_{[1, \infty) \times \Sigma} = r\alpha_M$, and $\widehat{\omega}=d \widehat{\lambda}$.

\begin{rem}In order to simplify notation, we will usually write $M$ and $\widehat{M}$ instead of $Y$ and $\widehat{Y}$, respectively, as the domain of Hamiltonian functions or the target space of Floer trajectories. This does not cause any confusion since the smooth manifolds $Y$ and $\widehat{Y}$ are part of  the data defining $M$ and $\widehat{M}$, respectively. Similarly, when we write $\Sigma=\partial M$ it should be understood as $\Sigma= \partial Y$. \end{rem}

Let $\mathsf{f}: \partial Y \rightarrow (0, \infty)$ be a smooth function. Let $Y_{\mathsf{f}} = \widehat{Y} \setminus \{(r,x) \, | \, r > \mathsf{f}(x), x \in \partial Y\}$. It is easy to see that $M_{\mathsf{f}} = (Y_{\mathsf{f}},\widehat{\omega}|_{Y_{\mathsf{f}}},\widehat{\lambda}|_{Y_{\mathsf{f}}})$ is a Liouville domain. For example, given $\delta >0$ we denote by $M_{1+\delta}$  the Liouville domain $(Y_{1+\delta},\omega_{1+\delta},\lambda_{1+\delta})$ embedded in $\widehat{M}$ defined by ${Y}_{1+\delta} = Y \cup_{\Sigma} \left([1, 1+\delta] \times \Sigma\right)$, 
$\omega_{1+\delta} =  \widehat{\omega}|_{Y_{1+\delta}}$, $\lambda_{1+\delta} =  \widehat{\lambda}|_{Y_{1+\delta}}$.

In our paper we only consider Liouville domains that have vanishing first chern class $c_1(M) \subset \Ho^2(M;\Z)$.

\ 

We consider Lagrangians $(L, \partial L)$ in $(M, \Sigma)$ that are exact, i.e. $\lambda|{L} = df$, and  that satisfy 
\begin{equation}\label{asympt}
\begin{split}
 &\Lambda = \partial L \text{ is a Legendrian submanifold in } (\Sigma,\xi_M), \\
 & L \cap [1-\epsilon, 1] \times \Sigma = [1-\epsilon, 1] \times \Lambda \text{ for a sufficiently small } \epsilon >0 .
\end{split} 
\end{equation}
We will call a Lagrangian that satisfies \eqref{asympt} \textit{asymptotically conical}. We can extend it naturally to an exact Lagrangian $\widehat{L} = L \cup_{\Lambda} ([1,\infty) \times \Lambda)$ in $\widehat{M}$.
We will refer to a Lagrangian in $\widehat{M}$ of this form also as \textit{asymptotically conical (with respect to $M$)}. More generally, given a subset $U \subset \widehat{M}$ we say that $L$ is \textit{conical in $U$} if the Liouville vector field is tangent to  $L\cap \mathrm{int}(U)$ in the interior $ \mathrm{int}(U)$ of $U$. 


\subsubsection{Wrapped Floer homology} \label{sec:wrapped}

For two asymptotically conical exact Lagrangians $L_0$ and $L_1$ in ${M}$ denote by $\mathcal{P}_{L_0 \to L_1} = \{\gamma : [0,1] \rightarrow \widehat{M} \, |\, \gamma(0) \in \widehat{L}_0, \,  \gamma(1) \in \widehat{L}_1\}$
the space of (smooth) paths from $\widehat{L}_0$ to $\widehat{L}_1$. 

Denote by $X_{\alpha_M}$ the Reeb vector field on the boundary $(\Sigma, \xi_{M} = \ker{\alpha_M})$.  
A \textit{Reeb chord of length $T$} of $\alpha_M$ from $\Lambda_0 = \partial L_0$ to $\Lambda_1 = \partial L_1$ is a map $\gamma: [0,T] \to \Sigma$ with $\dot{\gamma}(t) = X_{\alpha_M}(\gamma(t))$ with $\gamma(0) \in \Lambda_0$ and $\gamma(T) \in \Lambda_1$. Denote the set of Reeb chords of length  $<T$ by $\mathcal{T}^{T}_{\Lambda_0 \to \Lambda_1}(\alpha_M)$, and the set of all Reeb chords by $\mathcal{T}_{\Lambda_0 \to \Lambda_1}(\alpha_M)$. The Reeb chord $\gamma$ of length $T$ of $\alpha_M$ from $\Lambda_0$ to $\Lambda_1$ is said to be \textit{transverse} if the subspaces $T_{\gamma(1)}(\phi_{X_{\alpha_M}}^{T}(\Lambda_0))$ and $T_{\gamma(1)} \Lambda_1 $ of $T_{\gamma(1)}\Sigma$ intersect at only one point.
The \textit{spectrum} of the triple $(M,L_0 \to L_1)$, denoted by $\mathcal{S}(M,L_0 \to L_1)$, is the set of lengths of Reeb chords from $\Lambda_{0}$ to $\Lambda_{1}$ in $\Sigma$. It is a nowhere dense set in $[0,\infty)$. 

Given a  contact form $\alpha$ on $(\Sigma, \xi_{M})$ and a pair of Legendrian submanifolds $(\Lambda_0, \Lambda_1)$ on $(\Sigma, \xi_{M})$, we say that the triple $(\alpha,\Lambda_0 \to \Lambda_1)$ is \textit{regular} if all Reeb chords of $\alpha$ from $\Lambda_0$ to $\Lambda_1$ are transverse.
We say that $(M,L_0 \to L_1)$ is \textit{regular} if $(\lambda_{\Sigma}, \Lambda_0 \to \Lambda_1 )$ is regular and $L_0$ and $L_1$ intersect transversely. 

From now on, we assume that for the contact form $\alpha_M$ induced by $M$ on $(\Sigma,\xi_M)$ the triple $(\alpha_M,\Lambda_0 \to \Lambda_1)$  is regular.

An autonomous Hamiltonians $H: \widehat{M} \to \R$ is called  \textit{admissible} if
\begin{itemize} \label{ham}
\item $H < 0$  on $M$,
\item and  there exist constants  $\mu>0$  and  $b\leq - \mu$ such that $H(x,r) = h(r) = \mu r +  b$  on $[1, \infty) \times \partial M$.
\end{itemize}
If $H: \widehat{M} \to \R$ is admissible and satisfies $H(x,r) = \mu r +  b$  on $[1, \infty) \times \partial M$ we say that $H$ is admissible with \textit{slope} $\mu$ \textit{(at infinity)}.

Define the action functional $\mathcal{A}^{L_0 \to L_1}_H = \mathcal{A}_H : \mathcal{P}_{L_0 \to L_1} \rightarrow \R$ by 
\[
\mathcal{A}_H(\gamma) = f_0(x(0)) - f_1(x(1)) + \int_{0}^{1} \gamma^{*}\lambda - \int_{0}^{1} H(\gamma(t))dt, 
\]
where $f_0$ and $f_1$ are functions on $L_0$ and $L_1$ respectively with $df_i = \lambda|_{\widehat{L}_i}, i= 0,1$. 
The critical points of $\mathcal{A}_H$ are Hamiltonian chords from $\widehat{L}_0$ to $\widehat{L}_1$ that reach $\widehat{L}_1$ at time $1$. We define
\[
 \mathcal{T}_{L_0\to L_1}(H):=\crit{\mathcal{A}_H} = \{\gamma \in \mathcal{P}_{L_0 \to L_1} \mid \dot{\gamma}(t) = X_H(\gamma(t))\}, 
\]
and write $\mathcal{T}_{L}(H)$ instead of $\mathcal{T}_{L \to L}(H)$. 
Here $X_H$ is the Hamiltonian vector field defined by $\omega(X_H, \cdot ) = -dH$. 
We call an admissible Hamiltonian \textit{non-degenerate} for $L_0 \to L_1$ if 
all elements in $\mathcal{T}_{L_0\to L_1}(H)$ are non-degenerate, i.e. $\phi^1_{X_H}(\widehat{L_0})$ is transverse to $\widehat{L_1}$. Such a Hamiltonian must have slope $\mu \notin \mathcal{S}(M,L_0 \to L_1)$. Note that every admissible Hamiltonian can be made non-degenerate for $L_0 \to L_1$ after a generic perturbation (\cite[Lemma 8.1]{AbouzaidSeidel2010}). We denote by 
 \begin{equation}
 \mathcal{H}_{\reg}(M,L_0 \to L_1)
\end{equation}
the set of admissible Hamiltonians which are non-degenerate for $L_0 \to L_1$.
For a Hamiltonian $H\in \mathcal{H}_{\reg}(M,L_0 \to L_1) $ all elements in $\mathcal{T}_{L_0 \to L_1}(H)$ have their image contained in $M$.

 For admissible Hamiltonians $H$ with slope $\mu \notin \mathcal{S}$ that are constant in $M$ away from the boundary, depend on $r$ and increase sharply near $\partial M$, $\mathcal{T}_{L_0\to L_1}(H)$ corresponds to  $\mathcal{T}^{\mu}_{\Lambda_0 \to \Lambda_1}(\alpha_M)$ and intersection points of $L_0$ and $L_1$ in $M$. If $(M,L_0 \to L_1)$ is regular, such Hamiltonians belong to the set $ \mathcal{H}_{\reg}(M,L_0 \to L_1)$. 

 If  $(\alpha_M,  \Lambda_0 \to  \Lambda_1)$ is regular but $(M,L_0 \to L_1)$ is not,  we can take $H \in  \mathcal{H}_{\reg}(M,L_0 \to L_1)$ to be a $C^2$-small negative function away from the boundary  of $M$, and to depend only on $r$ and  increase sharply near $\partial M$.  Then, $\mathcal{T}_{L_0\to L_1}(H)$ will correspond to $\mathcal{T}^{\mu}_{\Lambda_0 \to \Lambda_1}(\alpha_M)$ and intersection points of $L_0$ and $L_1$ in $M$ that are not destroyed by the Hamiltonian flow of $H$.

An almost complex structure $J$ on $\left((0,\infty) \times \partial M, \lambda = r\alpha_M \right)$ is called \textit{cylindrical} if it preserves $\xi_{M}= \ker \alpha_M$, if $J|_{\xi_{M}}$ is independent of $r$ and compatible with $d(r\alpha_M)|_{\xi_{M}}$, and if $J X_{\alpha_M}=r\partial_{r} $. 
In the following we take almost complex structures $J$ on $\widehat{M}$ that are \textit{asymptotically cylindrical}, i.e. cylindrical on $[r, \infty) \times \partial M$ for some $r>1$. 
The $L^2$-gradient of the action functional with respect to the Riemannian metric given by $d\lambda(J \cdot,\cdot) = g(\cdot, \cdot)$
is given by
\[
\nabla \mathcal{A}_H(\gamma) = -J(\gamma)\left(\partial_t \gamma - X_H(\gamma) \right), 
\]
and we interpret the negative gradient flow lines as Floer strips
\begin{equation}
\begin{split}\label{Floer}
&u:\R \times [0,1] \rightarrow \widehat{M}, \\
&\overline{\partial}_{J,H}(u)=\partial_s u  + J(u)(\partial_t u - X_H(u) ) = 0, \\
&u(\cdot,0) \in \widehat{L}_0, \mbox{ and } u( \cdot ,1) \in \widehat{L}_1. 
\end{split}
\end{equation}
We define the moduli space of parametrized Floer strips connecting two critical points $x$ and $y$ of $\mathcal{A}_H$
\begin{equation}
\begin{split}
\widetilde{\mathcal{M}}(x,y, H, J) = \{ u: \R \times [0,1] \rightarrow \widehat{M} \, | \, u \text{ satisfies \eqref{Floer} }, \lim_{s\rightarrow  -\infty} = x \mbox{ and } \lim_{s\rightarrow  +\infty} = y \}. 
\end{split}
\end{equation}
There is a natural $\R$-action on $\mathcal{M}(x,y, H, J)$ coming from the translations in the domain.
Letting $\widetilde{\mathcal{M}}^1(x,y, H, J) $ be the set of elements of $\widetilde{\mathcal{M}}(x,y, H, J) $ that have Fredholm index 1 we write
\begin{equation}
\begin{split}
& \mathcal{M}^0(x,y, H, J) :=  \widetilde{\mathcal{M}}^1(x,y, H, J)  / \R. 
\end{split}
\end{equation}
where the quotient is taken with respect to the $\R$-action mentioned above.The \textit{energy} of an element $u$ is
\[
E(u) := \int_{-\infty}^{\infty}|\nabla \mathcal{A}_H|_{L^2}^2 \, ds 
= \mathcal{A}_H(x)-\mathcal{A}_H(y).
\]

For a generic $J$ and non-degenerate admissible $H$ define the wrapped Floer chain complex
\[
\CW(H,L_0 \to L_1) = \bigoplus_{x \in \crit(\mathcal{A}_H)} \Z_2 \cdot x,
\]
with differential $\partial : \CW(H,L_0 \to L_1) \rightarrow \CW(H,L_0 \to L_1)$ given by  
\[
\partial (x) = \sum_{y \in \crit(\mathcal{A}_H)} \#_{\Z_2} \mathcal{M}^0(x,y,H,J) \cdot y.
\]
For generic $J$ the differential is well-defined and moreover $\partial^2= 0$. 
For simplicity we will write $\CW(H)$ instead of $\CW(H,L_0 \to L_1)$ when there is no possibility of confusion.
In this paper we are not concerned with gradings in $\CW$.
The homology of $(\CW(H,L_0 \to L_1),\partial)$ is called the wrapped Floer homology of $(H,L_0 \to L_1)$ and is denoted by $\HW(H;L_0 \to L_1)$, or in short $\HW(H)$.

Next we consider continuation maps. 
Let $H_-$ and $H_+$ be non-degenerate admissible Hamiltonians with $H_+(x) \geq H_-(x)$ for all $x \in \widehat{M}$, in short $H_+ \succ H_-$. Take an increasing homotopy through admissible Hamiltonians $(H_s)_{s\in \R}$,  $\partial_s H_s \geq 0$,  with $H_s = H_{\pm}$ near $\pm \infty$.
For elements in $\mathcal{M}^0(x_-,x_+,H_s,J)$, i.e. Floer strips 
\begin{equation}
\begin{split}\label{contin}
&u:\R \times [0,1] \rightarrow \widehat{M}, \\
&\overline{\partial}_{J,H_s}(u) := \partial_s u + J(\partial_t u - X_{H_{s}}(u)) = 0, \\
&\lim_{s\rightarrow \pm \infty} u(s,t) = x_{\pm}, \\
&u(\cdot, 0) \in \widehat{L}_0, \mbox{ and } u(\cdot, 1) \in \widehat{L}_1, 
\end{split}
\end{equation}
with Fredholm index 0 connecting $x_- \in \crit(\mathcal{A}_{H_-})$ and $x_+ \in \crit(\mathcal{A}_{H_+})$, the action difference 
is 
\[\mathcal{A}_{H_-}(x_-) -\mathcal{A}_{H_+}(x_+) = E(u) + \int_{\R \times [0,1]}\partial_s H_s(u). 
\]
Hence the action decreases under the continuation maps
\[ 
\iota^{H_-,H_+} : \CW(H_-) \rightarrow \CW(H_+),
\]
given by 
\[
\iota^{H_-,H_+} (x_-) = \sum_{x^+ \crit(\mathcal{A}_{H_+})} \#_{\Z_2} \mathcal{M}^0(x_-,x_+,H_s,J)  \cdot x_+. 
\]
Define the wrapped Floer homology $\HW(M,L_0 \to L_1) := \dlim_H\HW(H;L_0 \to L_1)$, where the direct limit is taken over all $H\in \mathcal{H}_{\reg}(M,L_0 \to L_1) $. 

\begin{defn}
The homology $\HW(M,L_0 \to L_1)$ is the direct limit of the filtered directed system $\widetilde{\HW}(M,L_0 \to L_1) = \left(\HW^a(M,L_0 \to L_1)\right)_{a \in (0,\infty)}$. Here 
\[\HW^a(M,L_0 \to L_1) := \dlim_{H} {\HW^{a}(H;L_0 \to L_1)},
\]
where $\HW^{a}(H;L_0 \to L_1) $ is the homology of the Floer chain complex restricted to critical points of action less than $a$. The persistence maps $\iota_{a \rightarrow b}: {\HW^a(M,L_0 \to L_1)} \to {\HW^b(M,L_0 \to L_1)}$ are induced by the natural maps $\HW^a(H,L_0 \to L_1) \rightarrow \HW^b(H,L_0 \to L_1)$ that come from inclusions. We write $\iota_{a}: {\HW^a(M,L_0 \to L_1)} \to {\HW(M,L_0 \to L_1)}$ for the induced map from $\HW^a(M,L_0 \to L_1)$ to the direct limit ${\HW(M,L_0 \to L_1)}$. 
\end{defn}

\color{black}
Let $H'\geq H$ be Hamiltonians in $ \mathcal{H}_{\reg}(M,L_0 \to L_1) $ and $b\geq a$. Let $ \iota^{H,H'}_{a \to b}: \HW^a(H, L_0 \to L_1) \to \HW^b(H', L_0 \to L_1)$ by the continuation map induced by any non-decreasing homotopy from $H$ to $H'$.  In case $b=+\infty$ we write $ \iota^{H,H'}_{a}: \HW^a(H, L_0 \to L_1) \to \HW(H', L_0 \to L_1)$.

\begin{rem} \label{rem:stationary}
Notice that $ \iota^{H,H}_{a \to b}$ is the map induced by the chain level inclusion $\CW^a(H,L_0 \to L_1) \hookrightarrow \CW^b(H,L_0 \to L_1)$. For this reason we will also denote this inclusion also by $\iota^{H,H}_{a \to b}: \CW^a(H,L_0 \to L_1) \hookrightarrow \CW^b(H,L_0 \to L_1)$.
\end{rem}

By the construction of  $\widetilde{\HW}(M,L_0 \to L_1)$ presented above we have for every number $a\geq 0$ and Hamiltonian $H\in \mathcal{H}_{\reg}(M,L_0 \to L_1) $ a map 
\begin{equation}
\chi^H_{a \to a}:\HW^a(H, L_0 \to L_1) \to \HW^a(M, L_0 \to L_1).
\end{equation}
This allows us to define for every Hamiltonian $H\in \mathcal{H}_{\reg}(M,L_0 \to L_1) $ and numbers $b\geq a$ the map 
\begin{equation*}
\chi^H_{a \to b}:=\iota_{a\to b} \circ \chi^{H}_{a \to a} = \HW^a(H, L_0 \to L_1) \to \HW^b(M, L_0 \to L_1).
\end{equation*}
Using functoriality properties of continuation maps it is straightforward to check that $$ \chi^H_{a \to b} =\chi^H_{b \to b} \circ  \iota^{H,H}_{a \to b }. $$  For simplicity, in the case $b=+\infty$ we write $$  \chi^H_{a } = \iota_{a} \circ \chi^{H}_{a \to a} :  \HW^a(H, L_0 \to L_1) \to \HW(M, L_0 \to L_1).$$

For each $H\in  \mathcal{H}_{\reg}(M,L_0 \to L_1) $ we also have a map $$ \chi^H:  \HW(H, L_0 \to L_1) \to \HW(M, L_0 \to L_1).$$ To define it, we first notice that since $ \mathcal{T}_{L_0 \to L_1}(H)$ is a finite set we can choose a number $ a_H>\max_{x \in \mathcal{T}_{L_0 \to L_1}(H)}\{\mathcal{A}(x)\}$. For this choice of $a_H$ the chain complexes $ (\CW(H, L_0 \to L_1),\partial) $ and $ ( \CW^{a_H}(H, L_0 \to L_1),\partial)$ are identical, and we get  $ \HW(H, L_0 \to L_1) =  \HW^{a_H}(H, L_0 \to L_1)$. We then define $\chi^H := \chi^H_{a_H}$. It is an elementary exercise to check that the definition of  $\chi^H $ does not depend on the choice of $a_H>\max_{x \in \mathcal{T}_{L_0 \to L_1}(H)}\{\mathcal{A}(x)\}$. In the same way we can construct for each $b >\max_{x \in \mathcal{T}_{L_0 \to L_1}(H)}\{\mathcal{A}(x)\}$ a map $$ \chi^H_{\to b}:  \HW(H, L_0 \to L_1) \to \HW^b(M, L_0 \to L_1).$$

These maps are useful for the study of spectral numbers done in the next section. We will need the identity 
\begin{equation}\label{eq:naturality'}
 \chi^H_a = \chi^{H'} \circ \iota^{H,H'}_a,
 \end{equation}
which is established in an elementary way from the functoriality properties of continuation maps. In particular, we have
\begin{equation}\label{eq:lepolepo}
 \chi^H = \chi^{H'} \circ \iota^{H,H'},
 \end{equation}
 and
\begin{equation}\label{eq:naturality}
 \chi^H_a = \chi^{H} \circ \iota^{H,H}_a.
 \end{equation}

 We will now define the symplectic growth rate of $\HW$.
\begin{defn} \label{defi:Gamma_symp}
The exponential symplectic growth rate $\Gamma^{\symp}(M,L_0 \to L_1)$ is  defined by 
\begin{equation}
\Gamma^{\symp}(M,L_0 \to L_1) :=  \limsup_{a \to \infty} \frac{\log (\dim \im \ \iota_a)}{a} =  \widetilde{\Gamma}(\widetilde{\HW}(M,L_0 \to L_1)).
\end{equation}
Analogously, given a family $(L_i)_{i\in I}$ of asymptotically conical exact Lagrangians in $M$ we define $\Gamma_{i\in I}^{\mathrm{symp}}(M,L_0 \to L_{i}) := \widetilde{\Gamma}_{i\in I}(\widetilde{\HW}(M,L_0 \to L_i)_{i\in I})$, where $\widetilde{\Gamma}_{i\in I}(\widetilde{\HW}(M,L_0 \to L_i)_{i\in I})$  is defined as in Definition \ref{defi:fam_growth}.
\end{defn}

\subsubsection{Spectral numbers in $\HW$}

\begin{defn}
As ${\HW(M,L_0 \to L_1)}$ is the direct limit of the the f.d.s. $\widetilde{\HW}(M,L_0 \to L_1) $, we define the spectral number $c$ of elements of ${\HW(M,L_0 \to L_1)}$ via the recipe given in Section~\ref{subsec:algebra}.
\end{defn}

We now present an equivalent definition of $c$ which is more geometrical. Given $H\in \mathcal{H}_{\reg}(M,L_0 \to L_1) $, and a cycle $w \in \CW(H, L_0 \to L_1)$ we denote by $[w] \in \HW(H, L_0 \to L_1)$ the homology class of $w$ in $\HW(H, L_0 \to L_1)$. The cycle $w$ can be expressed in a unique way as a sum of elements of $\mathcal{T}_{L_0 \to L_1}(H)$ and we denote by $\mathcal{A}(w)$ the maximum of the actions of these elements. 

 If $w' \in \CW^a(H, L_0 \to L_1)$, then it can be expressed in a unique way as a sum of elements in  $\mathcal{T}^a_{L_0 \to L_1}(H)$. This expression is identical to the one of $\iota_a^{H,H}(w')$, from what we conclude
\begin{equation}
\mathcal{A}(\iota_a^{H,H}(w')) < a \mbox{ for all } w' \in \CW^a(H, L_0 \to L_1).
\end{equation}
 The right hand side in the following identity is often taken as the definition of the spectral number~$c(\mathfrak{h})$.
\begin{lem} \label{lem:minimax}
For a homology class $\mathfrak{h}\in\HW(M, L_0 \to L_1)$ we have
\begin{equation} \label{eq:minmax}
c(\mathfrak{h}) = \inf_{H\in \mathcal{H}_{\reg}(M,L_0 \to L_1)} \{\mathcal{A}(w) \mid w\in \CW(H, L_0 \to L_1) \mbox{ is a cycle with } \chi^H ([w]) = \mathfrak{h} \}.
\end{equation}
\end{lem}
\textit{Proof:} Let $H\in \mathcal{H}_{\reg}(M,L_0 \to L_1) $ and  $w\in \CW(H, L_0 \to L_1)$  be a cycle with  $\chi^H ([w]) = \mathfrak{h}$. For each $a > \mathcal{A}(w) $  we know that there exists a cycle $w' \in \CW^a(H, L_0 \to L_1)$ such that $\iota_a^{H,H}(w')=w$. Using \eqref{eq:naturality} we obtain  
\begin{equation*}
 \chi^{H}_a([w'])=  \chi^{H}\circ\iota^{H,H}_{a}( [w']) = \chi^{H} ( [w]) = \mathfrak{h}.
\end{equation*}
This implies that $\mathfrak{h}$ is in the image of $\chi^{H}_a$ and thus in the image of $\iota_a$, from what we get $c(\mathfrak{h})\leq a$. Since this is valid for each $a > \mathcal{A}(w) $ we obtain that $c(\mathfrak{h}) \leq \mathcal{A}(w)$, and it follows that
\begin{equation} \label{eq:minmaxll}
c(\mathfrak{h}) \leq \inf_{H\in \mathcal{H}_{\reg}(M,L_0 \to L_1)} \{\mathcal{A}(w) \mid w\in \CW(H, L_0 \to L_1) \mbox{ is a cycle with } \chi^H ([w]) = \mathfrak{h} \}.
\end{equation}

To obtain the reverse inequality let $a > c(\mathfrak{h})$. Then there exists $\beta \in \HW^a(M, L_0 \to L_1)$ such that $\iota_a(\beta) = \mathfrak{h}$. By the construction of $\HW^a(M, L_0 \to L_1)$ we know that there exists $H \in \mathcal{H}_{\reg}(M,L_0 \to L_1)$ and a cycle $w' \in \CW^a(H, L_0 \to L_1)$ such that $\chi^H_{a \to a} ([w'])=\beta$. It follows that $$ \chi^H_a ([w'])=\iota_a \circ \chi^H_{a \to a} ([w']) = \iota_a(\beta)=\mathfrak{h}.$$
Let $w:= \iota^{H,H}_{a}( w')$. By the observation we made before the lemma we have $\mathcal{A}(w) < a$. Using \eqref{eq:naturality} we obtain $$ \chi^H([w])= \chi^H(\iota^{H,H}_a([w'])= \chi^H_a([w'])= \mathfrak{h}. $$
We have shown that for each $a > c(\mathfrak{h})$ there exists $H \in \mathcal{H}_{\reg}(M,L_0 \to L_1)$ and a cycle $w \in \CW(H, L_0 \to L_1)$ such that $\mathcal{A}(w) < a$ and  $\chi^H([w]) = \mathfrak{h}$. It follows that 
\begin{equation}
c(\mathfrak{h}) \geq\inf_{H\in \mathcal{H}_{\reg}(M,L_0 \to L_1)} \{\mathcal{A}(w) \mid w\in \CW(H, L_0 \to L_1) \mbox{ is a cycle with } \chi^H ([w]) = \mathfrak{h} \}.
\end{equation} \qed

\subsubsection{A special type of Hamiltonians} \label{rem:rem1}

Given an  admissible Hamiltonian $H$ in $M$ and a number $a> 0$ we write $H\prec a$ if the slope of $H$ is $<a$.
We first define 
\begin{equation}\label{K(M,L)}
\mathsf{K}(M,L_0 \to L_1) := \max \{ \max\{f_0(x) - f_1(x) \, | \, x \in L_0 \cap L_1\},\,  0\}. 
\end{equation}

For $a > \mathsf{K}(M,L_0 \to L_1)$ a careful choice of a cofinal family of Hamiltonians shows that $\HW^a(M,L_0 \to L_1)$ is isomorphic to $\dlim_{H\prec a}\HW(H; L_0 \to L_1)$, where the direct limit is taken only over all non-degenerate admissible $H$ with slope less than $a$.

To explain  this we first take a collar neighbourhood $\mathfrak{V}=([1-\delta,1] \times \Sigma) \subset M$ of $\partial M$ on which $L_0$ and $L_1$ are conical, and $\lambda$ is given by $r\alpha_M$. Since $a> \mathsf{K}(M,L_0 \to L_1)$ we can choose $ \mathsf{K}(M,L_0 \to L_1)<\mu<a$, such that there is no element in $\mathcal{T}_{\Lambda_0 \to \Lambda_1}(\alpha_M)$ with length in the interval $[\mu,a)$. We now choose an admissible Hamiltonian $H^\mu$ in $\widehat{M}$ with slope $\mu$ such that
\begin{itemize}
\item $H^\mu$ is a negative constant $-\mathrm{k}$ in $M\setminus \mathfrak{V}$, with $\mathrm{k}$ small,
\item $H^\mu$ depends only on $r$ in $\mathfrak{V}$,  and is a convex function of $r$ that increases sharply close to $\partial M$.
\end{itemize} 
If $\mathrm{k}$ is small enough, and $H^\mu$ increases sharply enough close to $\partial M$ then we have
\begin{itemize}
\item[$\blacksquare$] the action of all elements of $\mathcal{T}_{L_0 \to L_1}(H^\mu)$ have action $<a$;
\end{itemize}
see for example \cite[Lemma 9.8]{Ritter2013}.

If $(M,L_0 \to L_1)$ is regular then $H^\mu \in \mathcal{H}_{\reg}(M, L_0 \to L_1)$. In this case we have 
\begin{equation}
\mbox{the set } \mathcal{T}_{L_0 \to L_1}(H^\mu) \mbox{ is in bijective correspondence with }\mathcal{T}^a_{\Lambda_0 \to \Lambda_1}(\alpha_M) \cup (L_0\cap L_1).
\end{equation}
 In case $(\alpha_M,\Lambda_0 \to L_1)$ is regular but $(M,L_0 \to L_1)$  is not, we can make a $C^\infty$-small perturbation of $H^\mu$ inside $M$ that still satisfies $\blacksquare$ and is in $\mathcal{H}_{\reg}(M, L_0 \to L_1)$; for simplicity we still denote this perturbation by $H^\mu$. In both cases,  the reasoning used to prove \cite[Lemma 1.5]{Viterbo1999} gives
\begin{equation} \label{muitoimportante}
\chi^{H^\mu}_a: \HW(H^\mu,L_0 \to L_1)  \to  \HW^a(M,L_0 \to L_1) \mbox{ is an isomorphism.}
\end{equation}
It follows that for  $a > \mathsf{K}(M,L_0 \to L_1)$ we have 
\begin{equation} \label{eq:coisabonita}
\HW^a(M,L_0 \to L_1)\cong \dlim_{H\prec a}\HW(H; L_0 \to L_1)
\end{equation}


\color{black}

\subsection{Algebra and module structures on wrapped Floer homology} \label{sec:productonHW}

\subsubsection{Algebra structure in $\HW$}
Let $L$ be an exact asymptotically conical Lagrangian on a Liouville domain $M$. We endow $M$ with an asymptotically cylindrical almost complex structure as in Section \ref{subsec:def}.  We recall the definition of the triangle product in the wrapped Floer homology $\HW(M,L)$, and follow the conventions of \cite{AS-product}. 

We first define the triangle $\Delta$. One first takes the disjoint union $\mathbb{R} \times [-1,0] \cup \mathbb{R} \times [0,1]$. We identify the points $(s,0^-) \in \mathbb{R} \times [-1,0] $ and $(s,0^+) \in \mathbb{R} \times [0,1] $ for all $s \geq 0$, and denote the resulting space by $\Delta$. Let $p_{\mathrm{sing}}$ be the point in $\Delta$ which comes from the points $(0,0^-) \in \mathbb{R} \times [-1,0]$ and $(0,0^+) \in \mathbb{R} \times [0,1] $. 

The interior of $\Delta$ coincides with $(\mathbb{R} \times (-1,1)) \setminus ((-\infty,0] \times \{0\})$. As $(\mathbb{R} \times (-1,1)) \setminus ((-\infty,0] \times \{0\})$ is a subset of $\mathbb{C}$ we can restrict  the complex structure of $\mathbb{C}$ to $(\mathbb{R} \times (-1,1)) \setminus ((-\infty,0] \times \{0\})$. We then obtain a complex structure $j$ in the interior of $\Delta$. This extends to a complex structure on $\Delta\setminus p_{\mathrm{sing}}$. 
Using again that the interior of $\Delta$ coincides with $(\mathbb{R} \times (-1,1)) \setminus ((-\infty,0] \times \{0\})$, we can define global coordinates $(s,t)$ on $\Delta \setminus p_{\mathrm{sing}}$.

For an admissible Hamiltonian $H$ on $\widehat{M}$, the solutions of the Floer equation on $\Delta$ are maps $u: \Delta \to \widehat{M}$ that satisfy
\begin{equation}
\overline{\partial}_{J,H}(u):= \partial_s u + J(u)(\partial_t u - X_H(t,u))=0.
\end{equation}
  We write $\widehat{H} =2H  \in C^{\infty}(M)$.

Given $x_1, x_2\in \mathcal{T}_{L}(H)$ and $y \in  \mathcal{T}_{L}(\widehat{H})$ we let 
$\mathcal{M}(x_1,x_2;y,L,J)$ be the space of maps $u: \Delta \to \widehat{M}$ that satisfy $\overline{\partial}_{J,H}(u)=0$ and such that 
$u(z) \in L$ for all $z\in \partial(\Delta)$, $\lim_{s\to -\infty} u(s,t-1)=x_1(t)$ for $t \in [0,1]$, $\lim_{s\to -\infty} u(s,t)=x_2(t)$ for $t \in [0,1]$, and $\lim_{s\to +\infty} u(s,2t-1)=y(t)$ for $t \in [0,1]$. Define $n(x_1,x_2;y)$ as the number of elements of $\mathcal{M}(x_1,x_2;y,L,J)$ which have Fredholm index $0$. If the moduli spaces $\mathcal{M}(x_1,x_2;y,L,J)$ are transversely cut out, something that can be achieved by perturbing $H$ and $J$, the numbers $n(x_1,x_2;y)$ are always finite.

 Define $\Upsilon_L: \CW(H,L) \otimes \CW(H,L) \to \CW(\widehat{H},L) $ by
\begin{equation}
\Upsilon_L(x_1,x_2)= \sum_{y \in \mathcal{T}_L(\widehat{H})}(n(x_1,x_2;y)\mod 2) y
\end{equation}
for  $x_1,x_2\in \mathcal{T}_{L}(H)$, and extending it linearly to $ \CW(H,L) \otimes \CW(H,L)$.
It is proved in \cite{AS-product} that $\Upsilon_L$ descends to a map $H \Upsilon_L: \HW(H,L) \otimes \HW(H,L) \to \HW(\widehat{H},L)$, that endows $ \HW(H,L)$ with a product which we denote by $\star$. 
 It is compatible with the continuation maps, as follows by the results in \cite[Chapter 5]{Schwarz-thesis}, and passing to the direct limit  $H\Upsilon_L$  endows $ \HW(M,L)$ with a product. For homology classes $\mathfrak{h},\mathfrak{h}' \in \HW(M,L)  $ we will also denote their product by $\mathfrak{h} \star \mathfrak{h}'$. 
The product $\star$ is associative: the proof is identical to the proof in \cite{Schwarz-thesis} that the pair of pants product in Floer homology is associative.
As $\star$ is distributive with respect to the vector space structure of $ \HW(M,L)$ it gives $ \HW(M,L)$ the structure of a ring. Since we defined $\HW(M,L)$ with coefficients in $\mathbb{Z}_2$ the product $\star$ actually endows $\HW(M,L)$ with the structure of an algebra.

It was proved in \cite{AS-product} that in the case where $M =T^*Q$ of a compact manifold $Q$ and $L=T_q Q$ for some point $q\in Q$, the triangle product coincides with the Pontrjagin product.

\color{black}
An important property of the triangle product is given by 
\begin{lem} \label{lem:subad}
The spectral numbers $c$ of ${\HW}(M,L)$ are subadditive with respect to $\star$.
\end{lem}
 \textit{Proof:} We will need the triangle inequality 
\begin{equation}\label{triangle_inequ}
\mathcal{A}_{\widehat{H}}(y) \leq \mathcal{A}_{H}(x_1) + \mathcal{A}_{H}(x_2),
\end{equation}
that must be satisfied by the actions of $x_1,x_2 \in \mathcal{T}_{L}(H)$ and $y \in \mathcal{T}_{L}(\widehat{H})$ if the moduli space $\mathcal{M}(x_1,x_2;y,L,J) \neq \emptyset$ (see \cite[Formula 3.18]{AS-product}).

Let $\mathfrak{h}_1,\mathfrak{h}_2 \in {\HW}(M,L)$. Given $\delta>0$, we know from Lemma \ref{lem:minimax} that there exits Hamiltonians $H_1, H_2 \in \mathcal{H}_{\reg}(M,L)$ and cycles $w'_i \in \CW(H_i,L)$ such that
\begin{eqnarray*}
  \chi^{H_i}([w'_i]) = \mathfrak{h}_i \ \mbox{ and }  \ \mathcal{A}(w'_i)< c(\mathfrak{h}_i) +\frac{\delta}{2}
\end{eqnarray*}
  for  $i=1,2$.
 Let now $H\in \mathcal{H}_{\reg}(M,L)$ such that $H\geq H_1$ and $H\geq H_2$. We define $w_i:= \iota^{H_i,H}(w'_i)$  for  $i=1,2$. Since the action decreases under the continuation maps $\iota^{H_i,H}$  we have $\mathcal{A}(w_i)< c(\mathfrak{h}_i) +  \frac{\delta}{2}$, and  using \eqref{eq:lepolepo} we obtain
 \begin{equation*}
 \chi^H([w_i])=\chi^{H} ( \iota^{H_i,H}([w'_i])) = \chi^{H_i}([w'_i])=\mathfrak{h}_i, 
 \end{equation*}
 for $i=1,2$.
   By \eqref{triangle_inequ} we have $\mathcal{A}(\Upsilon_L(w_1 \otimes w_2))\leq c(\mathfrak{h}_1) +c(\mathfrak{h}_2) +\delta$. By definition $[\Upsilon_L(w_1 \otimes w_2)]= [w_1]\star[w_2]$, and
  by our construction of $\star$ in $\HW(M,L)$ we have 
  \begin{equation*}
  \chi^{\widehat{H}}([w_1]\star[w_2]) = \chi^H([w_1])\star  \chi^H([w_2]) = \mathfrak{h}_1 \star \mathfrak{h}_2.
  \end{equation*}
    It then follows from Lemma \ref{lem:minimax} that $c( \mathfrak{h}_1 \star \mathfrak{h}_2) \leq \mathcal{A}(\Upsilon_L(w_1 \otimes w_2)) \leq c(\mathfrak{h}_1) +c(\mathfrak{h}_2) +\delta$. 
  
 Summing up, we have shown that for any $\delta>0$ we have $c( \mathfrak{h}_1 \star \mathfrak{h}_2) <c(\mathfrak{h}_1) +c(\mathfrak{h}_2) +\delta$, which implies
  \begin{equation}
  c( \mathfrak{h}_1 \star \mathfrak{h}_2) \leq c(\mathfrak{h}_1) +c(\mathfrak{h}_2).
  \end{equation}
 \qed

\color{black}

We are ready to define the algebraic growth of $\HW$.
\begin{defn} \label{defi:growthHW}
Let $S$ be a finite set of elements of $\HW(M,L)$. We define 
\begin{equation}
\Gamma^{\alg}_S(M,L):=\Gamma^{\alg}_S(\HW(M,L)).
\end{equation}
\end{defn}

Combining Lemma \ref{finite_growth} and Lemma \ref{lem:subad} we obtain: 
\begin{prop} \label{prop:alg_symp}
For every finite set $S$ of  $\HW(M,L)$ we have
\begin{equation}
\Gamma^{\symp}(M,L) \geq \frac{\Gamma^{\alg}_S(M,L)}{\rho(S)}.
\end{equation}
\end{prop}

\color{black}

\subsubsection{$\HW(M,L \to L')$ as a module over $\HW(M,L)$}

We start by picking two exact asymptotically conical Lagrangians $L$ and $L'$ on $(M,\omega,\lambda)$. 
The boundary $\partial(\Delta)$ contains three connected components: the component $\mathcal{D}_{\lef
}$ which is equal to $\mathbb{R}\times \{-1\}$, the component $\mathcal{D}_{\midde}$ which contains the singular point, and the component $ \mathcal{D}_{\rig}$ which is equal to $\mathbb{R}\times \{1\}$. 

Let $x \in \mathcal{T}_{L}(H)$, $z  \in \mathcal{T}_{L \to L'}(H)$ and $\widetilde{z}  \in \mathcal{T}_{L \to L'}(H)$. 
We let $\mathcal{M}(x;z,\widetilde{z},J,H)$ be the moduli space of maps $u : \Delta \to \widehat{M}$ which satisfy \eqref{Floer} and such that $u(\mathcal{D}_{\lef}) \subset L$, $u(\mathcal{D}_{\midde}) \subset L$, $u(\mathcal{D}_{\rig}) \subset L'$, and $\lim_{s\to -\infty}u(s,t-1) = x(t)$ for $t\in [0,1]$, $\lim_{s\to -\infty}u(s,t) = z(t)$ for $t\in [0,1]$, and $\lim_{s\to +\infty}u(s,2t-1) = \widetilde{z}(t)$ for $t\in [0,1]$. 
Let $n(x;z,\widetilde{z})$ be the number of elements in $\mathcal{M}(x;z,\widetilde{z},J,H)$ that have Fredholm index $0$. For non-degenerate $H$ and a generic choice of $J$, all the spaces $\mathcal{M}(x;z,\widetilde{z},J,H)$ are transversely cut out, and therefore the numbers $n(x;z,\widetilde{z})$ are all finite. 

We then define a map $\Theta_{L,L'}: \CW(H;L) \otimes \CW(H; L\to L') \to \CW(\widehat{H} ; L \to L')  $ by letting
\begin{equation}
\Theta_{L,L'}(x\otimes z) = \sum_{\widetilde{z}  \in \mathcal{T}_{L \to L'}(H)} (n(x;z,\widetilde{z})\mod 2 )\widetilde{z},
\end{equation}
for $x \in \mathcal{T}_{L}(H)$, $z  \in \mathcal{T}_{L \to L'}(H)$ , and extending it linearly to $\CW(H;L) \otimes \CW(H; L\to L')$.

 The map $\Theta_{L,L'}$ descends to a map $H\Theta_{L,L'}: \HW(H;L) \otimes \HW(H; L\to L') \to \HW(\widehat{H} ; L \to L')  $.
  The proof is again identical to the one used in \cite{Schwarz-thesis} to show that the pair of pants product descends to the Floer homology.
   Taking direct limits we obtain a product $H\Theta_{L,L'}: \HW(M,L) \otimes \HW(M,L \to L') \to \HW(M,L \to L')$. We will use the notation $H\Theta_{L,L'}(\mathfrak{h},\mathfrak{m}) = \mathfrak{h} \ast \mathfrak{m}$.
  
  In order to conclude that $\HW(M; L\to L')$ is a module over the algebra $\HW(M; L)$ we must prove that:
 \begin{enumerate}
\item[] $ \mathfrak{h}  \ast (\mathfrak{m}_1 + \mathfrak{m}_2) = \mathfrak{h}  \ast \mathfrak{m}_1  + \mathfrak{h}  \ast \mathfrak{m}_2$  for all $\mathfrak{h}  \in \HW(H;L)$ and $ \mathfrak{m}_1, \mathfrak{m}_2 \in \HW(H ; L \to L') $,
\item[] $ (\mathfrak{h}_1 + \mathfrak{h}_2)  \ast  \mathfrak{m} = \mathfrak{h}_1 \ast \mathfrak{m}  +\mathfrak{h}_2 \ast \mathfrak{m}$  for all $\mathfrak{h}_1, \mathfrak{h}_2 \in \HW(H;L)$ and $\mathfrak{m} \in \HW(H ; L \to L') $,
\item[] $( \mathfrak{h}_1 \star \mathfrak{h}_2)  \ast \mathfrak{m}  = \mathfrak{h}_1 \ast (\mathfrak{h}_1  \ast  \mathfrak{m} )$ for all $\mathfrak{h}_1, \mathfrak{h}_2 \in \HW(H;L)$ and $\mathfrak{m} \in \HW(H ; L \to L') $.
\end{enumerate}
The first two properties follow from the linearity of $H\Theta_{L,L'}$. The proof of the third one is a cobordism argument identical to the of \cite[Chapter 5]{Schwarz-thesis} that proves the associativity the associativity of the triangle product $\star$.  
\color{black}
An argument identical to one used to prove Lemma \ref{lem:subad} gives
\begin{lem} \label{lem:subad2}
The spectral numbers $c$ are subadditive with respect to $\ast$. 
\end{lem}

\color{black}
\section{Viterbo functoriality}\label{sec:Viterbo}

The Viterbo transfer map on $\HW$ will be described.  
As first applications we then deduce invariance properties under a graphical change of the boundary of the Liouville domain in the completion.   
\subsection{The Viterbo transfer map on $\HW$}

The Viterbo transfer map was first introduced as a map for symplectic homology in \cite{Viterbo1999}, see also \cite{Cieliebak2001, McLean2009}. The analogous map in wrapped Floer homology was studied by \cite{AbouzaidSeidel2010}, see also \cite{Ritter2013}. Our focus lies on its compatibility with the action filtration.  

Let $M:=(Y_M, \omega_M, \lambda_M)$ be a Liouville domain and let $j: W \rightarrow M$ be a codimension $0$ exact embedding of a Liouville domain $W:=(Y_W, \omega_W, \lambda_W)$ into $M$, i.e. $j^{*}\lambda_M = \lambda_W$. 
Let $L_0$ and $L_1$ be asymptotically conical exact Lagrangians in $M$, and assume $L^{'}_0 := L_0 \cap W$ and $L^{'}_1 := L_1 \cap W$ are asymptotically conical in $W$.
Additionally assume that $L_0$ is also conical on $M \setminus W$ and $L_1$ satisfies the property
\begin{align}
\begin{split}\label{transfer_admissible}
 &\lambda|_{L \setminus L^{'}} \text{ vanishes on the boundary }\partial (L\setminus L^{'}) = \partial L \cup \partial L^{'}, \text{ and}\\
&\text{one can write }\lambda|_{L \setminus L^{'}} = df, \text{ where }f \text{ vanishes near } \partial L \cup \partial L^{'}. 
\end{split}
\end{align}
We will call a Lagrangian with this property \textit{transfer admissible} for the pair $(M,W)$.
See \cite{AbouzaidSeidel2010} for a discussion of that condition and why the transfer map can in general not be defined if one removes this condition.

We give the construction of the \textit{Viterbo transfer map} as an asymptotical morphism of filtered directed systems ${j}(L_0,L_1): \widetilde{\HW}(M,L_0 \to L_1) \to \widetilde{\HW}(W,L^{'}_0 \to L^{'}_1)$.  More precisely we get for  $a > K = \mathsf{K}(M,L_0 \to L_1)$, defined in \eqref{K(M,L)},  homomorphisms 
\[
j_{!}(L_0,L_1)_a : \HW^a(M,L_0 \to L_1) \rightarrow \HW^a(W, L^{'}_{0} \to L^{'}_{1})
	\]
that are compatible with the persistence morphisms $\iota_{a\rightarrow b}$, for $K< a < b$. Moreover, the homomorphisms are functorial with respect to a composition of embeddings $W_1 \subset W_2 \subset M$ and the induced maps in the direct limit 
\[
\bar{j}_{!}(L_0) = \bar{j}_{!}(L_0,L_0): \HW(M,L_0) \rightarrow \HW(W,L^{'}_0), \text{ and} 
\] 
\[
\bar{j}_{!}(L_0,L_1): \HW(M,L_0 \to L_1) \rightarrow \HW(W, L^{'}_0 \to L^{'}_1)
\] 
are compatible with the algebra and module structure, i.e.
\begin{equation}\label{vitproduct}
\bar{j}_{!}(L_0) (x \star y) =  \bar{j}_{!}(L_0)(x) \star \bar{j}_{!}(L_0)(y)
\end{equation}
and
\begin{equation}\label{vitmodule}
 \bar{j}_{!}(L_0,L_1)(x \ast z) = \bar{j}_{!}(L_0)(x) \ast \bar{j}_{!}(L_0,L_1)(z)
\end{equation}
for all $x, y \in \HW(M,L_0)$ and $z \in \HW(M,L_0, L_1)$.

We first give the definition of ${j}_{!}(L_0,L_1)$. We may assume that $(M,L_0 \to L_1)$ and $(W,L^{'}_0 \to L^{'}_1)$ are regular. \color{black} Otherwise we can perform the construction considering suitable compactly supported  Hamiltonian perturbations of $L_0$ and $L_1$.  \color{black} Let $\mathcal{S} := \mathcal{S}(M,L_0 \to L_1) \cup \mathcal{S}(W,L^{'}_0 \to L^{'}_1)$.  We furthermore assume that actually $W \subset M_{\tau^2}$ for some $\tau < 1$, sufficiently close to $1$. One can get the maps for general $W \subset M$ by an inverse limit. 

First of all, for every $R > 1$ one can construct a compactly supported Hamiltonian isotopy $(\psi^R_t)_{t\in [0,1]}$ on $\widehat{M}$, ($\psi^R_0 = \id$, $\psi:= \psi^R_1$) that leaves $\widehat{L}_0$ invariant  and maps $\widehat{L}_1$ to a Lagrangian $\widehat{L}^R_1$ that is conical on $(\widehat{M} \setminus M_R) \cup (W_R \setminus W)$ and that is transfer admissible for the pair $(M_R, W_R)$ as follows. Map $L_1 \setminus W$  by the Liouville flow $(\phi_{\log t})_{t\in [1,R]}$ into $A_R = M_R \setminus W_R$. Since  $L_1$ is conical near $\partial W$, we can extend $\left(L^{'}_1 \cup \phi_{\log{t}}(L_1 \setminus W)\right)_{t\in [1,R]}$ to a 1-parameter family of exact Lagrangians interpolating between $\widehat{L}_1$ and a Lagrangian $\widehat{L}^R_1$. Therefore we can choose a Hamiltonian isotopy $(\psi^R_t)_{t \in [0,1]}$ in $\widehat{M}$ that realizes this Lagrangian isotopy and is supported in $M_{\frac{1}{\tau}R} \setminus W_{\tau}$. Since $\widehat{L}_0$ is conical outside $W$, we can choose the isotopy to leave $\widehat{L}_0$ invariant. We can choose the isotopy such that $(\psi\circ\zeta)^{*}\lambda = R\zeta^{*}\lambda$,
where $\zeta: {L}_1 \setminus W \hookrightarrow \widehat{M}$ is the embedding of $L_1$ restricted to $L_1 \setminus W$. The function
\[
f_R: \widehat{L}^R_1 \to \R, \text{ with } f_R(x) = \begin{cases} f_1(x), &\text {if } x\in L_1 = \widehat{L}^R_1\cap W, \\
       Rf_1(\psi^{-1}x), &\text{ elsewhere}
			\end{cases}
\]
is a primitive of $\lambda|_{\widehat{L}^R_1}$.
 

We now carefully choose for every $\mu \notin \mathcal{S}$ sufficiently large a step-shaped Hamiltonian $H^{step}_{\mu}$ on $\widehat{M}$. 
Let $k_W := \min\{f_0(x) - f_1(x) \, | \, x \in L_0 \cap L_1 \cap W \}$ where $f_i$ are the primitives of $\lambda|_{L_i}$, $i=0,1$. Let $\widetilde{k} = \max\{-k_W, 0\}$. Let 
\[
\widetilde{K} = \mathsf{K}(M,W,L_0 \to L_1) = \max\{ \max \{f_0(x) - f_1(x) \, | \, x \in L_0 \cap L_1 \cap M\setminus W\}, 0\}.
\]
Choose a small $\epsilon > 0$. 
Let $\mu > \widetilde{K}$, $\mu \notin \mathcal{S}$, and let $\delta_{\mu} = \min \{\dist(\mu,\mathcal{S}), \mu - \widetilde{K} \}$. Choose $R > \frac{\widetilde{k} + \mu + 4\epsilon}{\delta_{\mu}}$.

\input{figure1.tex}

We choose a smooth function $H^{step}_{\mu}: \widehat{M} \rightarrow \R$ that only depends on the radial coordinate $r= r_W$ in $(0,R) \times \partial W$ and only on the radial coordinate $r=r_M$ in $(\tau R,\infty) \times \partial M$, and such that
\begin{equation}\label{step}
H^{step}_{\mu}(x) = \begin{cases}
-\epsilon,&\text{ if }x \in W_{\tau} \\
\frac{\partial^2 H}{\partial r} \geq 0, &\text{ if } x = (r,y) \in W \setminus W_{\tau} \\
\mu r - \mu, &\text{ if } x=(r,y) \in W_{\tau R} \setminus W \\
\frac{\partial^2 H}{\partial r} \leq 0,  &\text{ if } x = (r,y) \in W_R \setminus W_{\tau R} \\
(R-1)\mu - \epsilon, &\text{ if } x \in M_{\tau R} \setminus W_R \\
\frac{\partial^2 H}{\partial r} \geq 0,  &\text{ if } x = (r,y) \in M_R \setminus M_{\tau R} \\
\mu r -\mu, &\text{ if }x = (r,y) \in \widehat{M} \setminus M_R.
\end{cases}
\end{equation}

We divide the critical points of the action functional $\mathcal{A} := \mathcal{A}_{H_{\mu}^{step}}^{\widehat{L}_0 \to \widehat{L}^R_1}$ of $H_{\mu}^{step}$ with respect to $\widehat{L}_0$ and $\widehat{L}^R_0$ into four classes: Intersections of $L_0$ and $L_1$ in $W_{\tau}$ denoted by $\mathfrak{A}^{*}$, Hamiltonian chords close to $\partial W$ denoted by $\mathfrak{A}^{**}$, intersections of $\widehat{L}_0$ and $\widehat{L}^R_1$ in $M_R \setminus W_R$ denoted by $\mathfrak{B}^{*}$, and chords close to $\partial W_R$ and $\partial M_R$ denoted by $\mathfrak{B}^{**}$. 
We can estimate the action values as follows. 
\begin{align}
\mathcal{A}(x) &\geq k_w - \epsilon \geq -\widetilde{k} - \epsilon, \text{ if } x \in \mathfrak{A}^{*},  \\
\mathcal{A}(x) &> -\epsilon \geq -\widetilde{k} - \epsilon, \text{ if } x  \in \mathfrak{A}^{**}, \\
\mathcal{A}(x) &\leq R\widetilde{K} - ((R-1)\mu -\epsilon) < - \widetilde{k} - 3\epsilon,  \text{ if } x \in \mathfrak{B}^{*}, \text{ and} \label{B*}\\
\mathcal{A}(x) &< (\mu - \dist(\mu, \mathcal{S}))R - ((R-1)\mu -\epsilon) < -\widetilde{k} - 3\epsilon, \text{ if } x \in \mathfrak{B}^{**}.
\end{align}

In \eqref{B*} we use that $f_0(x) - f_R(x) \leq \widetilde{K}R$ for every $x \in \mathfrak{B}^{*}$. 

Altogether we get that $\mathcal{A}(x) \geq -\widetilde{k} - \epsilon $, if $x\in \mathfrak{A} = \mathfrak{A}^{*} \cup \mathfrak{A}^{**}$ and $\mathcal{A}(x) < -\widetilde{k} - 3\epsilon$, if $x\in \mathfrak{B} = \mathfrak{B}^{*} \cup \mathfrak{B}^{**}$. Hence there are no Floer trajectories from $\mathfrak{B}$ to $\mathfrak{A}$. So
 $\CW^{(-\widetilde{k} - 2\epsilon, +\infty)}_{*}(H^{step}_{\mu}; \widehat{L}_0 \to \widehat{L}^R_1) = \CW_*(H^{step}_{\mu})/\CW^{(-\infty, -\widetilde{k} - 2\epsilon)}_{*}(H^{step}_{\mu})$ generated by elements of action larger then $-\widetilde{k} - 2\epsilon$ is a chain complex, and the projection $\CW(H^{step}_{\mu}) \rightarrow \CW^{(-\widetilde{k} - 2\epsilon, +\infty)}(H^{step}_{\mu})$ induces a map 
\begin{equation}\label{vit}
\HW(H^{step}_{\mu}; \widehat{L}_0 \to \widehat{L}^R_1) \rightarrow \HW^{(-\widetilde{k} - 2\epsilon,+\infty)}(H^{step}_{\mu}; \widehat{L}_0 \to \widehat{L}^R_1)
\end{equation}
on homology. 

Let now $H^M_{\mu}$ be a non-degenerate admissible Hamiltonian with respect to $M$  on $\widehat{M}$ with slope $\mu$, and $H^W_{\mu}$ a non-degenerate admissible Hamiltonian with respect to $W$ on $\widehat{W}$ with slope $\mu$. We have the isomorphisms 
\begin{align}
&\HW(H_{\mu}^M; L_0 \to L_1) \overset{\cong}\rightarrow \HW((\psi^{-1})^{*}H_{\mu}^M; \widehat{L}_0 \rightarrow \widehat{L}^R_1) \overset{\cong}\rightarrow \HW(H^{step}_{\mu}; \widehat{L}_0 \to \widehat{L}^R_1), \text{ and} \label{H^M}\\
&\HW^{(-\widetilde{k} - 2\epsilon,+\infty)}(H^{step}_{\mu}; \widehat{L}_0 \to \widehat{L}^R_1) \overset{\cong}\rightarrow \HW(H_{\mu}^W; L^{'}_0 \to L^{'}_1). \label{H^W} 
\end{align}

Here, the second isomorphism in \eqref{H^M} holds, since $(\psi^{-1})^{*}H_{\mu}^M$ and $H^{step}_{\mu}$ can be connected by a compactly supported homotopy of Hamiltonians. To get the isomorphism in $\eqref{H^W}$ we choose a conical almost complex structure near $\partial W$. By \cite[Lemma 7.2]{AbouzaidSeidel2010}, see also \cite[Appendix D]{Ritter2013} there are no Floer trajectories with asymptotics in $W$ that leave $W$ and hence the differential of $\CW^{(-\widetilde{k} - 2\epsilon,+\infty)}(H^{step}_{\mu},  \widehat{L}_0 \to \widehat{L}^R_1)$ only counts Floer trajectories that map into $W$.

Combining \eqref{vit}, \eqref{H^M}, and \eqref{H^W} gives maps
\begin{equation}\label{vit1}
 j_{\mu}: \HW(H^M_{\mu}; L_0 \to L_1) \rightarrow \HW(H^W_{\mu}; L^{'}_0 \to L^{'}_1) 
\end{equation}
for any  $\mu > \widetilde{K}$, $\mu \notin \mathcal{S}$. 
The isomorphisms  \eqref{vit}, \eqref{H^M}, and \eqref{H^W} are all compatible with Floer continuation maps induced by  monotone increasing homotopies of the corresponding Hamiltonians. We do not give the details here and refer the reader to \cite[Theorem 9.8]{Ritter2013}. We thus get commutative diagrams
\[
\begin{CD}
 \HW(H^M_{\mu}; L_0 \to L_1)  @>{j_{\mu}}>>  \HW(H^W_{\mu}; L^{'}_0 \to L^{'}_1) \\
 @V{\iota^{H^M_{\mu}, H^M_{\eta}}}VV          @V{\iota^{H^W_{\mu}, H^W_{\eta}}}VV \\
 \HW(H^M_{\eta}; L_0 \to L_1)  @>{j_{\eta}}>>  \HW(H^W_{\eta}; L^{'}_0 \to L^{'}_1)
\end{CD}
\]
for any $\eta > \mu > \widetilde{K}$, $\mu, \eta \notin \mathcal{S}$. 

Hence, for any $a > K=\mathsf{K}(M,L_0 \to L_1) \geq \widetilde{K}$ one obtains, because of the construction in Section \ref{rem:rem1}, a map  
\[
j_{!}(L_0,L_1)_a : \HW^a(M, L_0 \to L_1) \rightarrow \HW^a(W, L^{'}_0 \to L^{'}_1) 
\]
induced in the direct limit taken over all non-degenerate admissible Hamiltonians with slope $\mu$, $K <\mu < a$. By the construction these maps are compatible with the persistence morphisms $\iota_{a \rightarrow b}$, for $K < a < b$.

By a standard compactness-cobordism argument, and by using once again the non-escaping result \cite[Lemma 7.2]{AbouzaidSeidel2010} one can show the compatibility of the algebra and module structure with the Viterbo transfer maps \eqref{vitproduct} and \eqref{vitmodule}; for this see \cite{Ritter2013}.

\subsection{Change of the contact hypersurface $\partial M$}

From the Viterbo transfer one can deduce invariance properties of $\HW$ under a graphical change of $\partial M$ in $\widehat{M}$. This will be used to bound the growth rate of Reeb chords for different choices of contact forms on $(\partial M, \xi_{M})$.   
Let $M$ be a Liouville domain with asymptotically conical exact Lagrangians $L_0$ and $L_1$ as above, let $0 < \epsilon < 1$.
\begin{lem}\label{M_epsilon}
Assume that  $L_i$, $i=0,1$, are conical on $M\setminus M_{\epsilon}$. Then, for  $a > K = \mathsf{K}(M,L_0 \to L_1)$, we have $\HW^a(M_{\epsilon},L_0 \cap M_{\epsilon} \to L_1 \cap M_{\epsilon}) \overset{\varphi_a}{\cong} \HW^{\frac{1}{\epsilon}a}(M,L_0 \to L_1)$.
Moreover, the Viterbo map $\HW^a(M, L_0 \to L_1) \rightarrow \HW^a(M_{\epsilon},L_0 \cap M_{\epsilon} \to L_1 \cap M_{\epsilon})$ composed with $\varphi_a$ is the persistence morphism $\HW^a(M, L_0 \to L_1) \rightarrow \HW^{\frac{1}{\epsilon}a}(M,L_0 \to L_1)$. 
\end{lem}
\textit{Proof: }
Note, that adding a constant to any Hamiltonian $H$ or applying a compactly supported deformation to $H$ does not change its Floer homology. 
Let $H$ be an admissible Hamiltonian with slope $\mu$ with respect to $M_{\epsilon}$. Then $H-\mu(\frac{1}{\epsilon}-1)$ is an admissible Hamiltonian with slope $\frac{1}{\epsilon} \mu$ with respect to $M$. Moreover, if one chooses a cofinal sequence of Hamiltonians of the first kind with slopes $K< \mu < a$, there are compactly supported homotopies of the shifted Hamiltonians to a cofinal sequence with respect to $M$ with slopes $\frac{1}{\epsilon}\mu$. This gives the first statement.  

Observe, that both the Viterbo transfer map in the present situation and the persistence morphisms are given by a continuation map induced by a monotone homotopy. One can apply a usual chain homotopy argument in Floer homology to see the second statement.
\qed

Let $\mathsf{f}: \partial M \rightarrow [1, \infty)$ be a smooth function. Recall that $M_{\mathsf{f}} = \widehat{M} \setminus \{(r,x) \, | \, r > \mathsf{f}(x), x \in \partial M\}$. Let $\zeta = \max_{\partial M} \mathsf{f}$. 
\begin{lem}\label{lem:changeofhypersurface}
The filtered directed systems $(\HW^a(M,L_0 \to L_1))_{a\in (0,\infty)}$ and  $(\HW^a(M_{\mathsf{f}},\widehat{L}_0 \cap M_{\mathsf{f}} \to \widehat{L}_0 \cap M_{\mathsf{f}} ))_{a\in (0, \infty)}$ are $(\zeta, 1)$-interleaved. 
\end{lem}
\textit{Proof: }
The morphisms of filtered directed systems $f$ and $g$, with
$f_a: \HW^a(M,L_0 \to L_1) \cong \HW^{\zeta a}(M_{\zeta}, L_0 \to L_1) \to \HW^ {\zeta a}(M_{\mathsf{f}},\widehat{L}_0 \cap M_{\mathsf{f}}, \widehat{L}_1 \cap M_{\mathsf{f}})$ and $g_a:\HW^a(M_{\mathsf{f}},\widehat{L}_0 \cap M_{\mathsf{f}}, \widehat{L}_1 \cap M_{\mathsf{f}}) \to \HW^a(M,L_0 \to L_1)$, given by Viterbo maps, yield by functoriality of Viterbo maps and Lemma \ref{M_epsilon} the $(\zeta,1)$- interleaving of $(\HW^a(M,L_0 \to L_1))_{a\in (0,\infty)}$ and  $(\HW^a(M_{\mathsf{f}},\widehat{L}_0 \cap M_{\mathsf{f}} \to \widehat{L}_0 \cap M_{\mathsf{f}} ))_{a\in (0, \infty)}$.
\qed

\section{From algebraic growth to positivity of topological entropy}\label{sec:entropy}
In this section we prove Theorem \ref{theorementropy}.

\color{black}
\subsection{Legendrian isotopies, transfer admissible Lagrangians and growth} \label{subsec:technical}

We start by introducing some notation. Let $M=(Y,\omega,\lambda)$ be a Liouville domain and $L$ be an asymptotically conical exact Lagrangian disk in $M$. We denote by $\Lambda$ the Legendrian sphere $\partial L$.
 Letting $\Sigma := \partial M$ and  $\alpha_{M} := \lambda |_{\Sigma}$ be the contact form induced by $M$ on $\Sigma$ we assume that $(\alpha_{M}, \Lambda \to \Lambda)$ is regular. As usually, we denote by $\xi_{M}$ the contact structure $\ker \alpha_{M}$.
 
 

Our approach to prove invariance of the exponential symplectic growth of $\HW$ differs from the ones developed by \cite{MacariniSchlenk2011,McLean2012}. It makes extensive use of the module and algebra structures that exist on $\HW$.
 We will need the following
 \begin{defn}
Let $\mu>0$ and $\Lambda_0$ be a Legendrian sphere in $(\Sigma, \xi_{M})$. Assume that $\Lambda_1$ is Legendrian isotopic to $\Lambda_0$. We say that $\Lambda_1$ is \textit{$\mu$-close to $\Lambda_0$ in the $C^3$-sense} if there exists a Legendrian isotopy $\theta:[-1,1] \times S^{n-1} \to (\Sigma, \xi_{M})$ from $\Lambda_0$ to $\Lambda_1$ whose $C^3$-norm is $<\mu$, and which is stationary in the first coordinate outside a compact subset of $(-1,1)$.
\end{defn}

Recall that the symplectisation of a contact form $\alpha$ on $(\Sigma, \xi_{M})$ is the exact symplectic manifold $((0,+\infty)\times \Sigma, dr\alpha, r\alpha)$ where $r$ denotes the first coordinate in $(0,+\infty)\times \Sigma$.
The following lemma is essentially due to Chantraine \cite{Baptiste} and is proved in Appendix B.
\begin{lem} \label{lemmaBaptiste}
Fix a constant $\epsilon >0$, a contact form $\alpha$ on $(\Sigma,\xi)$, a Legendrian $\Lambda_0$ in $(\Sigma,\xi)$, and a tubular neighbourhood $U(\Lambda_0)$ of $\Lambda_0$ in $\Sigma$. Then there exists $\delta>0$ such that if $\Lambda_1$ is $\delta$-close to $\Lambda_0$ in the $C^3$-sense, then there exist exact Lagrangian cobordisms $\mathcal{L}^-$ from $\Lambda_1$ to $\Lambda_0$ and $\mathcal{L}^+$ from $\Lambda_0$ to $\Lambda_1$ in the symplectization of $\alpha$ satisfying:
\begin{itemize}
\item [a)] $\mathcal{L}^- $ is conical outside  $[1-\frac{\epsilon}{2}, 1-\frac{\epsilon}{4}] \times \Sigma$,
\item [b)] $\mathcal{L}^+$  is conical outside  $[1+\frac{\epsilon}{4}, 1+ \frac{\epsilon}{2}] \times \Sigma$,
\item [c)] the projections of $\mathcal{L}^+$ and $\mathcal{L}^-$ to $\Sigma$ are completely contained in $U(\Lambda_0)$,
\item [d)] the primitives $f^\pm$ of $(r\alpha) |_{\mathcal{L^\pm}}$ have support in $[1-\frac{\epsilon}{2}, 1-\frac{\epsilon}{4}] \times \Sigma$ and $[1+\frac{\epsilon}{4}, 1+ \frac{\epsilon}{2}] \times \Sigma$, respectively, and $| f^\pm |_{C^0} < \epsilon$.
\end{itemize}
Moreover if $\mathcal{L}$ is the exact Lagrangian cylinder obtained by gluing  $\mathcal{L}^+ \cap [1,+\infty) \times \Sigma)$ on top of $\mathcal{L}^- \cap ((0,1] \times \Sigma)$ we have that
\begin{itemize}
\item [e)] $\mathcal{L}$ is Hamiltonian isotopic to $\mathbb{R}\times \Lambda_0$ in the symplectization of $\alpha$, and the Hamiltonian producing the isotopy can be taken to have  support  in $[1- \frac{\epsilon}{2},1+  \frac{\epsilon}{2}] \times \Sigma$.
\end{itemize}
\end{lem}
 
 We now fix $\epsilon > 0$ such that $L$ is conical on $M \setminus M_{1 - 2\epsilon}$. We choose a Legendrian tubular neighbourhood $\mathcal{U}(\Lambda)$ of $\Lambda$ on $(\Sigma,\xi_M) $. For these choices of $\epsilon>0$ and $\mathcal{U}(\Lambda)$, we choose $\delta_1>0$ given by  Lemma \ref{lemmaBaptiste}.
 
 We then choose a Legendrian sphere $\Lambda_1$ which is $\delta_1$-close to $\Lambda$ in the $C^3$ sense, is disjoint from $\Lambda$, and satisfies that $(\alpha_{M}, \Lambda \to \Lambda_1)$ is regular.

It follows from Lemma \ref{lemmaBaptiste} that there exists an  exact Lagrangian cobordism $\mathcal{L}^-$  from $\Lambda_1$ to $\Lambda$ in the symplectization of $\alpha_{M}$ which is conical
outside $[1-\frac{\epsilon}{2},1-\frac{\epsilon}{4}] \times \Sigma$. We can then glue $\mathcal{L}^- \cap[1-\frac{\epsilon}{2},1] \times \Sigma$ to $L \cap M_{1-\frac{\epsilon}{2}}$ to obtain an exact Lagrangian submanifold $L_1$ in $M$. The Lagrangian $L_1$ is an exact filling of $\Lambda_1$. Let $f_L$ be the primitive of $\lambda\mid_{L}$ which vanishes in $\Lambda$. Using Lemma \ref{lemmaBaptiste} we can glue $f^-$ to the restriction of $f_L$ to  $L \cap M_{1-\frac{\epsilon}{2}}$ to obtain primitive of $f_{L_1}$ of $\lambda\mid_{L_1}$ which vanishes in $\Lambda_1$. 

Because of the control given by Lemma \ref{lemmaBaptiste} on the function $| f^- |_{C^0}$ on $\mathcal{L}^-$, and the facts that $L$ and $L_1$ coincide on $M_{1-\frac{\epsilon}{2}}$ and $f_L$ vanishes on $L  \cap(M\setminus M_{1-\frac{\epsilon}{2}})$ we have
\begin{equation} \label{eq:control:LtoL_1}
\mathsf{K}(M,L \to L_1) < \epsilon.
\end{equation}

By Lemma \ref{lemmaBaptiste} d) the Lagrangian $L_1$ is transfer admissible for the pair $(M,M_{1-\epsilon})$. Combining  this with \eqref{eq:control:LtoL_1} we obtain for each $a>\epsilon\geq \mathsf{K}(M,L \to L_1)$ a Viterbo map $\Psi^a_{\mathcal{L}^-} : \HW^a(M,L \to L_1) \to \HW^a(M_{1-\epsilon},L)$, where to simplify notation we keep denoting by $L$ and $L_1$ the restrictions of $L$ and $L_1$ to  $M_{1-\epsilon}$. Passing to the direct limit we obtain a map $\Psi_{\mathcal{L}^-} : \HW(M,L \to L_1) \to \HW(M_{1-\epsilon},L)$.

By  Lemma \ref{lemmaBaptiste} we also have an exact Lagrangian cobordism $\mathcal{L}^+$ from $\Lambda$ to $\Lambda_1$, which is diffeomorphic to $\mathbb{R}\times S^{n-1}$, and is conical over $\Lambda$ for $r\geq 1 + \frac{\epsilon}{2}$ and conical over $\Lambda_1$ for $r \leq 1 + \frac{\epsilon}{4} $. By gluing $\mathcal{L}^+ \cap ([1, 1+ \epsilon] \times \Sigma)$ to $L_1$ we obtain an exact Lagrangian $\overline{L}$ in $M_{1+\epsilon}$.  By Lemma \ref{lemmaBaptiste} d) the Lagrangian $\overline{L}$ is transfer admissible for the pair $(M_{1+\epsilon},M)$. By gluing $f^+$ to $f_{L_1}$ we obtain a primitive $f_{\overline{L}}$ of $\lambda\mid_{\overline{L}}$. Reasoning as in the proof of \eqref{eq:control:LtoL_1}  one obtains
\begin{equation} \label{eq:control:LtoL_1above}
\mathsf{K}(M_{1+\epsilon},L \to \overline{L}) < \epsilon.
\end{equation}
 We thus obtain for each $a>\epsilon$ a Viterbo map $\Psi^a_{\mathcal{L}^+} : \HW^a(M_{1+\epsilon},L \to \overline{L}) \to \HW^a(M,L \to  L_1 )$, where by abuse of notation we denote by $L$ the conical extension of $L$ to $M_{1+\epsilon}$. Passing to the direct limit we obtain a map $\Psi_{\mathcal{L}^+} : \HW(M_{1+\epsilon},L \to \overline{L}) \to \HW(M,L \to  L_1 )$.

By Lemma \ref{lemmaBaptiste}, $\overline{L}$ is Hamiltonian isotopic to the conical extension of $L$ to $M_{1+\epsilon}$, which we will still denote by $L$, for a Hamiltonian function which vanishes outside $M_{1+\frac{\epsilon}{2}} \setminus  M_{1- \frac{\epsilon}{2}}$. 
A continuation argument then implies that for each admissible Hamiltonian $H$ that is regular for both $(M_{1+\epsilon},L \to \overline{L})$ and $(M_{1+\epsilon}, L)$ and has slope $>\epsilon$ we have that $\HW(H,L \to \overline{L})$ and $\HW(H, L)$ are isomorphic. By Section \ref{rem:rem1} we conclude that for each $a>\epsilon$ 
 the wrapped Floer homologies 
\begin{equation} \label{eq:isocont}
\HW^a(M_{1+\epsilon},L \to \overline{L}) \mbox{ and } \HW^a(M_{1+\epsilon}, L) \mbox{ are isomorphic}. 
\end{equation}
This induces an isomorphism $\Phi: \HW(M_{1+\epsilon}, L) \to \HW(M_{1+\epsilon},L \to  \overline{L})$.


Since $L$ is conical on  $M_{1+\epsilon} \setminus M_{1-\epsilon}$, $M \setminus M_{1-\epsilon}$ and  $M_{1+\epsilon} \setminus M$, we have transfer maps
\begin{itemize}
\item[] $\Psi^{\pm}_{L} : \HW(M_{1+\epsilon}, L)  \to \HW(M_{1-\epsilon}, L)$, 
\item[] $\Psi^{-}_{L} : \HW(M, L)  \to \HW(M_{1-\epsilon}, L)$,
\item[] $\Psi^{+}_{L} : \HW(M_{1+\epsilon}, L)  \to \HW(M, L)$.
 \end{itemize}
We notice that the contact forms induced by $\lambda$ on $\{1-\epsilon\} \times \Sigma$ and $\{1+\epsilon\} \times \Sigma$ are 
$\frac{\alpha_{M}}{1-\epsilon}$ and $\frac{\alpha_{M}}{1+\epsilon}$, respectively. Thus, as explained in Lemma \ref{M_epsilon}, the maps $\Psi^\pm_L$, $\Psi^-_L$ and $\Psi^+_L$
 are induced by asymptotic isomorphisms of f.d.s. For this reason we will denote by $A_L$ the algebras $ \HW(M_{1+\epsilon}, L)$, $\HW(M_{1-\epsilon}, L)$ and $ \HW(M, L)$. More generally, the same reasoning shows that for any $\zeta > -\epsilon$ the algebra $\HW(M_{1+\zeta}, L)$ is isomorphic to $ \HW(M_{1-\epsilon}, L)$ by an asymptotic isomorphism. 

The homologies $\HW(M, L \to L_1)$, $\HW(M_{1+\epsilon}, L \to \overline{L})$ and $\HW(M_{1-\epsilon}, L )$ are modules over the algebras $\HW(M, L)$, $\HW(M_{1+\epsilon}, L)$  and  $\HW(M_{1-\epsilon}, L)$, respectively: they are therefore $A_L$-modules.
By this discussion and \eqref{vitmodule} in section \ref{sec:Viterbo} the maps $\Phi$, $\Psi_{\mathcal{L}^-}$ and $\Psi_{\mathcal{L}^+}$ are $A_L$-module homomorphisms.

By functoriality of continuation maps, the diagram 
\\
\[
\begin{CD}
 \HW(M_{1+\epsilon}, L\to \overline{L})  @<\Phi<<  \HW(M_{1+\epsilon},L ) \\
@V{\Psi_{\mathcal{L}^-}\circ\Psi_{\mathcal{L}^+}}VV                                    @V\Psi^{\pm}_{L}VV  \\
 \HW(M_{1-\epsilon}, L ) @<\mathrm{id}<<\HW(M_{1-\epsilon}, L) \\
\end{CD}
\]
\\
is commutative. It thus follows that the map ${\Psi_{\mathcal{L}^-}\circ\Psi_{\mathcal{L}^+}}$ is an $A_L$-module isomorphism. We thus conclude that $\Psi_{\mathcal{L}^+}$ is injective. 
Let $\mathbf{1}_L$ be the unit in  $\HW(M_{1+\epsilon},L )$. As $\Phi$ is an $A_L$-module isomorphism and  $\Psi_{\mathcal{L}^+}$ is an injective $A_L$-module homomorphism we know that the element $m_{L_1} :=\Psi_{\mathcal{L}^+} \circ \Phi( \mathbf{1}_L)$ in  $\HW(M,L \to L_1 )$ is a stretching element.
We have thus proved the following:
\begin{lem} \label{lemmaprelim} 
The wrapped Floer homology $\HW(M,L \to L_1 )$ is a stretched module over $\HW(M,L)$. It follows from  Lemma \ref{mod_growth}, Lemma \ref{lem:subad}, and Lemma \ref{lem:subad2} that 
\begin{equation} 
\Gamma^{\symp}(M,L \to  L_1 ) \geq {\Gamma^{\symp}(M,L )}. 
\end{equation} 
\end{lem}

Recall that our Legendrian sphere $\Lambda_1$ was chosen disjoint from $\Lambda$. It follows that intersections of the Lagrangian disk $L_1$ and $L$ are contained in $M_{1-\frac{\epsilon}{5}}$.
By a small Hamiltonian isotopy supported inside $M_{1-\frac{\epsilon}{5}}$  we can perturb  $L_1$ to an exact Lagrangian $L'_1$ that is  transverse to $L$. We take the perturbation to be small enough so that there is a primitive $f_{L'_1}$ of $\lambda\mid_{L'_1}$ which vanishes in $\partial L_1'$ and satisfies 
\begin{equation} \label{eq:control:LtoL'_1}
\mathsf{K}(M_{1},L \to L'_1) < \epsilon.
\end{equation}
A continuation argument identical to the one used in the proof of \eqref{eq:isocont} implies that for $a>\epsilon$ 
 the homologies  $\HW^a(M, L \to L_1)$ and $\HW^a(M, L \to L'_1)$ are isomorphic. 

We let 
\begin{equation}
\mathsf{C}_{\regium}:= \#( L'_1 \cap L).
\end{equation}
 This number will be useful later for estimates of the growth of the number of Reeb chords.

We now consider a tubular neighbourhood $\widetilde{\mathcal{U}}(\Lambda_1) $ which does not intersect $\Lambda$. By Lemma \ref{lemmaBaptiste}  there exists $\delta_2>0$ such that if a Legendrian sphere $\Lambda_2$ is $\delta_2$-close to $\Lambda_1$ in the  $C^3$-sense, then there exist  exact Lagrangian cobordisms $\mathcal{L}_{2 \to 1}$ from $\Lambda_2$ to $\Lambda_1$,  and $\mathcal{L}_{1 \to 2}$,  from $\Lambda_1$ to $\Lambda_2$, both contained in the symplectization of $\alpha_{M}$.
 It follows from Lemma \ref{lemmaBaptiste} that by taking $\delta_2>0$ smaller, if necessary, we can guarantee that
\begin{itemize}
\item $\mathcal{L}_{2 \to 1}$ is conical outside $[1-\frac{\epsilon}{5}, 1-\frac{\epsilon}{6}] \times \Sigma$,
\item $\mathcal{L}_{1 \to 2}$ is conical outside $[1+\frac{\epsilon}{6}, 1+\frac{\epsilon}{5}] \times \Sigma$,
\item  the projections of  $\mathcal{L}_{2 \to 1}$ and  $\mathcal{L}_{1 \to 2}$ to $\Sigma$ are contained in $\widetilde{\mathcal{U}}(\Lambda_1) $,
\item there exist primitives $f_{2 \to 1}$ and $f_{1 \to 2}$  of $r\alpha_{M} |_{\mathcal{L}_{2 \to 1}}$ and $r\alpha_{M} |_{\mathcal{L}_{1 \to 2}}$, respectively, with support in $[1-\frac{\epsilon}{5}, 1-\frac{\epsilon}{6}] \times \Sigma$ and $[1+\frac{\epsilon}{6}, 1+\frac{\epsilon}{5}] \times \Sigma$, respectively, such that $|f_{2 \to 1}|_{C^0}<\epsilon$ and $|f_{1 \to 2}|<\epsilon$,
\item the exact Lagrangian $\mathcal{L}_{1 \to 1}$ in the symplectisation   of  $\alpha_{M}$ obtained by gluing  $\mathcal{L}_{1 \to 2} \cap ({[1,+\infty)\times \Sigma})$ on top of $\mathcal{L}_{2 \to 1} \cap ((0,1]\times \Sigma)$ is Hamiltonian isotopic to $(0,+\infty)\times \Lambda_1$ for an isotopy which is stationary outside  ${(1-\frac{\epsilon}{5},1+\frac{\epsilon}{5})\times \Sigma}$.
\end{itemize}
 It is clear that one can glue $f_{1 \to 2}$ and $f_{2 \to 1}$ to obtain a primitive $f_{1 \to 1}$ of $r\alpha_{M} |_{\mathcal{L}_{1 \to 1}}$ which satisfies $|f_{1 \to 1}| < \epsilon$.

We then glue $\mathcal{L}_{2 \to 1} \cap ([1-\frac{\epsilon}{5}, 1] \times \Sigma)$ on top of $L_1 \subset M_{1-\frac{\epsilon}{5}}$ to obtain an asymptotically conical exact Lagrangian  $L_2$ with $L_2\cap \partial M= \Lambda_2$. Let $L'_2$ be the exact Lagrangian submanifold obtained from gluing $\mathcal{L}_{2 \to 1} \cap ([1-\frac{\epsilon}{5}, 1] \times \Sigma)$ on top of $L'_1 \subset M_{1-\frac{\epsilon}{5}}$. It is clear that $L'_2$ and $L_2$ are Hamiltonian isotopic for a Hamiltonian which has support contained in $M_{1 - \frac{\epsilon}{6}}$. 

 Notice that the intersection points of $L'_2 \cap L$ are the same as the intersection points of  $ L'_1 \cap L$. We thus conclude:
\begin{equation}
 \#(L'_2 \cap L) = \mathsf{C}_{\regium}.
\end{equation}
We can glue $f_{2 \to 1}$ to the restriction of $f_{L_1}$ to $L_1 \cap M_{1-\frac{\epsilon}{5}}$ to obtain a primitive $f_{L_2}$ of $\lambda\mid_{L_2}$ such that 
\begin{equation} \label{eq:controlLtoL_2}
\mathsf{K}(M,L \to L_2) < \epsilon.
\end{equation}
Similarly, one obtains a primitive $f_{L'_2}$ of $\lambda\mid_{L'_2}$ such that 
\begin{equation} \label{eq:controlLtoL'_2}
\mathsf{K}(M,L \to L'_2) < \epsilon.
\end{equation}
Assuming that $(\alpha_{M}, \Lambda \to \Lambda_2)$ is regular the Lagrangian $L_2$ is admissible for the pair $(M,M_{1-\frac{\epsilon}{5}})$.
We then obtain for each $a>\epsilon$ a transfer map $\Psi^a_{\mathcal{L}_{2\to 1}}: \HW^a(M, L \to L_2) \to  \HW^a(M_{1-\frac{\epsilon}{5}}, L \to L_1)$. These induce a map $ \Psi_{\mathcal{L}_{2\to 1}}: \HW(M, L \to L_2) \to  \HW(M_{1-\frac{\epsilon}{5}}, L \to L_1)$.

By \eqref{eq:controlLtoL_2} and \eqref{eq:controlLtoL'_2} and the fact that $L_2$ and $L'_2$ are Hamiltonian isotopic for an isotopy supported inside  $M_{1-\frac{\epsilon}{5}}$, we can apply the reasoning used to prove  \eqref{eq:isocont} to show that for each $a>\epsilon$ 
\begin{equation} \label{equsefulestimate}
 \HW^a(M, L \to L_2) \mbox{ and } \HW^a(M, L \to L'_2) \mbox{ are isomorphic.}
\end{equation}

Gluing $\mathcal{L}_{1 \to 2} \cap ([1,1+\frac{\epsilon}{5}] \times \Sigma)$ on top of $L_2 \subset M$ we obtain an asymptotically conical Lagrangian $\widetilde{L}_1$ in $M_{1+\frac{\epsilon}{5}}$ which is transfer admissible for the pair  $(M_{1+\frac{\epsilon}{5}},M)$. Reasoning as in the proof of \eqref{eq:control:LtoL_1} we obtain that 
\begin{equation}
\mathsf{K}(M_{1+\frac{\epsilon}{5}},L \to \widetilde{L}_1) < \epsilon.
\end{equation}
 We thus get for each $a>\epsilon$ a transfer map $\Psi^a_{\mathcal{L}_{1\to 2}}: \HW^a(M_{1+\frac{\epsilon}{5}}, L \to \widetilde{L}_1) \to  \HW^a(M, L \to L_2)$, and in the direct limit a homomorphism $\Psi_{\mathcal{L}_{1\to 2}}: \HW(M_{1+\frac{\epsilon}{5}}, L \to \widetilde{L}_1) \to  \HW(M, L \to L_2)$.

We finally glue $\mathcal{L}^+  \cap ([1 +\frac{\epsilon}{5} ,1 + {\epsilon}]\times \Sigma)$ on top of  $\widetilde{L}_1$  to obtain an asymptotically conical exact Lagrangian $\widetilde{L}$ on $M_{1+\epsilon}$. The Lagrangian $\widetilde{L}$ is an exact filling of $\Lambda$. It is clear from Lemma \ref{lemmaBaptiste} that $\widetilde{L}$ is Hamiltonian isotopic to $L$, for a Hamiltonian which has support contained in $M_{1+\frac{\epsilon}{2}} \setminus M_{1-\frac{\epsilon}{2}}$. Reasoning as in the proof of \eqref{eq:control:LtoL_1} we obtain a primitive $f_{\widetilde{L}}$ of $\lambda \mid_{\widetilde{L}}$ such that 
\begin{equation} \label{eq:controlLtoLhat}
\mathsf{K}(M_{1+\epsilon},L \to \widetilde{L}) < \epsilon.
\end{equation}
We claim that for every $a>\epsilon$ there exists an isomorphism
\begin{equation}\label{eq:isocrucial}
  \Psi^a_{L, \widetilde{L}}: \HW^a(M_{1+\epsilon}, L ) \to \HW^a(M_{1+\epsilon}, L \to \widetilde{L}).
 \end{equation}
 To establish this claim we first notice that if $H$ is a Hamiltonian in\footnote{By \cite[Lemma 8.1]{AbouzaidSeidel2010} any admissible Hamiltonian in $M_{1+\epsilon}$ can be perturbed to one in $\mathcal{H}_{\reg}(M_{1+\epsilon},L\to \widetilde{L}) \cap \mathcal{H}_{\reg}(M_{1+\epsilon},L)$.} $\mathcal{H}_{\reg}(M_{1+\epsilon},L\to \widetilde{L}) \cap \mathcal{H}_{\reg}(M_{1+\epsilon},L)$ it follows from the fact that  $\widetilde{L}$ is Hamiltonian isotopic to $L$ for a Hamiltonian which has support contained in $M_{1+\frac{\epsilon}{2}} \setminus M_{1-\frac{\epsilon}{2}}$ that there exists a continuation isomorphism  $\Psi_{H,L, \widetilde{L}}: \HW(H, L ) \to \HW(H, L \to \widetilde{L})$. Equation \eqref{eq:isocrucial} then follows from combining these isomorphisms and the identifications $\HW^a(M_{1+\epsilon}, L ) \cong \dlim_{H\prec a}\HW(H; L)$ and $\HW^a(M_{1+\epsilon}, L ) \cong \dlim_{H\prec a}\HW(H; L \to \widetilde{L})$ for  $a>\epsilon\geq \max\{\mathsf{K}(M_{1+\epsilon},L \to \widetilde{L}); \mathsf{K}(M_{1+\epsilon},L)\}$ which were established in \eqref{eq:coisabonita}. The maps  $\Psi^a_{L, \widetilde{L}}$ are compatible with the persistence morphisms of the f.d.s. $\widetilde{\HW}(M_{1+\epsilon}, L )$ and $\widetilde{\HW}(M_{1+\epsilon}, L \to \widetilde{L})$ and induce an asymptotic morphism between them. On the direct limit we get a map
 \begin{equation}\label{eq:asymp}
  \Psi_{L, \widetilde{L}}: \HW(M_{1+\epsilon}, L ) \to \HW(M_{1+\epsilon}, L \to \widetilde{L}).
 \end{equation}
 
 
 The succession of exact Lagrangian submanifolds we constructed is schematically presented in Figure 2.
 
  \input{figure3.tex} \label{figure3}
 

 Since  $\widetilde{L}$ is transfer admissible for $(M_{1+\epsilon}, M_{1+\frac{\epsilon}{5}})$ we also obtain for each $a>\epsilon\geq \mathsf{K}(M_{1+\epsilon},L \to \widetilde{L})$ a transfer map $\Phi^a_{\mathcal{L}^+} :  \HW^a(M_{1+\epsilon}, L \to \widetilde{L}) \to  \HW^a(M_{1+\frac{\epsilon}{5}}, L \to \widetilde{L}_1) $. This induces a homomorphism $\Phi_{\mathcal{L}^+} :  \HW(M_{1+\epsilon}, L \to \widetilde{L}) \to  \HW(M_{1+\frac{\epsilon}{5}}, L \to \widetilde{L}_1) $.
 
 Analogously, it follows from Lemma \ref{lemmaBaptiste} that $L_1$ is transfer admissible for the pair $(M_{1-\frac{\epsilon}{5}}, M_{1-\epsilon})$, which gives us for each $a>\epsilon\geq \mathsf{K}(M,L \to L_1) \geq \mathsf{K}(M_{1-\frac{\epsilon}{5}},L \to L_1)$ a map $\Phi^a_{\mathcal{L}^-} :  \HW^a(M_{1-\frac{\epsilon}{5}}, L \to {L}_1) \to  \HW^a(M_{1-\epsilon}, L)$.  These homomorphisms induce a homomorphism $\Phi_{\mathcal{L}^-} :  \HW(M_{1-\frac{\epsilon}{5}}, L \to {L}_1) \to  \HW(M_{1-\epsilon}, L)$.

 The following lemma will be important for the study of the growth rate of $\HW(M, L \to L_2)$.

\begin{lem} \label{lem:cacareco}
For $0<\delta_1$ and $0<\delta_2$ chosen as above we have that the spectral number of $\Psi_{L, \widetilde{L}}( \mathbf{1}_L)$ is $\leq \epsilon$.
\end{lem}
\textit{Proof:} 
We know from \cite{Ritter2013} that $ c(\mathbf{1}_L) = 0$. This implies that for every $a\geq 0$ the element $\mathbf{1}_L $ is in the image of $\iota_a: \HW^a(M_{1+\epsilon},L) \to \HW(M_{1+\epsilon},L)$. 

Let $a>\epsilon$. As remarked above, the maps $\Psi^a_{L, \widetilde{L}}$ are compatible with the persistence morphisms of $\widetilde{\HW}(M_{1+\epsilon}, L )$ and $\widetilde{\HW}(M_{1+\epsilon}, L \to \widetilde{L})$, which implies that the diagram
\\
\[
\begin{CD}
 \HW^a(M_{1+\epsilon}, L ) @>\Psi^a_{L, \widetilde{L}}>> \HW^a(M_{1+\epsilon}, L \to \widetilde{L})  \\
  @V\iota_aVV @V\iota_aVV \\  
 \HW(M_{1+\epsilon}, L) @>\Psi_{L, \widetilde{L}}>>\HW(M_{1+\epsilon}, L \to \widetilde{L})
\end{CD}
\]
\\
is commutative. It follows that $\Psi_{L, \widetilde{L}}(\mathbf{1}_L)$ is in the image of $\iota_a: \HW^a(M_{1+\epsilon}, L \to \widetilde{L}) \to \HW(M_{1+\epsilon}, L \to \widetilde{L})$, from what we obtain that $c(\Psi_{L, \widetilde{L}}(\mathbf{1}_L)) \leq a$. Since this is true for every $a>\epsilon$ we conclude that $c(\Psi_{L, \widetilde{L}}(\mathbf{1}_L)) \leq \epsilon$. \qed

By our discussion so far we have transfer maps
\\
\[
\begin{CD}
\HW(M_{1+\epsilon}, L \to \widetilde{L}) @>\Phi_{\mathcal{L}^+}>>  \HW(M_{1+\frac{\epsilon}{5}}, L \to \widetilde{L}_1) 
@>\Psi_{\mathcal{L}_{1 \to 2}}>> \HW(M, L \to {L}_2)  \\ @. @. @V\Psi_{\mathcal{L}_{2\to 1}}VV \\ @. \HW(M_{1-\epsilon}, L) @<\Phi_{\mathcal{L}^-}<< \HW(M_{1-\frac{\epsilon}{5}}, L \to L_1)
\end{CD}
\]
\\
  Using the fact that $\widetilde{L}$ is Hamiltonian isotopic to $L$ by a Hamiltonian with support contained in $M_{1+\frac{\epsilon}{2}} \setminus M_{1-\frac{\epsilon}{2}}$ and reasoning identically as in the proof of Lemma \ref{lemmaprelim} we conclude that the composition $\Phi_{\mathcal{L}^-} \circ \Psi_{\mathcal{L}_{2\to 1}} \circ \Psi_{\mathcal{L}_{1 \to 2}} \circ \Phi_{\mathcal{L}^+} \circ \Psi_{L, \widetilde{L}}$  is induced by an asymptotic isomorphism  from $\HW(M_{1+\epsilon}, L )$ to $\HW(M_{1-\epsilon}, L )$.
It follows that  $\Psi_{\mathcal{L}_{1 \to 2}} \circ \Phi_{\mathcal{L}^+} :  \HW(M_{1+\epsilon}, L \to \widetilde{L}) \to  \HW(M, L \to {L}_2)$ is an injective $A_L$-module homomorphism. 
We define $m_{L_2}:= \Psi_{\mathcal{L}_{1 \to 2}} \circ \Phi_{\mathcal{L}^+}(\Psi_{L, \widetilde{L}}( \mathbf{1}_L))$.
The element  $m_{L_2} \in \HW(M, L \to {L}_2)$ is stretching since it is the image of a stretching element by an injective $A_L$-module homomorphism.


 By the behaviour of spectral numbers under transfer maps, combined with \eqref{eq:controlLtoLhat} and Lemma \ref{lem:cacareco} we conclude that 
\begin{equation} \label{uniformitymodules}
c(m_{L_2}) \leq \max\{c(\Psi_{L, \widetilde{L}}( \mathbf{1}_L)),\mathsf{K}(M_{1+\epsilon}, L \to \widetilde{L}) \} \leq \epsilon.
\end{equation}

We denote by  $\mathcal{V}_{\alpha_{M}}(\Lambda_1)$ the set of Legendrian spheres $\Lambda_2$ in the same Legendrian isotopy class of $\Lambda_1$ that are $\delta_2$-close to $\Lambda_1$ is the $C^3$-sense.
Let $\mathcal{V}_{\alpha_{M}}^{  \reg}(\Lambda_1) \subset \mathcal{V}_{\alpha_{M}}(\Lambda_1)$ be the subset of these $\Lambda_2$ for which, in addition,  $(\alpha_{M}, \Lambda \to \Lambda_2)$ is regular. We denote by $L_2$ the filling of $\Lambda_2$ constructed above.
 Our discussion so far implies the following
\begin{prop} \label{uniformmodules2} 
The family  $(\HW(M,L \to L_2))_{\Lambda_2 \in \mathcal{V}_{\alpha_{M}}^{\reg}(\Lambda_1)}$ of $A_L$-modules is uniformly stretched. It follows from Lemma \ref{mod_fam_growth}, Lemma \ref{lem:subad}, and Lemma \ref{lem:subad2} that 
\begin{equation}
\Gamma_{\Lambda_2 \in \mathcal{V}_{\alpha_{M}}^{\reg}(\Lambda_1)}^{\symp}(M,L \to L_2) \geq  \Gamma^{\symp}(M, L).
\end{equation} 

\end{prop}
\textit{Proof:} 
The proposition follows directly from the fact that the element $m_{L_2} \in \HW(M,L \to L_2)$ is stretching and from  \eqref{uniformitymodules}. \qed

Let $\alpha$ be a contact form on $(\Sigma,\xi_{M})$. We assume that the function $\mathsf{f}_\alpha$ defined by $\alpha= \mathsf{f}_\alpha \alpha_{M}$ satisfies $\mathsf{f}_\alpha \geq 1$. We thus have the inclusions $ M_{\mathsf{f}_\alpha} \subset M_{\max \mathsf{f}_\alpha}$ and $M \subset  M_{\mathsf{f}_\alpha}$.

We denote by $\mathcal{V}_{\alpha_{M}}^{ \alpha-\reg}(\Lambda_1) \subset \mathcal{V}_{\alpha_{M}}^{  \reg}(\Lambda_1)$ the set of $\Lambda_2 \in \mathcal{V}_{\alpha_{M}}^{  \reg}(\Lambda_1)$ such  that $(\alpha, \Lambda \to \Lambda_2)$ is regular.

Let $W_\alpha^+:= M_{\max \mathsf{f}_\alpha} \setminus M_{\mathsf{f}_\alpha}$ and $W_\alpha^-:= M_{ \mathsf{f}_\alpha} \setminus M$.
Since the Lagrangians $L_1$ and $L_2$ are conical in $M_{\max \mathsf{f}_\alpha} \setminus M$ we obtain for elements $\Lambda_2 \in \mathcal{V}_{\alpha_{(\Sigma,\lambda)}}^{ \alpha-\reg}(\Lambda_1)$ transfer maps
\\
\[
\begin{CD}
 \HW(M_{\max \mathsf{f}_\alpha}, L\to L_2)  @>\Phi_{W^+,L \to L_2}>>  \HW(M_{\mathsf{f}_\alpha},L  \to L_2) @>\Phi_{W^-,L \to L_2}>>  \HW(M,L  \to L_2).
\end{CD}
\]
\\
By Lemma \ref{lem:changeofhypersurface}, the composition $\Phi_{W^-,L \to L_2}\circ \Phi_{W^+,L \to L_2}$  is induced by asymptotic morphisms,  and the f.d.s. $ \widetilde{\HW}(M_{\mathsf{f}_\alpha},L  \to L_2)$ and $\widetilde{\HW}(M,L  \to L_2)$ are $(\max \mathsf{f}_{\alpha}, 1)$- interleaved.

The following proposition then follows from combining this observation and Proposition \ref{uniformmodules2}.
\begin{prop} \label{uniformgrowth}
Let $\alpha$ be a contact form on $(\Sigma,\xi_{M})$ and assume that the function $\mathsf{f}_\alpha$ defined by $\alpha=\mathsf{f}_\alpha \alpha_{M}$ is $\geq 1$.
Then, the family of f.d.s.  $(\widetilde{\HW}(M_{\mathsf{f}_\alpha},L \to L_2))_{\Lambda_2 \in \mathcal{V}_{\alpha_{M}}^{\alpha-\reg}(\Lambda_1)}$ satisfies
\begin{equation} \label{equniformgrowth1}
\Gamma^{\symp}_{\Lambda_2 \in \mathcal{V}_{\alpha_{M}}^{\alpha-\reg}(\Lambda_1)}(M_{\mathsf{f}_\alpha},L \to  L_2) \geq   \frac{\Gamma^{\symp}(M, L)}{\max \mathsf{f}_\alpha}.
\end{equation} 
\end{prop}

A reasoning identical to the one used to  establish  \eqref{equsefulestimate} shows that for every $\Lambda_2 \in  \mathcal{V}_{\alpha_{M}}^{\alpha-\reg}(\Lambda_1)$ and for the exact filling $L'_2$ of $\Lambda_2$ constructed above we have 
\begin{equation}
\HW^a(M_{\mathsf{f}_\alpha}, L \to L_2) \mbox{ and } \HW^a(M_{\mathsf{f}_\alpha}, L \to L'_2)  \mbox{ are isomorphic.}
\end{equation}
Combining this with Proposition \ref{uniformgrowth} we have
\begin{cor} \label{coro:uniformgrowth}
Let $\alpha$ be a contact form on $(\Sigma,\xi_{M})$ and assume that the function $\mathsf{f}_\alpha$ defined by $\alpha=\mathsf{f}_\alpha \alpha_{M}$ is $\geq 1$.
Then, the family of f.d.s.  $(\widetilde{\HW}(M_{\mathsf{f}_\alpha},L \to L'_2))_{\Lambda_2 \in \mathcal{V}_{\alpha_{M}}^{\alpha-\reg}(\Lambda_1)}$ satisfies
\begin{equation} \label{equniformgrowth1}
\Gamma^{\symp}_{\Lambda_2 \in \mathcal{V}_{\alpha_{M}}^{\alpha-\reg}(\Lambda_1)}(M_{\mathsf{f}_\alpha},L \to  L'_2) \geq   \frac{\Gamma^{\symp}(M, L)}{\max \mathsf{f}_\alpha}.
\end{equation} 
\end{cor}

Recall that for every $\Lambda_2 \in  \mathcal{V}_{\alpha_{M}}^{\alpha-\reg}(\Lambda_1)$ the exact filling $L'_2$ of $\Lambda_2$ satisfies
\begin{equation} \label{bound}
 \#(L'_2 \cap L) = \mathsf{C}_{\regium}.
\end{equation}

Now, given a Legendrian $\Lambda_2 \in  \mathcal{V}_{\alpha_{M}}^{\alpha-\reg}(\Lambda_1)$ let $\mathsf{N}_\alpha^a(\Lambda \to \Lambda_2)= \# \mathcal{T}^a_{\Lambda \to \Lambda_2}(\alpha)$.
We define 
\begin{equation}
\mathsf{N}^a_\alpha( \Lambda \to \mathcal{V}_{\alpha_{M}}^{\alpha-\reg}(\Lambda_1)):= 
\inf_{\Lambda_2 \in  \mathcal{V}_{\alpha_{M}}^{\alpha-\reg}(\Lambda_1)} \{\mathsf{N}^a_\alpha(\Lambda \to \Lambda_2)   \}.
\end{equation}

Let $a > \epsilon$. By the results of Section \ref{rem:rem1} there exists a Hamiltonian $H^a\in \mathcal{H}_{\reg}(M_{\mathsf{f}_\alpha},L \to L'_2)$ with slope $<a$ such that 
\begin{itemize}
\item[(\textbf{p.1})] all elements in $\mathcal{T}_{L \to L'_2}(H^a)$ have action $<a$,
\item[(\textbf{p.2})] there is a bijection  between $\mathcal{T}_{L \to L'_2}(H^a)$ and $\mathcal{T}_{\Lambda \to \Lambda_2}(\alpha)\cup (L\cap L'_2)$,
\item[(\textbf{p.3})] the map $\chi^{H^a}_{ \to a}:\HW(H^a,L\to L'_2) \to  \HW^{a}(M_{f_\alpha},L\to L'_2)$ is isomorphism.
\end{itemize} 
Combining  (\textbf{p.3}) and Corollary \ref{coro:uniformgrowth} we obtain that 
\begin{equation} \label{eq:tatatadada}
\limsup_{a \to +\infty} \frac{\log (\inf_{\Lambda_2 \in  \mathcal{V}_{\alpha_{M}}^{\alpha-\reg}(\Lambda_1)}\{\dim\CW(H^{a} ,L\to L'_2)\})}{a} \geq \frac{\Gamma^{\symp}(M, L)}{\max \mathsf{f}_\alpha}.
\end{equation}
Since $\#(L\cap L'_2)=C_{\regium}$ it follows from (\textbf{p.2}) that $\dim(\CW^{a}(H^{{a}},L\to L'_2)) - \mathsf{C}_{\regium} = \mathsf{N}_\alpha^a(\Lambda \to \Lambda_2)$. This together with \eqref{eq:tatatadada} gives
\begin{cor} \label{corocrucial}
The sequence of numbers $\mathsf{N}^a_\alpha( \Lambda \to \mathcal{V}_{\alpha_{M}}^{\alpha-\reg}(\Lambda_1))$ satisfies
\begin{equation}
\limsup_{a \to +\infty} \frac{\log \mathsf{N}^a_\alpha( \Lambda \to \mathcal{V}_{\alpha_{M}}^{\alpha-\reg}(\Lambda_1))}{a} \geq \frac{\Gamma^{\symp}(M, L)}{\max \mathsf{f}_\alpha}.
\end{equation}
\end{cor}
This corollary will be crucial for the estimate of the topological entropy of $\phi_\alpha$ proved next.  \color{black}

\subsection{From the growth of Reeb chords to topological entropy}

Let $\alpha$ be a contact form on a contact manifold $(\Sigma,\xi)$, and $X_\alpha$ be its Reeb vector field. Recall that a Riemannian metric $g$ on $X$ is said to be compatible with $\alpha$ if $g(X_\alpha,X_\alpha)=1$ and $X_\alpha$ is orthogonal to $\xi$ with respect to $g$.

We proceed by fixing some more notation. We denote by $\mathbb{D}^n(\rho)$ the $n$-dimensional disk of radius $\rho>0$ around the origin. We endow $\mathbb{D}^n(\rho)$ with the Euclidean metric, and consider on $T^*_1\mathbb{D}^n(\rho) = \mathbb{D}^n(\rho) \times S^{n-1}$ the contact form $\alpha_{\euc}$ associated to the Euclidean metric. For each $z\in  \mathbb{D}^n(\rho)$ the sphere $S^{n-1}_z:= \{z\} \times S^{n-1}$ is Legendrian in $(\mathbb{D}^n(\rho) \times S^{n-1},\ker \alpha_{\euc}) $. 
Let $g_{\round}$ be the metric with constant curvature $1$ on $S^{n-1}$ and $g_{\euc}$ be the Euclidean metric on $\mathbb{D}^{n}(\rho)$. The metric $\widetilde{g}= g_{\euc} \oplus g_{\round}$ on $\mathbb{D}^n(\rho) \times S^{n-1}$ is compatible with the contact form $\alpha_{\euc}$; see \cite{calvaruso}.

\begin{prop} \label{propgrowth}
Let $\alpha$ be a contact form on $(\Sigma,\xi_{M})$ and assume that we have $\Gamma^{\alg}_S(M_{f_\alpha}, L) >0$. Then there exists a Riemannian metric $g$ on $(\Sigma,\xi_{M})$ adapted to the $\alpha$, such that
\begin{equation}
\limsup_{t \to +\infty} \frac{ \log \Vol_g^{n-1}(\phi_{\alpha}^t(\Lambda))}{t} \geq \frac{\Gamma^{\symp}(M, L)}{\max f_\alpha}>0,
\end{equation}
where $\Vol_g^{n-1}$ is the $(n-1)$-dimensional volume with respect to $g$, and $f_\alpha$ is the function such that $\alpha=f_\alpha \alpha{(\Sigma,\lambda)}$.
\end{prop}
\textit{Proof:} 
The proof will consist of several steps. 

\textbf{Step 1.} It suffices to prove the proposition for all contact forms $\alpha$ for which $\mathsf{f}_\alpha\geq1$. 
Indeed assume that the proposition holds  for all such contact forms.

Take a contact form ${\alpha'}$ on  $(\Sigma,\xi_{M})$. For the contact form $\widehat{\alpha}:= \frac{\alpha'}{\min \mathsf{f}_{\alpha'}}$ we have $\mathsf{f}_{\widehat{\alpha}}=  \frac{\mathsf{f}_{\alpha'}}{\min \mathsf{f}_{\alpha'} }\geq 1$. By assumption there is a Riemannian metric $g$ on $\Sigma$ compatible with $\widehat{\alpha}$ and such that 
\begin{equation}
\limsup_{t \to +\infty} \frac{ \log \Vol_g^{n-1}(\phi_{\widehat{\alpha}}^t(\Lambda))}{t} \geq \frac{\Gamma^{\symp}(M, L)}{\max \mathsf{f}_{\widehat{\alpha}}}.
\end{equation}

The Riemannian metric $g':=(\min \mathsf{f}_{{\alpha'}})^2g$ is compatible with $\alpha'$. 
A simple computation shows that 
$\limsup_{t \to +\infty} \frac{ \log \Vol_{g'}^{n-1}(\phi_{{\alpha'}}^t(\Lambda))}{t} \geq \frac{\Gamma^{\symp}(M, L)}{\max \mathsf{f}_{{\alpha'}}}$, as claimed.
We thus fix from now on a contact form $\alpha$ on $(\Sigma,\xi_{M})$  with $\mathsf{f}_\alpha\geq1$.

\textbf{Step 2. A tubular neighbourhood of $\Lambda_1$ and construction of the metric $g$} \\
It follows from the Legendrian neighbourhood theorem (see \cite[Proposition 43.18]{MichorKriegl}) that there exists a tubular neighbourhood  $\mathcal{V}(\Lambda_1)$ of $(\Lambda_1)$ and a contactomorphism 
$\Upsilon : (\mathcal{V}(\Lambda_1), \xi_{M}) \to (\mathbb{D}^n(\rho) \times S^{n-1},\ker \alpha_{\euc}) $ that satisfies
\begin{eqnarray} \label{neighbourhood}
\Upsilon^* \alpha_{\euc} =  \alpha, \\
\Upsilon(\Lambda_1) = \{0\} \times S^{n-1}.
\end{eqnarray}

We extend the Riemannian metric $\Upsilon^*\widetilde{g}$, which is compatible with $\alpha$ on $\mathcal{V}(\Lambda_1)$, to a metric $g$ on $\Sigma$ which is compatible with the contact form $\alpha$.

After shrinking the neighbourhood $\mathcal{V}(\Lambda_1)$ and $\rho>0$, we can assume that for every $z \in \mathbb{D}^{n}(\rho)$ the Legendrian $\Lambda^z := \Upsilon^{-1}(\{z\} \times S^{n-1})$ is  in the neighbourhood $\mathcal{V}_{\alpha_{M}}(\Lambda_1)$ constructed in Section~\ref{subsec:technical}.

\textbf{Step 3.} 
For each $a>0$ we define the map $F^a_{\Lambda}: \Lambda \times [0,a] \to \Sigma$ by
\begin{equation}
F^a_{\Lambda}(q,t) = \phi^t_\alpha(q).
\end{equation}
Let $\Cyl^a_\alpha(\Lambda)$ be the image $F^a_{\Lambda}(\Lambda \times [0,a] )$.
We want to estimate from below the $n$-dimensional volume $\Vol^{n}_g(\Cyl^a_\alpha(\Lambda))$ of $\Cyl^a_\alpha(\Lambda)$ with respect to the Riemannian metric $g$. For this we define $\mathfrak{B}^a_\alpha(\Lambda):= \Upsilon(\Cyl^a_\alpha(\Lambda) \cap \mathcal{V}(\Lambda_1)) $. We have
\begin{equation} \label{tititi}
\Vol^{n}_g(\Cyl^a_\alpha(\Lambda)) \geq \Vol^{n}_{{g}}(\Cyl^a_\alpha(\Lambda)\cap \mathcal{V}(\Lambda_1)) =  \Vol^{n}_{\widetilde{g}}(\mathfrak{B}^a_\alpha(\Lambda)).
\end{equation}

Let $\Pi: \mathbb{D}^{n}(\rho) \times S^{n-1} \to \mathbb{D}^{n}(\rho)$ be the projection to the first coordinate. 
Applying Sard's Theorem to the map $\Pi \circ \Upsilon \circ F^a_\Lambda: (\{a\} \times \Lambda) \cap (F^a_\Lambda)^{-1}( \mathcal{V}(\Lambda_1)) \to   \mathbb{D}^n(\rho)$ we conclude that the set $ \mathbb{D}^{n}(\rho) \setminus \Pi\circ \Upsilon(\phi^a_\alpha(\Lambda))$ is an open set of full Lebesgue measure in $ \mathbb{D}^{n}(\rho)$.
We define the set $\mathfrak{U}^a_\alpha(\Lambda) \subset \mathbb{D}^{n}(\rho) \setminus \Pi\circ \Upsilon(\phi^a_\alpha(\Lambda) )$ by the property
\begin{itemize}
\item $z \in \mathfrak{U}^a_\alpha(\Lambda)$ if all $\alpha$-Reeb chords from $\Lambda$ to $\Lambda^z$ with length  $<a$ are transverse.  
\end{itemize}

The proof of the next lemma is identical to the one of \cite[Lemma 4]{Alves-Legendrian}.
\begin{lem}
The set $\mathfrak{U}^a_\alpha(\Lambda)$ is an open subset of $\mathbb{D}^{n}(\rho)$ of full Lebesgue measure. The set  $\widetilde{\mathfrak{U}}^a_\alpha(\Lambda) \subset \mathfrak{U}^a_\alpha(\Lambda)$ of elements $z \in  \mathfrak{U}^a_\alpha(\Lambda)$ such that  $\Lambda^z \in  \mathcal{V}_{\alpha_{M}}^{\alpha-\reg}(\Lambda_1)$ is a dense subset of full Lebesgue measure in ${\mathfrak{U}}^a_\alpha(\Lambda)$.
\end{lem}

\textbf{Step 4. A volume estimate} \\
The function $h^a: \mathfrak{U}^a_\alpha(\Lambda) \to [0, +\infty)$ defined by
$h^a(z) := \#(\mathcal{T}^a_{\Lambda \to \Lambda^z}(\alpha)) $
is locally constant  on $\mathfrak{U}^a_\alpha(\Lambda)$ since it is continuous and takes only integer values. 

We define $\mathfrak{R}^a_\alpha(\Lambda) :=  \Pi^{-1}(\mathfrak{U}^a_\alpha(\Lambda))\cap  \mathfrak{B}^a_\alpha(\Lambda)$. Since $\mathfrak{R}^a_\alpha(\Lambda) \subset \mathfrak{B}^a_\alpha(\Lambda)$
we have $\Vol^{n}_{\widetilde{g}}(\mathfrak{B}^a_\alpha(\Lambda)) \geq \Vol^{n}_{\widetilde{g}}(\mathfrak{R}^a_\alpha(\Lambda))$.
As the map $\Pi: \mathbb{D}^{n}(\rho) \times S^{n-1} \to \mathbb{D}^{n}(\rho)$ is a Riemannian submersion we have that $\Vol^{n}_{\widetilde{g}}(\mathfrak{R}^a_\alpha(\Lambda)) \geq \Vol^{n}_{g_{\euc}}(\Pi(\mathfrak{R}^a_\alpha(\Lambda)))$,
where $\Vol^{n}_{g_{\euc}}(\Pi (\mathfrak{B}^a_\alpha(\Lambda)))$  is computed with multiplicities. If an open set is covered $k$-times by $\Pi  : \mathfrak{R}^a_\alpha(\Lambda) \to \mathfrak{U}^a_\alpha(\Lambda)$, then its volume contributes $k$-times to  $\Vol^{n}_{g_{\euc}}(\Pi (\mathfrak{R}^a_\alpha(\Lambda)))$. 

For each $z \in \mathfrak{U}^a_\alpha(\Lambda)$ the number of times  $\Pi: \mathfrak{R}^a_\alpha(\Lambda) \to \mathfrak{U}^a_\alpha(\Lambda)$ covers $z$ is $h^a(z) =\#(\mathcal{T}^a_{\Lambda \to \Lambda^z}(\alpha)) $. We thus obtain 
\begin{equation}
\Vol^{n}_{g_{\euc}}(\Pi(\mathfrak{R}^a_\alpha(\Lambda))) = \int_{ \mathfrak{U}^a_\alpha(\Lambda)} h^a(z)  \ d \mathsf{vol}_{g_{\euc}},
\end{equation}
where $d \mathsf{vol}_{g_{\euc}}$ is the volume form generated by $g_{\euc}$ on $\mathbb{D}^n(\rho)$.

Since $\frac{\Gamma^{\symp}(M, L)}{\max \mathsf{f}_\alpha}>0$, we can fix $0<\eta <\frac{\Gamma^{\symp}(M, L)}{\max \mathsf{f}_\alpha}$. It follows from Corollary \ref{corocrucial} that there exists a sequence $a_j \to +\infty$ such that $h^{a_j}(z) \geq e^{\eta a_j}$ for all $z \in \widetilde{\mathfrak{U}}^{a_j}_\alpha(\Lambda)$. Since $\widetilde{\mathfrak{U}}^{a_j}_\alpha(\Lambda)$ is dense in ${\mathfrak{U}}^{a_j}_\alpha(\Lambda)$ and $h^{a_j}$ is locally constant on ${\mathfrak{U}}^{a_j}_\alpha(\Lambda)$ we obtain $h^{a_j}(z) \geq e^{\eta a_j} \mbox{ for all } z \in {\mathfrak{U}}^{a_j}_\alpha(\Lambda)$
and all $a_j$.
With \eqref{tititi} it follows that 
\begin{equation} \label{volgrowth}
\Vol^{n}_g(\Cyl^{a_j}_\alpha(\Lambda)) \geq \int_{ \mathfrak{U}^{a_j}_\alpha(\Lambda)} h^{a_j}(z)  \ d \mathsf{vol}_{g_{\euc}} \geq e^{\eta a_j} \uppi \rho^2
\end{equation}
for every $a_j$.

\textbf{Step 5. A Fubini type equality}
We define $\widehat{g}:= (F^a_\Lambda)^*g$. Then 
\begin{equation}\label{changecoord}
\Vol^{n}_g(\Cyl^{a}_\alpha(\Lambda))=  \int_{ \Lambda \times [0,a]} d \mathsf{vol}_{\widehat{g}},
\end{equation}
where $d \mathsf{vol}_{\widehat{g}}$ is the volume form associated to $\widehat{g}$.
Since the metric $g$ is adapted to the contact form $\alpha$ the Reeb vector field has length $1$ and is orthogonal to the Legendrian spheres $F^a_\Lambda(t,\Lambda) = \phi^t_\alpha(\Lambda)$ for every $t\in [0,a]$.
Letting $\partial_t$ be the tangent vector field on $[0,a] \times \Lambda$ associated to the first coordinate $t\in [0,a]$, and using the definition of $F^a_\Lambda$, it follows that $D(F^a_\Lambda) \partial_t = X_\alpha$. Therefore  $\partial_t$ has $\widehat{g}$-norm equal to $1$ at every point in $[0,a] \times \Lambda$, and is orthogonal to the spheres $\{t\}\times \Lambda$. We thus conclude that 
\begin{equation} \label{fubini}
 \Vol^{n}_g(\Cyl^{a}_\alpha(\Lambda))=\int_{ \Lambda \times [0,a]} d \mathsf{vol}_{\widehat{g}} = \int_0^{a} \Vol_{\widehat{g}}^{n-1}(\{t\} \times \Lambda) dt = \int_0^a \Vol_{g}^{n-1}( \phi^t_\alpha(\Lambda)) dt,
\end{equation}
where $\Vol_{\widehat{g}}^{n-1}$ is the $(n-1)$-dimensional volume associated to $\widehat{g}$.

\textbf{Step 6. End of the proof.}
To finish the proof we argue by contradiction and assume that 
$\limsup_{t \to +\infty } \frac{\log \Vol_{g}^{n-1}( \phi^t_\alpha(\Lambda))}{t} < \eta$. In this case, there exist $ a_0 >0$ and $\varepsilon>0$ such that for all $t\geq a_0$ we have $\Vol_{g}^{n-1}( \phi^t_\alpha(\Lambda)) \leq e^{t(\eta-\varepsilon)}$. 
Integrating both sides of this inequality from $0$ to $a\geq a_0$ and invoking \eqref{fubini} we obtain
\begin{equation} \label{cacaca}
\Vol^{n}_g(\Cyl^{a}_\alpha(\Lambda))\leq \frac{e^{a (\eta-\varepsilon)} - e^{a_0 (\eta-\varepsilon)}}{\eta-\varepsilon} +  \int_0^{a_0} \Vol_{g}^{n-1}( \phi^t_\alpha(\Lambda)) dt.
\end{equation}
For $a$ large enough the right hand side of \eqref{cacaca} is smaller than $e^{\eta a}\uppi \rho^2$, contradicting \eqref{volgrowth}. We  thus conclude that 
\begin{equation}
\limsup_{t \to +\infty } \frac{\log \Vol_{g}^{n-1}( \phi^t_\alpha(\Lambda))}{t} \geq \eta.
\end{equation}

Since this is valid for any $\eta < \frac{\Gamma^{\symp}(M, L)}{\max \mathsf{f}_{\alpha}}$, the proof of the proposition is completed. \qed

\textit{\textbf{Proof of Theorem \ref{theorementropy}:}} \\
From Proposition \ref{propgrowth} and Yomdin's theorem (see \eqref{Yomdin}) it follows that if $\Gamma^{\symp}(M,L)>0$, then for every contact form $\alpha$ on $(\Sigma,\xi_{M})$ we have
\begin{equation} \label{entropyineq}
h_{\topo}(\phi_{\alpha}) \geq \frac{\Gamma^{\symp}(M,L)}{\max(\mathsf{f}_{\alpha})}.
\end{equation}
We then obtain Theorem \ref{theorementropy} by combining \eqref{entropyineq} with the inequality
\begin{equation*}
\Gamma^{\symp}(M,L) \geq \frac{\Gamma^{\alg}_S(M,L)}{\rho(S)}.
\end{equation*} \color{black}
 from Lemma \ref{prop:alg_symp}. \color{black}
\qed

\color{black}

\section{Algebras in loop space homology}\label{sec:ring_module} 

\color{black}

Let $V$ be a compact manifold and fix a point $q \in V$.
We denote by $\Omega_q(V)$ the based loop space of $V$ with basepoint in $q$, which is the space of continuous maps from $[0,1]$ to $V$ that map $0$ and $1$ to $q$. 

The concatenation of based loops gives $\Omega_q(V)$ the structure of an $H$-space (see \cite{Hatcher}). 
More precisely, the concatenation induces the so-called Pontrjagin product on the singular homology $\Ho_*(\Omega_{q}(V))$ of $\Omega_{q}(V)$ with $\mathbb{Z}_2$ coefficients. \color{black}
The Pontrjagin product $[a_1]\cdot[a_2]$ of two homology classes $[a_1],[a_2] \in \Ho_*(\Omega_{q}(V)$ is well-known to be associative.
As it is distributive with respect to the vector space structure of $\Ho_*(\Omega_{q}(V))$, it makes $\Ho_*(\Omega_{q}(V))$ into a ring. Because the homology $\Ho_*(\Omega_{q}(V))$ is considered with coefficients in $\Z_2$ it actually has the structure of an algebra.

\subsection{Relation between the algebra structure of the singular homology of loop spaces and the algebra structures of the Floer homology of cotangent bundles} \color{black}
 Given a manifold $V$ and  $q \in V$ we denote by $L_q \subset T^* V$ the cotangent fibre over $q$.
The singular homology $\Ho_*(\Omega_{q}(V))$ of the based loop space $\Omega_{q}(V)$ is isomorphic to the wrapped Floer homology $\HW(T^*V,L_q)$; see Viterbo \cite{Viterbo1999}, Salamon-Weber \cite{SW} and Abbondandolo-Schwarz \cite{AS-iso} for different proofs. 

The Floer homology $\HW(H_g,L_q)$ is isomorphic to the wrapped Floer homology $\HW(L_q)$ we use in  this paper. The key point is that the Hamiltonian $H_g$ is quadratic in the fibres. This isomorphism is proven in \cite{Ritter2013}, and it  preserves the triangle product and the spectral value of homology classes.

Let $\Psi_{AS,q}: \Ho_*(\Omega_{q}(V)) \to \HW(T^*V,L_q)$ be the isomorphism constructed in  \cite{AS-iso}.
In \cite{AS-product} the authors proceed to study more properties of the map $\Psi_{AS,q}$. \color{black}
 They show that $\Psi_{AS,q}$ is also algebra isomorphism if we consider $\Ho_*(\Omega_{q}(V))$ as an algebra  with the Pontrjagin product and $\HW(H_g,L_q)$ as an algebra with the triangle product. Combining this with the isomorphism $\HW(H_g,L_q) \cong\HW(L_q)$ we obtain the following 
\begin{thm}[Abbondandolo-Schwarz \cite{AS-product}] \label{ASring}
The singular homology $\Ho_*(\Omega_{q}(V))$ and the wrapped Floer homology $\HW(L_q)$ are isomorphic as algebras. 
\end{thm} \color{black}
For simplicity we will still denote by $\Psi_{AS,q}$ the isomorphism between $\Ho_*(\Omega_{q}(V))$ and $\HW(L_q)$.

\section{Topological operations}\label{sec:top_operations}

\subsection{Subcritical surgery}\label{subsec:subcritical}

Here we study the Viterbo transfer maps under subcritical handle attachment in the situation that is sufficient for our purpose, that is we assume that the Lagrangians do not intersect the handle. 

\color{black} 
Let $W=(Y_W,\omega,\lambda)$ be a Liouville domain, $\Sigma = \partial W$, $\lambda|_{\Sigma} = \alpha$ and $\xi = \ker \alpha$. We recall some notions using the terminology of \cite[Section 2.5.2]{GeigesBook}. The form $d\alpha$ endows $\xi$ with a natural conformal symplectic bundle structure.
Let $S$ be an isotropic submanifold of $(\Sigma,\xi)$. We write $TS^{\bot}$ for the sub-bundle of $\xi$ that is $d\alpha$-orthogonal to $TS$. Because $S$ is isotropic $TS \subset TS^{\bot}$.
We can therefore write the normal bundle of $S$ in $\Sigma$ as 
\[
T \Sigma /T S = T\Sigma/\xi \oplus \xi/TS^{\bot} \oplus TS^{\bot} / TS.
\]  \color{black}
The \textit{conformal symplectic normal bundle} $\mathrm{CSN(S)} = TS^{\bot} / TS$ has a natural conformal symplectic structure via $d\alpha$. If $S$ is a sphere, $T\Sigma/\xi \oplus \xi/TS^{\bot}$ has a trivialization. The following theorem is due to Weinstein.
\begin{thm}\cite{Weinstein1990}
Let $S^n$ be an isotropic sphere in $\Sigma$ with a trivialization of $\mathrm{CSN}(S)$. Then there is a Liouville domain $M$ with an exact embedding $W \subset M$, such that $\partial M$ is obtained from $\Sigma$ by surgery on $S$. 
\end{thm}

The Liouville domain $M$ is obtained by attaching an $(n+1)$-handle to $W$ and the Liouville vector field $X$ can by chosen such that there is exactly one point $p \in M\setminus W$ where $X$ vanishes. The integral lines of $X$ that are asymptotic to $p$ intersect $\Sigma$ in $S$ and $\partial M$ in the co-core sphere $B \subset \partial M$. (See \cite{Weinstein1990, Cieliebak2001} or \cite[Chapter 6]{GeigesBook} for details.) 

Let now $L^{'}_0$, $L^{'}_1$ be  two asymptotically conical exact Lagrangians in $ W$  whose boundaries $\Lambda^{'}_0$ and $\Lambda^{'}_1$ in $\Sigma$ do not intersect $S$.
Outside $S$ the integral lines of the Liouville vector field starting at $\partial W$ intersect $\partial M$ and so the completed Lagrangians $\widehat{L^{'}_{i}} \subset \widehat{M}$  intersect $\partial M$. Moreover, $L_{i} = \widehat{L^{'}_{i}} \cap M \subset M$ for $i=0,1$ are exact and conical in the complement of $W$. We say that $(M,L_0,L_1)$ is \textit{obtained by surgery} from $(W,L^{'}_0,L^{'}_1).$

As described in section \ref{sec:Viterbo} we get a Viterbo transfer map
\[
{j}_{!}(L_0,L_1): \widetilde{\HW}(M,L_0 \to L_1) \rightarrow \widetilde{\HW}(W,L^{'}_0 \to L^{'}_1).  
\]

Assume that the isotropic sphere $S$ has the property that there is no Reeb chord from $\Lambda^{'}_0$ to $S$.
If $S$ is subcritical, i.e. $\dim(S) < n-1$, this can be achieved by a generic choice of $S$. 

The following proposition was proved by Cieliebak (\cite{Cieliebak2001}) for symplectic homology. The proof in our situation is analogous and even simpler. We give it here for the convenience of the reader.  
\begin{prop}\label{Viterbo_iso}
The Viterbo transfer map in the direct limit,
\[
\bar{j}_{!}(L_0,L_1): \HW(M,L_0,L_1) \to \HW(W,L^{'}_0, L^{'}_1), 
\]
is an isomorphism.
\end{prop}

\input{figure2.tex}

For the proof of Proposition \ref{Viterbo_iso} it is convenient to introduce the following weaker form of interleaving of f.d.s. 
Let $\sigma: [0,\infty) \rightarrow [0,\infty)$ be a monotone increasing function and $V$ a filtered directed system. \color{black}  Analogously to the notation in \ref{subsubsec:fds} let $(V(\sigma), \pi(\sigma))$ be given by $V({\sigma})_t = V_{\sigma (t) t}$, $\pi(\sigma)_{s\rightarrow t}= \pi_{\sigma(s) s \rightarrow \sigma(t) t}$ and $\pi[\sigma]_t = \pi_{\sigma(t) t}$. If $f$ is a morphism from $(V,\pi)$ to another f.d.s. we write $f(\sigma)_t = f_{\sigma(t) t}$ for the induced morphism with domain $(V(\sigma), \pi(\sigma))$.   \color{black} 
Call two f.d.s. $(V,\pi_V)$ and $(W,\pi_W)$  \textit{weakly interleaved} if there are morphisms
$f:V\rightarrow W(\sigma_1)$ and $g:W \rightarrow V(\sigma_2)$ for monotone increasing functions $\sigma_1, \sigma_2 \geq 1$ 
such that  
\[f(\sigma_2)\circ g = \pi_{W}[\widetilde{\sigma}_1] \text{ and  } g(\sigma_1)\circ f = \pi_{V}[\widetilde{\sigma}_2],\]
where $\widetilde{\sigma}_1$, and $\widetilde{\sigma}_2$ are suitably chosen.
The fact that the map $\bar{j}_{!}(L_0,L_1)$ in Proposition \ref{Viterbo_iso} is an isomorphism will follow from a weak interleaving of the corresponding f.d.s., which is in general not an interleaving. This is the reason why we cannot directly prove lower bounds for $\Gamma^{\symp}(M,L_0 \to L_1)$ in terms of $\Gamma^{\symp}(W,L^{'}_0 \to L^{'}_1)$ and this was originally our motivation to introduce the algebraic growth of wrapped Floer homology. 

\textit{Proof of Proposition \ref{Viterbo_iso}: }
Let $U = \widetilde{\HW}(M,L_0 \to L_1)$, and $V =\widetilde{\HW}(W,L^{'}_0 \to L^{'}_1)$. 
We will construct a filtered directed system $Q$ that is isomorphic to $V$ and weakly interleaved with $U$. 

For convenience we may assume $\mathsf{K}(M,L_0 \to L_1) = 0$. 
Let $S \subset \partial W$ be the attaching sphere and $B \subset \partial M$ be the co-core sphere. 
For $a > 0$, choose a tubular neighbourhood $U_a \subset \partial W$ of $S$ such that there is no Reeb trajectory starting at $\Lambda^{'}_0$ that intersects $U_a$ at a time less than $a$, and such that $U_{b} \subset U_a$ if $a<b$.  
Denote the Liouville flow on $\widehat{M}$ by $\varphi_t$ and let $g: \partial M \setminus B \rightarrow (0,1]$ given by $g(x) = t$ where $t$ is the unique number such that $\varphi_{\log t}(x) \in \partial W$. Note that $g$ tends to $0$ as $x$ tends to $B$. Define $N_{a} := \{x \in \partial M  \, |\, \varphi_{\log g(x)}(x) \in \partial W \setminus U_a \}$. 
Choose a family of smooth functions $f_a: \partial M \rightarrow (0,1]$, $a \in (0,\infty)$, with the property
\begin{align*}
f_a|_{N_a} = g, \text{ and} \text{ for all }x \in \partial M, \, f_a(x) \text{ is monotone decreasing in }a. 
\end{align*}

Note that $W \subset M_{f_b} \subset M_{f_a}$, for $b>a$ and $\partial W \setminus U_a \subset \partial M_{f_a}$.
Define $\sigma(a) = \frac{1}{\min_{\partial M} f_a}$.
Define $Q_a = \HW^a(M_{f_a}, L_0 \to L_1)$, where by abuse of notation we write $L_i$ instead of $L_i \cap M_{f_a}$, $i= 0,1$. 
For $a < b$ define $\pi_{a\to b} : Q_a \to Q_b$ as the composition of the Viterbo map $\HW^a(M_{f_a}, L_0 \to L_1) \to \HW^a(M_{f_b}, L_0 \to L_1)$ and the persistence map $\HW^a(M_{f_b}, L_0 \to L_1) \to \HW^b(M_{f_b}, L_0 \to L_1)$.
By the commutativity of the Viterbo map with persistence maps and by functoriality of the Viterbo map
 it follows that $\pi_{a\to c} = \pi_{b\to c} \circ \pi_{a \to b}$, for $a<b<c$, and hence $(Q,\pi)$ is a filtered directed system.
 Furthermore, $\phi:U \rightarrow Q$, with $\phi_a:\HW^a(M) \rightarrow \HW^a(M_{f_a})$ is a morphism of f.d.s. 
 We define $\psi: Q \rightarrow U({\sigma})$ by the Viterbo transfer $\HW^a(M_{f_a}) \rightarrow \HW^a(M_{\min f_a}) = \HW^{\sigma(a)\cdot a}(M)$. It is clear that $(\phi, \psi)$ is a weak interleaving of $U$ and $Q$.

It remains to show that $Q$ and $V$ are isomorphic. 
Let $a > 0$. Assume that $L_0$ and $L_1$ are conical in the complement of $W_{\frac{1}{2}}$. 
Let $H_{\mu}$ be an admissible Hamiltonian with slope $\mu$ with respect to $W_{\frac{1}{2}}$.  Consider a Hamiltonian $K_{\mu}$ such that
\begin{align}
K_{\mu}(x) = &H_{\mu}(x), \text{ if } x \in W_{\frac{3}{4}}. \\
K_{\mu}(x) = &H_{\mu}(x), \text{ if } x = (r,y) \in (0, +\infty) \times \partial W\setminus U_a. \label{K2}
\end{align}
It follows that $K_{\mu}(x) = 2\mu r + b$, for some $b\in \R$, where $x$ is written in the coordinates $(r,y) \in (1, \infty) \times \partial W \setminus U_a$. 
Hence we can assume additionally that 
\begin{align}
K_{\mu}(x) = &2\mu r + b,  \text{ where } x = (r,y) \in (1, \infty) \times \partial M_{f_a}. \label{K3} 
\end{align}
By definition of $U_a$, $\mathcal{A}_{H_{\mu}}$ and $\mathcal{A}_{K_{\mu}}$ have the same critical points, and so it follows from \cite[Lemma 7.2]{AbouzaidSeidel2010} that we actually have  $\HW(K_{\mu}) = \HW(H_{\mu})$. On the other hand $K_{\mu} - {\frac{1}{2}}\mu$ is admissible with respect to $M_{f_a}$ with slope $2\mu$. \color{black} 
 One  concludes, reasoning as in Lemma \ref{M_epsilon}, that 
$Q_a = \HW^{a}(M_{f_a}, L_0 \to L_1)  \cong \HW^{\frac{1}{2}a}(W_{\frac{1}{2}},L^{'}_0 \to L^{'}_1)$ which is, by Lemma \ref{M_epsilon}, isomorphic to  $\HW^{a}(W,L^{'}_0 \to L^{'}_1) = V_a$.
That this identification respects the persistence morphisms of $Q$ and $V$ is again deduced from the functoriality of the Viterbo maps and the fact that the Viterbo maps are themselves morphisms of filtered directed systems.
Denote the isomorphism from $Q$ to $V$ by $\tau$. We have obtained a weak interleaving $(\tau \circ \phi, \psi \circ \tau^{-1})$. Moreover $\tau \circ \phi = {j}_{!}$ by construction. \color{black} 
\qed

\subsection{Plumbing}

Let $Q_1$  and $Q_2$ be closed orientable  $n$-dimensional manifolds. We let $D^*Q_i$ be the unit cotangent bundle of $Q_i$. We choose balls $B_i \subset Q_i$ in each $Q_i$.
The plumbing $N$ of $D^*Q_1$ and $D^*Q_2$ is obtained by identifying $D^*B_1$ and $D^*B_2$ via a symplectomorphism that swaps the momentum and position coordinates of these manifolds; see \cite{AbouzaidSmith2012,GeigesBook} for the details. 
There are obvious embeddings of $D^*(Q_1 \setminus B_1)$ and  $D^*(Q_2 \setminus B_2)$ into $N$.
 It is shown in \cite[Section 4]{AbouzaidSmith2012} that $N$ admits a Liouville structure which coincides with those of  $D^*(Q_i \setminus B_i)$ on the image of these embeddings. This implies that for points $q_1 \in Q_1\setminus B_1$ the cotangent disc fibre $L_{q_1}$ over $q_1$ survives as a conical exact Lagrangian in the Liouville domain $N$.

This construction can be generalised in the following way. Let $\{Q_i \mid 1 \leq i \leq k\}$ be a finite collection of  orientable  $n$-dimensional manifolds. Let $\mathsf{T}$ be a tree with $k$ vertices and use a bijection to associate to each vertex a manifold $Q_i$. For each edge $\eta$ leaving the ``vertex'' $Q_i$ we choose an embedded open ball $B_i(\eta)$ in $Q_i$. We assume that these balls are chosen to be disjoint and do not cover $Q_i$.  For all $i\neq j$ and every edge $\eta$ connecting $Q_i$ and $Q_j$ (there can be at most one such edge as $\mathsf{T}$ is a tree) we identify $D^*(B_i(\eta))$ and $D^*(B_j(\eta))$ by the recipe explained in the previous paragraph. The resulting manifold $N$ can be given a Liouville structure as explained in  \cite[Section 4]{AbouzaidSmith2012} and \cite{GeigesBook}. Let $\dot{Q}_1$ be the complement of the ``edge balls'' in $Q_1$, and $q_1 \in \dot{Q}_1$.
  In \cite[Section 4]{AbouzaidSmith2012} the following result is proved.
 \begin{thm} \cite{AbouzaidSmith2012} \label{thm:plumbing}
 There exists an injective algebra homomorphism from the group algebra $\mathbb{Z}_2[\pi_1(Q_1)]$ to $\HW(N,L_{q_1})$.
 \end{thm}
In fact the injective algebra homomorphism obtained in \cite{AbouzaidSmith2012} is for the respective homologies with $\mathbb{Z}$ coefficients, and applying the Universal Coefficient Theorem one obtains the homomorphism mentioned above. Thus if $\pi_1(Q_1)$ grows exponentially then $\HW(N,L_{q_1})$ has exponential algebraic growth; see Section \ref{preliminaries}.

\textit{Proof of Proposition \ref{prop:operations}}: 
Part A) follows from Proposition \ref{Viterbo_iso} and Part B) follows from Theorem \ref{thm:plumbing}. \qed

\section{Construction of  contact structures with positive entropy }\label{sec:constructions}

In this  section we prove Theorem \ref{spheres} and Theorem \ref{maincorollary}. 

\subsection{Preliminaries} \label{preliminaries}

Let $Q$ be a closed connected smooth manifold and $g$ a Riemannian metric on $Q$.  Let $(D_{g}^*Q,\lambda_{geo})  \subset  (T^*Q,\lambda_{geo})$ be the unit disk bundle with respect to the Riemannian metric $g$ where  $\lambda_{geo}$ is  the canonical Liouville form on $T^*Q$. 
By  Theorem  \ref{ASring} of Abbondandolo and Schwarz the map  
\begin{equation}
\Psi_{AS,q_1} : \Ho_*(\Omega_{q_1}(Q))\to \HW(D_g^*Q , L_{q_1})
\end{equation} \color{black}
is an algebra isomorphism. 
It is well-known that there is an algebra isomorphism
\begin{equation}
\Phi: \mathbb{Z}_2[\pi_1(Q,q_1)] \to \Ho_0(\Omega_{q_1}(Q)).
\end{equation}
Composing these two maps  we obtain an injective algebra homomorphism  \color{black}
\begin{equation}
\widetilde{\Phi}: \mathbb{Z}_2[\pi_1(Q,q_1)]  \to \HW(D_g^*Q , L_{q_1}).
\end{equation}

For a finitely generated group $G$ and a finite set $\upsigma$ of generators of $G$, let $\widehat{\Gamma}_{\upsigma}(G)$ be the usual exponential growth of the group $G$ with respect to the set $\upsigma$; see \cite[Section VI.C]{delaHarpe}. 
To a finite set $\upsigma$ of generators of  $\pi_1(Q,q_1)$, we associate the finite set $S \subset \mathbb{Z}_2[\pi_1(Q,q_1)]$ that  is formed by the elements of  $\upsigma$ and its inverses.
It is immediate to see that 
\begin{equation}
 \widehat{\Gamma}_{\upsigma}(\pi_1(Q,q_1)) = {\Gamma_S^{\alg}}(\mathbb{Z}_2[\pi_1(Q,q_1)]).
  \end{equation}
 Using that $\widetilde{\Phi}$ is injective we obtain
 \begin{equation}
 \widehat{\Gamma}_{\upsigma}(\pi_1(Q,q_1)) = {\Gamma_S^{\alg}}(\mathbb{Z}_2[\pi_1(Q,q_1)])  \leq  {\Gamma_{\widetilde{\Phi}(S)}^{\alg}}( \HW(D_g^*Q , L_{q_1})),
\end{equation}
We have shown the following 
\begin{lem} \label{fundalgrowth}
It $\pi_1(Q,q_1)$ has exponential growth then there exists a finite set $S \subset \HW(D_g^*Q , L_{q_1})$ such that  $  {\Gamma_S^{\alg}}(D_g^*Q , L_{q_1}) >0$.
\end{lem}

 \subsection{ Proof of statement (A) of Theorem \ref{spheres}  and statement $\clubsuit$ of Theorem \ref{maincorollary} } \
 
 \textit{Proof of statement (A) of Theorem \ref{spheres} }
 
 \color{black}
 Let $G$ be a finitely presented group such that \color{black}
 \begin{itemize}
 \item $\Ho_1(G) = \Ho_2(G) =0$,
 \item $G$ has exponential growth,
 \item $G$ admits a presentation on which the number of relations does not exceed the number of generators.
 \end{itemize}
 Then, it follows from \cite{Kervaire}, that for every $n\geq 4$ there exists a manifold $Q^n$ which is an integral homology sphere and which satisfies $\pi_1(Q^n) = G$. We denote by $\varrho(G)$ the minimal number of generators of $G$.
 
 We denote by $D^*Q^n$ the unit disk bundle of $Q^n$, with respect to a Riemannian metric $g$ in $Q^n$, endowed with the canonical symplectic and Liouville forms. \color{black} We choose a point $q\in Q^n$ and $g$ generically so that $q$ is not conjugate to itself. Let $S^*Q^n= \partial D^*Q^n $ be the unit cotangent bundle of $Q^n$. In order to prove our result we consider two distinct cases.  \color{black}

 \textbf{Case 1: $n$ is odd and $\geq 5$.} \\
 In this case the Euler characteristic of $Q^n$ vanishes.  Because $G$ grows exponentially, we know that $\HW_0(D^*Q^n,L_q)$ has exponential algebraic growth. Let $N^1$ be the plumbing  of $D^*Q^n$ and $D^* S^n$ performed far from $L_q$. \color{black}
By Proposition \ref{prop:operations}, $\HW_0(N^1,L_q)$ has exponential algebraic growth. \color{black}
  
 It is a result of Milnor that the boundary of the plumbing of the unit disk bundles of two odd-dimensional homology spheres of dimension $\geq 3$ is a homology sphere; see \cite[Chapter VI - Section 18]{Bredon}. Applying this to the pair
   $D^*Q^n$ and $D^* S^n$  we conclude that $\partial N^1 $ is a homology sphere. 
   Since $N^1$ retracts to the one point union of $Q$ and $S^n$ we know that the homology of $N^1$ \color{black} is zero in  every degree different from $0$ and $n$, where we have $\Ho_0(N^1) = \mathbb{Z}$ and $\Ho_n(N^1) = \mathbb{Z} \oplus \mathbb{Z}$. \color{black}

  \textbf{Case 2: $n$ is even and $\geq 4$.} \\
 In this case the Euler characteristic of $Q^n$ is $2$. We consider the plumbing associated to the E8 tree; see \cite[Chapter VI - Section 18]{Bredon}. To each vertex of the E8 tree we associate a disk bundle in the following way:   
 \begin{itemize}
 \item to the leftmost vertex we associate $D^*Q^n$,
 \item to every other vertex we associate $D^* S^n$.
 \end{itemize}
 We let $N^1$ be the plumbing associated to the E8 tree determined by this choice of disk bundles at each vertex, and assume that the plumbing is done away from a cotangent fibre $L_q\subset D^*Q^n$ .  It was shown by Milnor (see \cite[Chapter VI - Section 18]{Bredon} ) that $\partial N^1$ is a homology sphere. Since $N^1$ retracts to the wedge sum of $Q$ and seven copies of $S^n$ determined by the E8 tree, we know that the homology of $N^1$ is zero \color{black} in every degree  different from $0$ and $n$, where we have $\Ho_0(N^1) = \mathbb{Z}$ and $\Ho_n(N^1) = \oplus_{i=1}^8 \mathbb{Z} $. \color{black} By Proposition \ref{prop:operations}, $\HW_0(N^1,L_q)$ has exponential algebraic growth. \color{black}
\color{black}
 
\textbf{We now treat both cases simultaneously.} \color{black}
By attaching $2$-handles to $N^1$ away from $L_q$ we can obtain a simply connected Liouville domain $N^2$ such that $\HW(N^2,L_q)$ has exponential algebraic growth. We choose the framing of these handle attachments so that the first Chern class of $N^2$ vanishes.

 The effect of the handle attachment on the homology of the boundary can  be read from the surgery formula in \cite[Section X.1]{Kosinsky}.  \color{black}
 One concludes that the homology of $\partial N^2$ coincides with that of $\partial N_1$ except in degree $2$, and $\Ho_2(\partial N^2)$ is the  direct sum of $\varrho(G)$ copies of $\mathbb{Z}$.  \color{black}

By Hurewicz' Theorem there is a basis of $\Ho_2(\partial N^2)$ which is composed of embedded $S^2$. Since  the first Chern class of $N^2$ vanishes, it follows from \cite[Lemma 2.19]{McLean2011} that these $S^2$ can be made isotropic and disjoint from $L_q$ by an isotopy and that their symplectic normal bundle is trivial. We can thus perform the Weinstein handle attachment over these spheres. The resulting Liouville domain $N^3$ still contains the Lagrangian $L_q$ and it follows from Proposition \ref{prop:operations} that $\HW(N^3,L_q)$ has exponential algebraic growth. By the surgery formula in \cite[Section X.1]{Kosinsky}, the effect of these handle attachments on the homology of the boundary implies that $\Ho_2(\partial N^3)=0$ and that the homologies of $\partial N^3$ and $\partial N^2$ coincide in all other degrees.  \color{black}
 Therefore, $\partial N^3$ is a simply connected homology sphere. It follows from Whitehead's Theorem for homology \cite[Corollary 4.33]{Hatcher} that $\partial N^3$ also has the  homotopy groups of a sphere.   \color{black}
 Since the dimension of $\partial N^3$ is $>5$ the h-cobordism theorem tells us that $\partial N^3$ is homeomorphic to a sphere.  \color{black}
 Since the smooth spheres under connected sum form a finite group, we can take the connected sum of finitely many copies of $\partial N^3$ to get the sphere $\partial N^4$ with the standard smooth structure such that $\HW(N^4,L_q)$ has exponential algebraic growth.  \color{black}
   This proves statement (A) of Theorem \ref{spheres}.
\qed

\

\textit{Proofs of statement $\clubsuit$ of Theorem \ref{maincorollary} }
 Let $V$ be a $(2n-1)$-dimensional manifold where $n\geq 4$, and assume that there exists an exactly fillable contact structure $\xi$ on $V$. Denote by $M_V$ a Liouville domain whose boundary is $(V,\xi)$.
 \color{black}
  Let $N^4$ be the Liouville domain constructed in the proof of statement (A) of Theorem \ref{spheres}. 
By Proposition \ref{prop:operations}, the Liouville domain $N^5=N^4 \# M_V$ has an asymptotically conical exact Lagrangian $L$ such that $\HW(N^5,L)$ has exponential algebraic growth.   \color{black}
The statement then follows from Theorem \ref{theorementropy}. \qed 

\subsection{Proof of statement (B) of Theorem \ref{spheres} and statement $\diamondsuit$ of Theorem \ref{maincorollary}} 

\

\textit{Proof of statement (B) of Theorem \ref{spheres}:}   \\
We will consider a carefully chosen $3$-manifold $Q$.
Consider the Brieskorn manifolds of dimension $3$, $M(p,q,r) = \{(z_1,z_2,z_3) \in \C^3 \, |\, {z_1}^p + {z_2}^q + {z_3}^r =0\} \cap S^5$. $M(p,q,r)$ is a $\Z$-homology sphere if $p,q,r$ are relatively prime (see for example \cite{Saveliev2002}).
It was shown by Milnor \cite{Milnor1975} that its fundamental group $\pi_1(M(p,q,r))$ is the commutator subgroup of the group 
$G = G(p,q,r) = \langle \gamma_1, \gamma_2, \gamma_3 \, | \, \gamma_1^p = \gamma_2^q = \gamma_3^r = \gamma_1\gamma_2\gamma_3 \rangle$ , see also \cite{Seade2006}.
The groups $\Sigma = G/Z(G)$ are the triangle groups, where $Z(G)$ is the center of $G$.   
Consider the case $p=2$, $q=3$, $r=7$.
 A short computation shows that $G(2,3,7) = [G(2,3,7), G(2,3,7)] = \pi_1(M(2,3,7))$.  One has $\widehat{\Gamma}(G(2,3,7)) \geq \widehat{\Gamma}(\Sigma(2,3,7))$, and the  exponential growth of $\Sigma(2,3,7)$ is $\log(x)$, where $x \approx 1.17628$ is equal to Lehmer's Salem number  (see \cite{Hironaka2003} or \cite{Breuillard2014}). 
We take $Q=M(2,3,7)$.
The integral homology of $D_g^*Q$ is the same as that of $Q$, which is $\mathbb{Z}$ in degrees $0$ and $3$ and vanishes in all other degrees. Moreover it is clear that $\pi_1(S^*Q)= \pi_1(Q\times S^2)= \pi_1(Q)$ is generated by the elements $\gamma_1$ and $\gamma_2$.
 \color{black}

Let $N^1$ be the Liouville domain obtained by plumbing $D_g^*Q$ with the unit disk bundle $D^*S^3$ of $S^3$. 
 \color{black} We assume that the plumbing is performed away from the cotangent fibre $L_q$ over a point $q \in Q$. 
Therefore  $L_q$ survives as a conical exact Lagrangian in $N^1$. 
By Proposition \ref{prop:operations} we know that  $\HW_*(N^1,L_{q})$ has exponential algebraic growth.

Since $N^1$ is the plumbing of $D_g^*Q$ and $D^*S^3$, and $Q$ and $S^3$ are both homology spheres we obtain that $\partial N^1$ is a homology sphere; see  \cite[Chapter VI - Section 18]{Bredon}.
 \color{black}
 Combining this with the fact that  $N^1$ retracts to the one point union of $S^3$ and $Q$  we conclude that   \color{black}
\begin{itemize}
\item $\Ho_0(N_1) = \mathbb{Z}$, $\Ho_3(N_1) = \mathbb{Z} \oplus \mathbb{Z}$ and $\Ho_i(N_1) =0$ for $i\neq 0,3$,
\item $\Ho_0( \partial N_1) = \mathbb{Z}$,  $\Ho_5( \partial N_1) = \mathbb{Z}$, and  $\Ho_i( \partial N_1) =0$ for  $i\neq 0,5$.
\end{itemize}

Let now $\{\overline{\sigma}_1,\overline{\sigma}_{2}, \overline{\sigma}_3\}$ be generators of $\pi_1(\partial N_1)= \pi_1(Q)$ corresponding to $\gamma_1$, $\gamma_2$ and $\gamma_3$ respectively.  By the h-principle for subcritical isotropic submanifolds of contact manifolds \cite{GeigesBook} we can isotope the curve $\overline{\sigma}_3$ to a curve ${\sigma}_3$ which is isotropic in $(S_g^*Q,\xi_{geo})$. We can also assume that ${\sigma}_3$  does not intersect $\Lambda_{q} := \partial L_{q}$. Since ${\sigma}_3$ is isotropic and has trivial normal bundle we can apply the Weinstein handle attachment \cite{Weinstein1990} and attach a  2-handle to $ N^1$ over  ${\sigma}_3$, obtaining a new Liouville domain $N^2$. From the presentation of $\pi_1(Q)$ that we used, it is clear that $\partial N^2$ is simply connected, and so is $N^2$ by \cite[Lemma 2.9]{McLean2011}. We choose the framing of the handle attachment so that $\partial N^2$ is spin. 
Using the Mayer-Vietoris sequence we obtain that
$\Ho_0(\partial N^2) = \mathbb{Z}$,  $\Ho_2(\partial N_2) = \mathbb{Z} $, and $\Ho_1(\partial N_2) =0$.  \color{black}
By Smale's classification of spin simply-connected five manifolds \cite{Smale} it follows that $\partial N_2$ is  $S^3\times S^2$.  \color{black}

Since $N^2$ is obtained from $N^1$ via a subcritical handle attachment and the Lagrangian $L_q$ is far from the attaching locus of this handles, we know that $L_q$ survives as a conical exact Lagrangian in $N^2$. Moreover Proposition \ref{prop:operations} implies that  $\HW_*(N^2,L_{q})$  has exponential algebraic growth, and it follows from Theorem \ref{theorementropy} that the contact manifold $\partial N_2$ has positive entropy. \qed

\textit{Proof of statement $\diamondsuit$ of Theorem \ref{maincorollary}:} \color{black}
The  statement is proved by a connected sum argument identical to the one in the proof of statement $\clubsuit$. \qed   \color{black}

\begin{rem} To guarantee the vanishing of the second Stiefel-Whitney class of $\partial N^2$ one must only  guarantee the vanishing of the first Chern class of $N^2$. As observed in the proof of \cite[Lemma 2.10]{McLean2011}, one can choose the framing when performing the attachments of the 2-handles so as to guarantee the vanishing of the first Chern class of $N^2$.
 \end{rem}

\appendix

\section*{Appendix}

\subsection*{Construction of exact Lagrangian cobordisms.}

\
Before proving the Lemma we recall that the symplectization of the contact form $\alpha$ on $(\Sigma,\xi)$ can also be given by $(\mathbb{R} \times \Sigma,e^s\alpha)$, where $s$ denotes the $\R$-coordinate. It is straightforward to see that the diffeomorphism $\mathrm{F}: ((0,+\infty) \times \Sigma,r\alpha ) \to (\mathbb{R} \times \Sigma,e^s\alpha)$ given by $\mathrm{F}(r,p) = (\log r, p)$ is an exact symplectmorphism.

It follows that an exact Lagrangian submanifold $\mathcal{L}^-$ is conical in $ ((0,+\infty) \times \Sigma,r\alpha )$ outside $[1-2\epsilon,1-\epsilon] \times\Sigma$ if, and only if, $\widetilde{\mathcal{L}}^-:=\mathrm{F}(\mathcal{L}^-)$ is conical in $(\mathbb{R} \times \Sigma,e^s\alpha)$ outside $[\log (1-2\epsilon),\log(1-\epsilon)] \times \Sigma$. Analogously, an exact Lagrangian submanifold $\mathcal{L}^+$ is conical in $ ((0,+\infty) \times \Sigma,r\alpha )$ outside $[1+\epsilon,1+2\epsilon] \times\Sigma$ if, and only if, $\widetilde{\mathcal{L}}^+:=\mathrm{F}(\mathcal{L}^+)$ is conical in $(\mathbb{R} \times \Sigma,e^s\alpha)$ outside $[\log (1+\epsilon),\log(1+ 2\epsilon)] \times \Sigma$. 

\textit{Proof of Lemma  \ref{lemmaBaptiste}:} 

We use the technique presented in \cite[Lemma 6.3]{EHK}.

\textbf{Step 1:}
We first apply the Legendrian neighbourhood Theorem \cite[Proposition 43.18]{MichorKriegl} to find a neighbourhood $\widetilde{U}(\Lambda_0)$ of the Legendrian   of  $\Lambda_0$ such that there exists a strict contactmorphism $\Upsilon: (\widetilde{U}(\Lambda_0), \alpha) \to (V(\widetilde{\Lambda}_0) \subset \mathbb{R}^{2n-1}, dz + \sum^{n-1}_{i=1}x_i dy_i) $ that satisfies $\Upsilon(\Lambda_0)= \widetilde{\Lambda}_0$,
where $\widetilde{\Lambda}_0$ is the standard Legendrian unknot in $\mathbb{R}^{2n-1}$ (see \cite[Example 3.1]{EES}) and $V(\widetilde{\Lambda}_0)$ is a tubular neighbourhood of $\widetilde{\Lambda}_0$ in $\mathbb{R}^{2n-1}$.
 Given this identification, it suffices to establish the lemma for the case of $\widetilde{\Lambda}_0$ since it will follow if we can establish it for  $U(\Lambda_0) \subset \widetilde{U}(\Lambda_0) $.

\textbf{Step 2:} Clearly, it suffices to establish the lemma for all $0<\epsilon<\frac{1}{e^{1000}}$.
We thus fix $0<\epsilon<\frac{1}{e^{1000}}$ and a tubular neighbourhood $U(\widetilde{\Lambda}_0)$ of $\widetilde{\Lambda}_0$.
To establish the lemma for $\widetilde{\Lambda}_0$ we assume that $\Lambda_1$ is a Legendrian sphere in $(\mathbb{R}^{2n-1},\alpha_{\can} = dz + \sum^{n-1}_{i=1}x_i dy_i)$ that is $\mu$-close to $\widetilde{\Lambda}_0$ in the $C^3$-sense, with $\mu>0$ so small that
\begin{itemize}
\item $\Lambda_1\subset U(\widetilde{\Lambda}_0)$, 
\item there exists a Legendrian isotopy $\theta: [-1,1] \times S^{n-1} \to \mathbb{R}^{2n-1}$ which is $\mu$-small in the $C^3$-topology and satisfies $\theta(\{-1\}\times S^{n-1}) = \widetilde{\Lambda}_0$ and $\theta(\{1\}\times S^{n-1}) = {\Lambda}_1$.
\end{itemize}
Moreover if $\mu>0$ is chosen sufficiently small we can also assume that 
\begin{itemize}
\item the isotopy $\theta$ is constant in the first coordinate for $t \notin [\log (1-2\epsilon),\log(1-\epsilon)] $.
\end{itemize}
Extend $\theta$ to $\R\times S^{n-1}$ by $\theta(t,p) =\theta(-1,p)$ for $t\leq -1$ and  $\theta(t,p) =\theta(1,p)$ for $t\geq 1$.

We write $\theta(t,p)= (x(t,p),y(t,p),z(t,p))$ for coordinates $(x,y,z)\in \mathbb{R}^{n-1} \times \mathbb{R}^{n-1}\times \mathbb{R}$, set $F(t,p):= \alpha_{\can}(\partial_t \theta(t,p)) $, and  define the cylinder  $\Theta: \mathbb{R}\times S^{n-1}\to (\mathbb{R}\times \mathbb{R}^{2n-1},e^s\alpha_{\can})$ in the symplectization $(\mathbb{R}\times \mathbb{R}^{2n-1},e^s\alpha_{\can})$ of $\alpha_{\can}$ by
\begin{equation}
\Theta(t,p) = (t,x(t,p),y(t,p),z(t,p)+F(t,p)).
\end{equation}
It is clear that if $\mu>0$ is chosen sufficiently small then $\Theta$ will be an embedding, since it will be a small compact perturbation of the embedding $(t,x(t,p),y(t,p),z(t,p))$. Let $\Pi :\mathbb{R}\times \mathbb{R}^{2n-1} \to  \mathbb{R}^{2n-1}$  be the projection of the symplectization to the contact manifold.

A direct computation shows that 
\begin{equation}
 \Theta^*(e^s\alpha_{\can})= d(e^tF(t,p)).
\end{equation}

\textbf{Step 3.}
Step 2 implies that the cylinder $\widetilde{\mathcal{L}}^-= \Theta( \mathbb{R}\times S^{n-1})$ is an admissible exact Lagrangian submanifold 
of $(\mathbb{R}\times \mathbb{R}^{2n-1},e^s\alpha_{\can})$. By taking $\mu$ even smaller we can guarantee that 
\begin{itemize}
\item $\widetilde{\mathcal{L}}^-$ is conical over  $\Lambda_1$ in $[\log(1-\epsilon), +\infty) \times \mathbb{R}^{2n-1}$ and is conical over  $\widetilde{\Lambda}_0$ in $(0,\log(1-2\epsilon)] \times \mathbb{R}^{2n-1}$,
\item and the projection of $\widetilde{\mathcal{L}}^-$ to $ \mathbb{R}^{2n-1}$ is contained in  $U(\widetilde{\Lambda}_0)$.
\end{itemize}

In order to construct $\widetilde{\mathcal{L}}^+$ we use the inverse isotopy $\theta^+(t,p) = \theta(-t,p)$, and apply the construction above. We have established statements a), b) and c) of Lemma \ref{lemmaBaptiste}.

To conclude d), notice that  $f^-:= e^tF(t,p)$ has support in $[\log(1-2\epsilon),\log(1 -\epsilon)] \times S^{n-1}$ and if $\mu>0$ is chosen small enough then $f^-:= e^tF(t,p)$ will satisfy  $| f^- |_{C^0} \leq \epsilon$ since the $C^0$-norm of $F(t,p)$ will be small. Applying the same argument to $\widetilde{\mathcal{L}}^+$ implies d).

\textbf{Step 4.}
Statement e) is obtained by performing this construction for a smooth $1$-parameter family of Legendrian isotopies which starts with the isotopy $\theta$ and ends at the stationary isotopy from $\widetilde{\Lambda}_0$ to itself. As the construction above depends $C^\infty$-smoothly on the parameter we obtain a smooth $1$-parameter family of exact Lagrangian embeddings $\Psi: \mathbb{R}\times S^{n-1}\to (\mathbb{R}\times \mathbb{R}^{2n-1},e^s\alpha_{\can})$ which starts at $\mathcal{L}$ and ends at $\mathbb{R}\times \widetilde{\Lambda}_0$, and which is constant in $(\mathbb{R}\setminus [ \log (1-2\epsilon),\log(1+2\epsilon)]) \times S^{n-1}$. This is an exact Lagrangian isotopy from $\widetilde{\mathcal{L}} := \mathrm{F}(\mathcal{L})$ to  $\mathbb{R}\times \widetilde{\Lambda}_0$ which is constant in $(\mathbb{R}\setminus [\log(1-2\epsilon),\log(1+2\epsilon)]) \times \widetilde{\Lambda}_0$. Statement e) then follows from \cite[Lemma 5.6]{Mclean2015}.
\qed

\bibliographystyle{amsalpha}
\bibliography{References}

\end{document}

%% file: figure1.tex
\begin{figure}\label{fig:1}
 \begin{tikzpicture}
\draw[-] (0,-0.1) -- (0.8,-0.1);
\draw[-] (0.8,-0.1) to [out=0,in=225] (1,0);
\draw[-] (1.0,0) -- (2.8,1.8);
\node[above left] at (2,1) { $\mu$};
\draw[-] (2.8,1.8) to [out=45,in=180] (3,1.9);
\draw[-] (3,1.9) -- (5.8,1.9);
\draw[-] (5.8,1.9) to [out=0,in=225] (6,2);
\draw[-] (6,2) -- (7,3);
\node[above left] at (6.5,2.5) { $\mu$};
\draw[-|] (0,0) -- (1,0)node[below] { $1$}; 
\draw[-|] (1,0) -- (3,0)node[below] { $R$}; 
\draw[-] (3,0) -- (4,0) node[anchor=north] {$r=r_W$};
\draw[dotted] (4,0) -- (5.5,0);
\draw[-] (5.5,0) -- (6,0);
\draw[|->] (6,0)node[below] { $R$}-- (8,0) node[anchor=north] { $r=r_M$};

\draw[-|] (0,-1) -- (0,-0.1)node[left] { $-\epsilon$};
\draw[-|] (0,-0.1) -- (0,1.9)node[left] { $\mu(R-1) -\epsilon$};
\draw[->] (0,1.9) -- (0,4) node[anchor=east] { $H^{step}_{\mu}$};

\draw[dotted] (0,1.9) -- (3,1.9);
\draw[dotted] (3,0) --  (3,1.9);
\draw[dotted] (6,0) -- (6,2);

\node[below] at (0.2,-0.5) { $\mathfrak{A}^{*}$};
\node[below] at (1,-0.5)   { $\mathfrak{A}^{**}$};
\node[below] at (3,-0.5)   { $\mathfrak{B}^{**}$};
\node[below] at (4.5,-0.5) { $\mathfrak{B}^{*}$};
\node[below] at (6,-0.5)   { $\mathfrak{B}^{**}$};
\end{tikzpicture}
\caption{}
\end{figure}

%% file: figure3.tex
 \begin{figure}\label{fig:3}
 \begin{tikzpicture}
\draw[-|] (1,1) -- (4,1);
\draw[-]  (4,1) -- (6,1);
\draw[-] (6,1) to [out=0, in=180] (7,1.5);
\draw[-|] (7,1.5) -- (8,1.5);
\draw[-] (8, 1.5) -- (8.5,1.5); 
\draw[-] (8.5,1.5) to [out=0, in=180] (9.5,2);
\draw[-|] (9.5,2) -- (10,2);
\draw[-] (10,2) -- (10.5,2);
\draw[-] (10.5,2) to [out=0, in=180] (11.5,1.5);
\draw[-|] (11.5,1.5) -- (12,1.5);
\draw[-] (12,1.5) -- (13,1.5);
\draw[-] (13,1.5) to [out=0, in=180] (14,1);
\draw[-|] (14,1) -- (16,1);

\node[above] at (4,1) { $L$};
\node[above] at (8,1.5)   { $L_1$};
\node[above] at (10,2)   { $L_2$};
\node[above] at (12,1.5) { $\widetilde{L_1}$};
\node[above] at (16,1)   { $\widetilde{L}$};

\draw[<-] (3,1.25) -- (3.75,1.25);
\draw[<-] (6,1.8)  -- (7.7,1.8);
\draw[<-] (8,2.3)  -- (9.7,2.3);
\draw[<-] (11.3,1.8)  -- (11.7,1.8);
\draw[<-] (14,1.3)  -- (15.75,1.3);

\draw[-|] (1,0) -- (4,0) node[below] { $1-\epsilon$}; 
\draw[-|] (4,0) -- (8,0) node[below] { $1-\frac{\epsilon}{5}$}; 
\draw[-|] (8,0) -- (10.0,0) node[anchor=north] {$1$};
\draw[-|] (10.0,0) -- (12,0) node[below] {$1 + \frac{\epsilon}{5}$};
\draw[-|] (12,0) -- (16,0) node[below] {$1 + \epsilon$}; 
\draw[|->] (16,0) -- (17,0) node[anchor=north] { $r$};


\draw[dotted] (4,0) -- (4,1);
\draw[dotted] (8,0) --  (8,1.5);
\draw[dotted] (10,0) -- (10,2);
\draw[dotted] (12,0) -- (12,1.5);
\draw[dotted] (16,0) -- (16,1);

\end{tikzpicture}
\caption{}\end{figure}

%% file: figure2.tex
\begin{figure}\label{fig:2}
 \begin{tikzpicture}

\draw[-] (-5,5) -- (0,5); 
\draw[-] (0,5) -- (5,5) node[below] { $\partial M$};

\draw[-|] (-6,1) -- (-6,5) node[left] { $1$};
\draw[->] (-6,5) -- (-6,6) node[anchor=east] { $r=r_M$};  

\draw[fill] (0,5) circle [radius=0.05] node[above right] { $B$};
\draw[dotted] (-1,3) -- (-1,5);
\draw[dotted] (1,3) -- (1,5);
\draw[decoration={brace,mirror,raise=5pt},decorate]
  (-1,5) -- node[below=6pt] {$\partial M \setminus N_a$} (1,5);

\begin{scope}
\clip (-5,1) rectangle (5,5);
\draw[dotted] (0,0) to [out=90,in=246] (0.3,2);
\draw[dotted] (0,0) to [out=90,in=294] (-0.3,2);
\end{scope}
\draw[dashed] (0.3,2) to [out=66,in=225] (1,3);
\draw[-] (1,3) to [out=45, in=180] (5,4) node[below] { $\partial M_{f_a}$}; 
\draw[dashed] (-0.3,2) to [out=114,in=315] (-1,3);
\draw[-] (-1,3) to [out=135,in=360] (-5,4) ;

\draw[-] (-1,3) to [out=315, in=225] (1,3); 
\draw[dashed] (-0.3,2) to [out=294, in=246] (0.3,2);
\node[right]  at (0.25,2) { $\partial M_{f_b}$};
\node[right] at (0,1.3) { $\partial W\setminus S$};


\end{tikzpicture}
\caption{}
\end{figure}